\documentclass[letterpaper,12pt,reqno]{amsart}
\RequirePackage[utf8]{inputenc}
\usepackage[portrait,margin=3cm]{geometry}
\usepackage{mathrsfs,soul}
\usepackage{hyperref}
\usepackage[foot]{amsaddr}
\usepackage{amssymb,amsthm,amsfonts,amsbsy,latexsym,dsfont}
\usepackage[textsize=small]{todonotes}       
\usepackage{graphicx}
\usepackage[numeric,initials,nobysame]{amsrefs}
\usepackage{upref,setspace}
\usepackage{xcolor,colortbl}

\usepackage{fourier}

\usepackage{xcolor}

\usepackage{amsbsy,latexsym,dsfont}
\usepackage{mathrsfs,soul}

\usepackage{enumerate}
\newenvironment{enumeratei}{\begin{enumerate}[\upshape (i)]}{\end{enumerate}}
\newenvironment{enumeratea}{\begin{enumerate}[\upshape (a)]}{\end{enumerate}}

\usepackage{paralist}

\newenvironment{inparaenuma}{\begin{inparaenum}[\upshape (a) ]}{\end{inparaenum}}

\newenvironment{inparaenumi}{\begin{inparaenum}[\upshape  (i) ]}{\end{inparaenum}}

\definecolor{refkey}{gray}{.75}
\definecolor{labelkey}{gray}{.75}

\newtheorem{thm}{Theorem}[section]
\newtheorem{lem}[thm]{Lemma}
\newtheorem{cor}[thm]{Corollary}
\newtheorem{prop}[thm]{Proposition}

\newtheorem{ass}[thm]{Condition}

\theoremstyle{remark}

\theoremstyle{definition}
\newtheorem{rem}{Remark}

\renewcommand{\leq}{\leqslant}
\renewcommand{\geq}{\geqslant}

\newcommand{\ind}{\mathds{1}}
\newcommand{\eps}{\varepsilon}

\newcommand{\set}[1]{\left\{#1\right\}}

\newcommand{\ie}{\emph{i.e.,}}

\newcommand{\equald}{\stackrel{\mathrm{d}}{=}}
\newcommand{\probc}{\stackrel{\mathrm{P}}{\longrightarrow}}

\newcommand{\weakc}{\stackrel{\mathrm{d}}{\longrightarrow}}
\newcommand{\convas}{\stackrel{\mathrm{a.s.}}{\longrightarrow}}

\def\qed{ \hfill $\blacksquare$}

\newcommand{\cA}{\mathcal{A}}\newcommand{\cB}{\mathcal{B}}\newcommand{\cC}{\mathcal{C}}
\newcommand{\cD}{\mathcal{D}}\newcommand{\cE}{\mathcal{E}}\newcommand{\cF}{\mathcal{F}}
\newcommand{\cG}{\mathcal{G}}\newcommand{\cI}{\mathcal{I}}

\newcommand{\cP}{\mathcal{P}}\newcommand{\cQ}{\mathcal{Q}}\newcommand{\cR}{\mathcal{R}}
\newcommand{\cS}{\mathcal{S}}\newcommand{\cT}{\mathcal{T}}\newcommand{\cU}{\mathcal{U}}
\newcommand{\cV}{\mathcal{V}}\newcommand{\cX}{\mathcal{X}}

\newcommand{\vzero}{\mathbf{0}}\newcommand{\vone}{\mathbf{1}}

\newcommand{\vD}{\mathbf{D}}\newcommand{\vE}{\mathbf{E}}

\newcommand{\vJ}{\mathbf{J}}
\newcommand{\vM}{\mathbf{M}}

\newcommand{\vT}{\mathbf{T}}

\newcommand{\vX}{\mathbf{X}}\newcommand{\vY}{\mathbf{Y}}
\newcommand{\va}{\mathbf{a}}\newcommand{\vb}{\mathbf{b}}
\newcommand{\vd}{\mathbf{d}}\newcommand{\ve}{\mathbf{e}}
\newcommand{\vg}{\mathbf{g}}

\newcommand{\vp}{\mathbf{p}}\newcommand{\vr}{\mathbf{r}}
\newcommand{\vs}{\mathbf{s}}\newcommand{\vt}{\mathbf{t}}\newcommand{\vu}{\mathbf{u}}
\newcommand{\vv}{\mathbf{v}}\newcommand{\vw}{\mathbf{w}}\newcommand{\vx}{\mathbf{x}}
\newcommand{\vy}{\mathbf{y}}

\newcommand{\mvw}{\boldsymbol{w}}

\newcommand{\mvalpha}{\boldsymbol{\alpha}}

\newcommand{\mvmu}{\boldsymbol{\mu}}

\newcommand{\mvpi}{\boldsymbol{\pi}}

\newcommand{\mvxi}{\boldsymbol{\xi}}\newcommand{\mvXi}{\boldsymbol{\Xi}}

\newcommand{\fA}{\mathfrak{A}}\newcommand{\fB}{\mathfrak{B}}\newcommand{\fC}{\mathfrak{C}}
\newcommand{\fE}{\mathfrak{E}}
\newcommand{\fI}{\mathfrak{I}}

\newcommand{\fN}{\mathfrak{N}}
\newcommand{\fR}{\mathfrak{R}}
\newcommand{\fS}{\mathfrak{S}}\newcommand{\fT}{\mathfrak{T}}
\newcommand{\fX}{\mathfrak{X}}
\newcommand{\fZ}{\mathfrak{Z}}
\newcommand{\fb}{\mathfrak{b}}

\newcommand{\fh}{\mathfrak{h}}

\newcommand{\fq}{\mathfrak{q}}\newcommand{\fs}{\mathfrak{s}}

\newcommand{\bD}{\mathbb{D}}\newcommand{\bE}{\mathbb{E}}
\newcommand{\bG}{\mathbb{G}}

\newcommand{\bN}{\mathbb{N}}
\newcommand{\bP}{\mathbb{P}}\newcommand{\bR}{\mathbb{R}}
\newcommand{\bT}{\mathbb{T}}

\newcommand{\bZ}{\mathbb{Z}}

\newcommand{\sC}{\mathfrak{C}}

\newcommand{\sP}{\mathscr{P}}

\newcommand{\sS}{\mathscr{S}}

\newcommand{\sW}{\mathscr{W}}

\newcommand{\current}{{\sf crnt}}
\newcommand{\CM}{\mathrm{CM}}

\newcommand{\vdd}{\mvalpha}

\DeclareMathOperator{\E}{\mathbb{E}}
\DeclareMathOperator{\pr}{\mathbb{P}}

\DeclareMathOperator{\var}{Var}
\DeclareMathOperator{\cov}{Cov}
\DeclareMathOperator{\dis}{dis}

\DeclareMathOperator{\GH}{GH}
\DeclareMathOperator{\GHP}{GHP}

\DeclareMathOperator{\ord}{ord}

\DeclareMathOperator{\con}{con}
\DeclareMathOperator{\cumu}{cumu}

\DeclareMathOperator{\ERRG}{ERRG}
\DeclareMathOperator{\IRG}{IRG}

\DeclareMathOperator{\bbp}{BP}
\DeclareMathOperator{\bbf}{BF}

\DeclareMathOperator{\dsc}{dsc}
\DeclareMathOperator{\scl}{scl}

\DeclareMathOperator{\mass}{mass}

\DeclareMathOperator{\spls}{spls}

\DeclareMathOperator{\free}{\sss [alv]}

\DeclareMathOperator{\project}{\sss [project]}

\DeclareMathOperator{\modi}{mod}

\DeclareMathOperator{\intra}{intra}
\DeclareMathOperator{\multi}{multi}
\DeclareMathOperator{\crit}{Crit}

\DeclareMathOperator{\Perc}{Perc}
\DeclareMathOperator{\FR}{AL}
\DeclareMathOperator{\pois}{Pois}

\newcommand{\sss}{\scriptscriptstyle}
\newcommand{\erdos}{Erd\H{o}s-R\'enyi }
\newcommand{\ldown}{l^2_{\downarrow}}

\newcommand{\td}{d^*}
\newcommand{\tR}{\tilde{R}}

\newcommand{\diam}{\mathrm{diam}}
\newcommand{\tH}{\tilde{H}}
\newcommand{\vCrit}{\mathbf{Crit}}

\newcommand{\girth}{\mathrm{girth}}

\newcommand{\bars}{\bar{s}}

\newcommand{\barg}{\bar{g}}
\newcommand{\barcd}{\bar{\cD}}

\usepackage[textsize=small]{todonotes}       

\definecolor{jrnl}{rgb}{0.0, 0.5, 0.0}
\definecolor{jrnl1}{rgb}{0.0, 0.72, 0.6}
\definecolor{aqua}{rgb}{0.0, 1.0, 1.0}
\definecolor{webbrown}{rgb}{.6,0,0}
\definecolor{pinegreen}{rgb}{0.0, 0.47, 0.44}
\definecolor{ultramarineblue}{rgb}{0.25, 0.4, 0.96}
\definecolor{jrnl}{rgb}{0.0, 0.5, 0.0}
\definecolor{lincolngreen}{rgb}{0.11, 0.35, 0.02}
\definecolor{green(html/cssgreen)}{rgb}{0.0, 0.5, 0.0}
\definecolor{airforceblue}{rgb}{0.36, 0.54, 0.66}
\definecolor{azure}{rgb}{0.0, 0.5, 1.0}
\definecolor{bleudefrance}{rgb}{0.19, 0.55, 0.91}
\definecolor{cobalt}{rgb}{0.0, 0.28, 0.67}

\hypersetup{colorlinks=true, linktocpage=true, pdfstartpage=1, pdfstartview=FitV,
	breaklinks=true, pdfpagemode=UseNone, pageanchor=true, pdfpagemode=UseOutlines,
	plainpages=false, bookmarksnumbered, bookmarksopen=true, bookmarksopenlevel=1,
	hypertexnames=true, pdfhighlight=/O,
	urlcolor=webbrown, linkcolor=cobalt, citecolor=jrnl
}

\newcommand{\bmu}{\boldsymbol{\mu}}
\newcommand{\om}{\overline{M}}

\newcommand{\ch}[1]{{#1}}
\newcommand{\chh}[1]{{#1}}

\newcommand{\chhh}[1]{\textcolor{black}{{#1}}}
\newcommand{\css}[1]{\textcolor{black}{{#1}}}

\newif\iffinal
\finalfalse	
\finaltrue	
\iffinal\else
\usepackage[notref,notcite]{showkeys}\fi

\begin{document}

	\title[Random graph models at criticality]{Scaling limits of random graph models at criticality: Universality and the basin of attraction of the Erd\H{o}s-R\'enyi random graph}
	
	\date{}
	
	\author[Bhamidi]{Shankar Bhamidi$^1$}
	\address{$^1$Department of Statistics and OR, University of North Carolina, Chapel Hill}
	\author[Broutin]{Nicolas Broutin$^2$}
	\address{$^2$Sorbonne Universit\'{e}, France,}
	\author[Sen]{Sanchayan Sen$^3$}
	\address{$^3$Department of Mathematics, Indian Institute of Science}
	\author[Wang]{Xuan Wang$^4$}
	\address{$^4$Databricks}
	\maketitle
\begin{abstract}
A wide array of random graph models have been postulated to understand properties of observed networks. Typically these models have a parameter $t$ and a critical time $t_c$ when a giant component emerges. It is conjectured that for a {large} class of models,  the nature of this emergence {is similar to that of} the \erdos random graph, in the sense that 
(a) the sizes of the maximal components in the critical regime scale like $n^{2/3}$, and 
(b) the structure of the maximal components at criticality (rescaled by $n^{\sss -1/3}$) converges to random fractals. 
To date, (a) has been proven for a number of models using different techniques. 
This paper develops a general program for proving (b) that requires three ingredients: (i) in the critical scaling window, components merge approximately like the multiplicative coalescent, (ii) scaling exponents of susceptibility functions are the same as that of the \erdos random graph, and (iii) macroscopic averaging of distances between vertices in the barely subcritical regime.
We show that these apply to two fundamental random graph models: the configuration model and inhomogeneous random graphs with a finite ground space.
For these models, we also obtain new results for component sizes at criticality and structural properties in the barely subcritical regime.
\end{abstract}

\tableofcontents

\section{Introduction}\label{sec:intro}
Over the last few years, motivated both by questions in fields such as combinatorics and statistical physics as well as the explosion in data on real networks, an array of random graph models have been proposed. 
Writing $[n] = \set{1,2,\ldots, n}$ for the vertex set, most of these models have a parameter $t$ (related to the edge density) and a model dependent critical time $t_c$ such that for $t< t_c$ (subcritical regime), size of the largest component is $o_P(n)$ (i.e., there exists no giant component),
while for $t> t_c$ (supercritical regime), the size of the largest component scales like $f(t)n$ where $f(t) > 0$ and is model dependent. The canonical example of such phenomenon is the \erdos random graph $\ERRG(n,t/n)$ where the critical time $t_c(\ERRG)=1$.
One fundamental question in the study of these models is understanding the nature of this phase transition.

Asymptotics in the critical regime are conjectured to be ``universal'' \cite{braunstein2003optimal}. 
More precisely, under moment assumptions on the average degree, a wide array of models are conjectured to exhibit the same behavior in the critical regime as the \erdos random graph in the sense that for any fixed $i\geq 1$, the size of the $i$-th maximal component scales like $n^{2/3}$, and the typical disctance in the component scales like $n^{1/3}$. 
To date, for \emph{component sizes} in the critical regime, such behavior has been proven for a number of models including the rank-one \chhh{inhomogeneous} random graph \cite{SBVHJL10},  configuration model \cite{joseph2010component,riordan2012phase}, and bounded size rules \cite{bhamidi2014augmented,bhamidi2015aggregation}. 
Turning to the intrinsic geometry of the components, it was shown in \cite{addario2012continuum} that rescaling edge lengths in the $i$-th maximal component $\cC_i(\lambda)$ of the \erdos random graph $\ERRG(n , 1/n+\lambda/n^{4/3})$ and denoting the resulting metric space by $n^{-1/3}\cC_i(\lambda)$, one has
\begin{equation}
	\label{eqn:lim-add-br-go}
	\big(n^{-1/3}\cC_i(\lambda);\, i\geq 1\big) \weakc \vCrit(\lambda):=(\crit_i(\lambda);\, i\geq 1),
\end{equation}
for a sequence of limiting random fractals that are described in more detail in Section \ref{sec:prelim}.

Invariance and universality principles play a fundamental role in probability; for example, 
Donsker's invariance principle studies the convergence of processes to Brownian motion under \emph{uniform asymptotic negligibility} conditions. This paper establishes general conditions for random graph models in the critical regime to satisfy asymptotics such as \eqref{eqn:lim-add-br-go}. General theory coupled with model specific analysis enables us to establish the metric space scaling limit for two major families of random graph models: 
the configuration model ({\bf CM}) and inhomogeneous random graphs ({\bf IRG}) with a finite ground space. 
We further develop new techniques to study component sizes of inhomogeneous random graphs \cite{BBSJOR07} at criticality. 
As part of the proof, we derive various estimates in the \emph{barely subcritical regime} for these two models, which are of independent interset.  
In future work, we will show how the techniques developed in this paper extend to a number of other models including \ch{Achlioptas processes with} bounded size rules \cite{spencer2007birth}, which are known to exhibit \erdos type phase transition \cite{bhamidi2014augmented, SBABXW14, sen2013largest}.  
As shown in \cite{addario2013scaling}, metric structure of maximal components at criticality is the first step in analyzing more complicated objects such as the minimal spanning tree on the giant component in the supercritical regime.  The techniques developed in this paper would be the first step in establishing universality for such objects as well.    

{\bf Organization of the paper:}
In Section \ref{sec:models} we recall the definitions of the configuration model and the inhomogeneous random graph model.
In Section \ref{sec:prelim}, we define mathematical constructs needed to state the results. Section \ref{sec:res} describes the general universality result.  Section \ref{sec:res-models} contains results for the two families of random graph models.  \chh{The relevance of these results are discussed in Section \ref{sec:disc}}. Section \ref{sec:proofs-universality} contains proofs of the universality results. Section \ref{sec:proof-irg} and Section \ref{sec:proof-cm}  contain proofs for IRG and CM  respectively.

\subsection{Models}
\label{sec:models}
The vertex set for both models will be $[n]$.  Suppressing dependence on $n$, write $\cC_i$ for the $i$-th maximal component and $\cC_i(t)$ {for the corresponding object} in a dynamic random graph process at time $t$.

\subsubsection{Inhomogeneous random graph \cite{BBSJOR07}}
\label{sec:irg-def}
Start with a ground space or type space $(\cX, \mu)$, where $\cX$ is a separable metric space, and $\mu$ is a Borel probability measure $\mu$ on $\cX$. 
Let $\kappa:\cX \times \cX\to \bR_+$ be a symmetric nonnegative kernel. 
Further, assume that we are given a sequence $(\vx_n; n\geq 1)$, where $\vx_n=(x_1,\ldots, x_n)\in\cX^n$ such that $n^{-1}\sum_{i\in[n]}\delta_{\{x_i\}}\to\mu$ as $n\to\infty$.
Consider the random graph with vertex set $[n]$ constructed by connecting $i,j\in [n]$ with probability $\min\{1, \kappa(x_i,x_j)/n\}$ independent across edges. 
In our regime, this model is asymptotically equivalent \cite{janson2010asymptotic} to the model with connection probabilities
$
1-\exp\big(-\kappa(x_i, x_j)/n\big)
$.
We use \chhh{both versions in the course of the paper}. 
We restrict ourselves to the case where the type space is finite: $\cX = [K] := \set{1,2,\ldots, K}$ for some $K\geq 1$. 
Thus, the kernel $\kappa$ is represented by a $K\times K$ symmetric matrix, and $\mu = (\mu(x), x\in [K])$ is a probability mass function. To avoid irreducibility issues, we assume $\kappa(x,y)> 0$ and $\mu(x)> 0$ for all $x, y\in [K]$. Define {the} operator  $T_{\kappa}: L^2([K],\mu) \to L^2([K],\mu)$:
\[(T_\kappa f)(x):= \sum_{y=1}^K  \kappa(x, y)f(y)\mu(y), \qquad x\in [K].\]
Write $||\kappa||$ for the operator norm of $T_{\kappa}$.
{In} \cite{BBSJOR07} it is shown that if $||\kappa|| < 1$, $|\cC_1| = {O_P}(\log{n})$, while if $||\kappa||> 1$, then $|\cC_1|\sim \rho(\kappa,\mu) n$, where $\rho(\kappa,\mu)$ is the survival probability of an associated supercritical multitype branching process. 
We will work in the critical regime where $||\kappa||=1$.

\noindent {\bf Dynamic version}: For each unordered pair of distinct vertices $u,v\in [n]$, generate {i.i.d.} rate one exponential random variables $\xi_{uv}$. For fixed $\lambda\in \bR$, form the graph $\cG_n^{\sss \IRG}(\lambda)$ by connecting each pair of vertices $u,v$ if
$\xi_{uv}\leq \left(1+{\lambda}{n^{-1/3}}\right){\kappa(x_u, x_v)}/{n}$,
where $x_u, x_v \in [K]$ denote the types of $u$ and $v$ respectively. 
By construction, $\cG_n^{\sss \IRG}(\lambda)$ is a random graph where edges are placed independently with probability
\begin{equation}
	\label{eqn:puv-def}
	p_{uv}(\lambda):= 1-\exp\Big(-\big(1+\lambda n^{-1/3}\big)\kappa(x_u, x_v)/n\Big), \qquad u\neq v\in[n].
\end{equation}
Further, the entire dynamic process $\big(\cG_n^{\IRG}(\lambda), -\infty < \lambda < \infty\big)$  ``increases'' in the sense that $\cG_n^{\IRG}(\lambda_1)$ is a subgraph of $\cG_n^{\IRG}(\lambda_2)$ if $\lambda_1<\lambda_2$.

\subsubsection{ Configuration model \cite{B80,BC78,MR98}:}
\label{sec:cm-def}
Our results apply to {\bf critical percolation} on the supercritical configuration model. 
We now describe the model.
Start with a given degree sequence $\vd_n = (d_i, i\in [n])$, where $d_i\in \bN\cup\{0\}$ with $\sum_{i=1}^n d_i$ even. 
Think of each vertex $i\in [n]$ as having $d_i$ many half-edges associated with it.  
Construct the random  multigraph $\CM_n(\vd_n)$ by performing a uniform matching of these half-edges (two half-edges forming a complete edge). 
Special cases include:
\begin{inparaenuma}
	\item {\it Random $r$-regular \chhh{multigraph:}} Fix $r\geq 3$ and let $d_i = r$ for all $i\in [n]$.
	\item {\it Uniform simple random graph:} Conditioned on simplicity, $\CM_n(\vd_n)$ follows the uniform distribution on the space of simple graphs with degree sequence $\vd_n$.
\end{inparaenuma}
The results in this paper apply directly to the multigraph setting; arguing as in \cite{joseph2010component}, one can extend these results to the uniform simple random graph setting. 
We {do not include this argument here to keep this paper under a manageable length}. 

Now define $\nu_n:= {\sum_{i=1}^n d_i(d_i-1)}/{\sum_{i=1}^n d_i}$. 
Under the assumption $\nu_n \to \nu > 1$ and some additional regularity conditions on the sequence $\big(\vd_n;\, n\geq 1\big)$,
\cite{MR98} proves existence of a unique giant component with size $\sim \rho n$ for some $\rho > 0$. Call this the supercritical regime. 
Let $\Perc_n(p)$ be the random graph obtained by performing {\bf bond percolation} on $\CM_n(\vd_n)$ with edge retention probability $p$.
It is well-known \cite{janson-percolation-CM,fountoulakis2007,MR98} that the critical value for the existence of a giant component is $p_c = 1/\nu$. 
We derive the metric space structure of the maximal connected components in the critical window, i.e., when $ p=p(\lambda)$ is given by
\begin{equation}\label{eqn:p-lambda-cm}
	p(\lambda) := \frac{1}{\nu} +\frac{\lambda}{n^{1/3}}
\end{equation}
for $\lambda\in\bR$ fixed.

\noindent {\bf Dynamic version:} Start with the vertex set $[n]$, where the vertex $i$ has $d_i$ many half-edges attached to it.
Assign {to} every half-edge an exponential rate one clock.  At time $t=0$, all half-edges are designated {\bf alive}.  
When the clock of an alive half-edge rings, this half-edge selects another alive half-edge uniformly at random and forms a full edge. 
Both half-edges are then considered {\bf dead} and {are} removed from the collection of alive half-edges. 
This construction is related but not identical to the dynamic construction in \cite{janson2009new}.  Let $(\CM_n(t), t\geq 0)$ denote this dynamic random graph process. {Then $\CM_n(\infty)$ has the same law as the random graph $\CM_n(\vd_n)$ constructed above.} Further, {we shall see that} this dynamic process stopped at the ``critical'' times
\[t(\lambda):= \frac{1}{2}\log\frac{\nu}{\nu-1} + \frac{\lambda \nu}{2(\nu-1)}\frac{1}{n^{1/3}}\]
gives information about the random graph $\Perc_n(p(\lambda))$.

\section{Preliminaries}
\label{sec:prelim}


\subsection{GHP convergence of compact metric measure  spaces}
\label{sec:gh-mc}

We follow \cite{EJP2116,addario2013scaling,burago2001course}. 
For two metric spaces $X_1 = (X_1,d_1)$ and $X_2 = (X_2, d_2)$ and $C\subseteq X_1 \times X_2$, define the distortion of $C$ with respect to $X_1$ and $X_2$ as
\begin{equation}
	\label{eqn:def-distortion}
	\dis(C) =\dis(C; X_1, X_2)
	:= \sup \big\{|d_1(x_1,y_1) - d_2(x_2, y_2)|: (x_1,x_2) , (y_1,y_2) \in C\big\}.
\end{equation}
A correspondence $C$ between $X_1$ and $X_2$ is a measurable subset of $X_1 \times X_2$ such that for every $x_1 \in X_1$ there exists at least one $x_2 \in X_2$ such that $(x_1,x_2) \in C$ and vice-versa. Let $\text{Corr}(X_1,X_2)$ be the set of all such correspondences.  The Gromov-Hausdorff distance between $(X_1,d_1)$ and $(X_2, d_2)$ is defined as
\begin{equation}
	\label{eqn:dgh}
	d_{\GH}(X_1, X_2) = \frac{1}{2}\inf \big\{\dis(C): C \in \text{Corr}(X_1, X_2)\big\}.
\end{equation}

The Gromov-Hausdorff-Prokhorov distance extends the above to keep track of measures on the corresponding spaces.
A compact metric measure space is a compact metric space equipped with a finite measure on its Borel sigma algebra.
Let $X_1 = (X_1, d_1, \mu_1)$ and $X_2=(X_2,d_2, \mu_2)$ be compact metric measure spaces.
Suppose $\pi$ is a measure on the product space $X_1\times X_2$ with marginals $\pi_1, \pi_2$. 
The discrepancy of $\pi$ with respect to $X_1$ and $X_2$ is defined as
\begin{equation}
	\label{eqn:def-discrepancy}
	\dsc(\pi)=\dsc(\pi; X_1, X_2):= ||\mu_1-\pi_1|| + ||\mu_2-\pi_2||,
\end{equation}
where $||\cdot||$ denotes the total variation of signed measures. {Now define:}
\begin{equation}
	\label{eqn:dghp}
	\bar d_{\GHP}(X_1, X_2):= \inf\Big\{ \max\big(\dis(C)/2, ~\dsc(\pi),~\pi(C^c)\big) \Big\},
\end{equation}
where the infimum is taken over all correspondences $C$ and measures $\pi$ on $X_1 \times X_2$. The function $\bar d_{\GHP}$ is a pseudometric and defines an equivalence relation $X_1\sim X_2 \Leftrightarrow \bar d_{\GHP}(X_1,X_2) = 0$ on. Let $\sS $ be the space of all equivalence classes of compact metric measure spaces and let $d_{\GHP}$ be the induced metric. Then by \cite{EJP2116}, $(\sS, d_{\GHP})$ is a complete separable metric space. 
We will continue to use $X = (X, d, \mu)$ to denote both the metric space and {the} corresponding equivalence class. In this paper we will typically deal with an
infinite sequence of metric spaces. The relevant space in this case is $\sS^{\bN}$ equipped with the product topology inherited from $d_{\GHP}$. 

\textbf{The scaling operator:} We will need to rescale both the metric as well as associated measures of the components in the critical regime. For $\alpha, \beta > 0$, let $\scl(\alpha,\beta)$ be the scaling operator
\begin{align*}
	\scl(\alpha, \beta) : \sS &\to \sS, \qquad \scl(\alpha, \beta)[(X , d , \mu)]:= (X, d',\mu'),
\end{align*}
where $d'(x,y) := \alpha d(x,y)$ for all $x,y \in X$, and $\mu'(A) := \beta \mu(A)$ for all $A \subseteq X$. For simplicity, write $\scl(\alpha, \beta) X := \scl(\alpha, \beta)[(X , d , \mu)]$ and $\alpha X := \scl(\alpha, 1) X$.


\subsection{Notation}
\label{sec:gr-constr}
We write $|A|$ or $\# A$ for the cardinality of a set $A$.
For {a} finite graph $\cG$, we write $V(\cG)$ and $E(\cG)$ for the set of vertices and the set of edges in $\cG$ respectively.
We write $|\cG|$ to mean $|V(\cG)|$, and
$\spls(\cG) := |E(\cG)| - |\cG|+1$ for the number of surplus edges in $\cG$.
We view a connected component $\cC$ of $\cG$ as a metric space using the graph distance $d_{\cG}$.
When the graph $\cG$ is clear from the context, we will supress the subscript and simply write $d$ for the graph distance. 
Often in our analysis, associated to a graph $\cG$ there will be a collection of vertex weights 
$\vw = \big\{w_v: v\in V(\cG)\big\}$.
Natural measures on $\cG$ are
\begin{inparaenuma}
	\item {\it Counting measure:} $\mu_{\text{ct}}(A): = |A|$, for $A \subseteq V(\cG)$; and
	\item {\it Weighted measure:} $\mu_\vw(A) :=  \sum_{v \in A} w_v$, for $A \subseteq V(\cG)$. 
\end{inparaenuma}
For $\cG$ finite and connected, we will often use $\cG$ for both the graph and its associated metric measure space.	
If the measure to be considered is not the counting measure then that will be explicitly mentioned.
If nothing is specified, the measure will be understood to be the counting measure.
For two graphs $\cG_1$ and $\cG_2$, we write $\cG_1\subseteq \cG_2$ to mean that $\cG_1$ is a subgraph of $\cG_2$.
For a metric space $(\fX, d)$, write $\diam(\fX):=\sup_{u, v\in\fX}d(u, v)$.
For a finite graph $\cG$, we let 
$\diam(\cG):=\max\big\{\diam(\cC)\, :\, \cC\text{ connected component of }\cG\big\}$, and 
$\girth(\cG)$ to be the length of a shortest cycle in $\cG$.
For two sequence $X_n$ and $Y_n$ of random variables, we will write $X_n\sim Y_n$ to mean that $X_n/Y_n\weakc 1$, as $n\to\infty$.
(We also write $X\sim\mu$ to mean that the random variable $X$ has law $\mu$; this, however, should not cause any confusion.)
We will say that a sequence $\cE_n,\, n\geq 1$, of events occurs with high probability (whp) if $\pr(\cE_n)\to 1$.
We will freely omit ceilings and floors when there is no confusion in doing so.

\subsection{Real trees with shortcuts}
\label{sec:cont-limit-descp}
For $l>0$, write
$\cE_l$ for the space of continuous excursions on $[0,l]$. 
Fix $h,g \in \cE_l$, a locally finite set $\cP \subseteq \bR_+\times \bR_+$, and set
\begin{equation*}
	g \cap \cP := \big\{(x,y) \in \cP: 0 \leq x \leq l, \; 0 \leq y < g(x)  \big\}.
\end{equation*}
Let $\cT(h)$ be the real tree encoded by $h$ (see, e.g., \cite{evans-book,legall-book}) equipped with the push forward of the Lebesgue measure on $[0,l]$. 
For $(x,y) \in g \cap \cP$, let $r(x,y) := \inf\set{x': x' \geq x, \; g(x') \leq y}$.
Write $\cG(h,g,\cP)$ for the metric measure space obtained by identifying the pairs of points in $\cT(h)$ corresponding to the pairs of points $\{(x,r(x,y)) : (x,y) \in g \cap \cP\}$.
Thus $\cG(h,g,\cP)$ is obtained by adding a finite number of shortcuts to $\cT(h)$.

\noindent \textbf{Tilted Brownian excursions:} 
Let $(\ve_l(s), s \in [0,l])$ be a Brownian excursion of length $l$. 
Fix $\theta >0$ and let $\tilde \ve_l^\theta$ be an $\cE_l$-valued random variable such that for any bounded continuous function $f : \cE_l \to \bR$,
\begin{equation}
	\label{eqn:tilt-exc-def}
	\E[f(\tilde \ve_l^\theta)] = { \E\Big[f(\ve_l) \exp\Big(\theta\int_0^l \ve_l(s)ds\Big) \Big] }\big/{\E\Big[\exp\Big(\theta\int_0^l \ve_l(s)ds \Big)\Big]}.
\end{equation}
For simplicity, write $\ve(\cdot) := \ve_1(\cdot)$, $\tilde \ve^\theta (\cdot) := \tilde \ve_1^{\theta} (\cdot)$, and $\tilde \ve_l (\cdot) := \tilde \ve_l^1 (\cdot)$.
In Section \ref{sec:1}, we will use the random objects defined in this section to describe a construction of the limiting random metric measure spaces of interest in this paper.

\subsection{Multiplicative coalescent and the random graph $\cG(\vx,q)$}
\label{sec:smc-def}
A multiplicative coalescent $\big(\vX(t);\, t\in \cA\big)$ is a Markov process with state space 
$
\ldown := \{(x_1,x_2,\ldots): x_1\geq x_2\geq \cdots \geq 0, \sum_i x_i^2< \infty\}
$
endowed with the metric inherited from $l^2$;
here, either $\cA=\bR$ or $\cA=[t_0, \infty)$ for some $t_0\in\bR$.
Its evolution can be described in words as follows:
Fix $\vx\in \ldown$.
Then conditional on $\vX(t) = \vx$, each pair of clusters $i$ and $j$ merge at rate $x_i x_j$ to form a new cluster of size $x_i+x_j$. While this description makes sense for a finite collection of clusters (i.e., $x_i = 0$ for $i> K$ for some finite $K$), Aldous \cite{aldous1997brownian} showed that \chh{it also} makes sense for $\vx\in \ldown$.
We refer the reader to \cite{aldous1997brownian} for a more detailed study of how the multiplicative coalescent is connected to the evoloution of the \erdos random graph.

An object closely related to the multiplicative coalescent is the random graph $\cG(\vx,q)$:  Fix vertex set $[n]$, a collection of positive vertex weights $\vx=(x_i, i\in [n])$, and parameter $q>0$. Construct the random graph $\cG(\vx, q)$ by placing an edge between $i\neq j\in [n]$ with probability
$1 - \exp(-q x_i x_j)$,
independently across pairs $\{i, j\}$ with $i\neq j$.
For a connected component $\cC$ of $\cG(\vx, q)$, define 
$\mass (\cC) = \sum_{i \in \cC} x_i=\chhh{\mu_{\vx}(\cC)}$ 
\chhh{(recall the notation from Section \ref{sec:gr-constr})}. 
Rank the components in terms of their masses and let $\cC_i$ be the $i$-th maximal component.

\begin{ass}\label{ass:aldous-basic-assumption}
	Consider a sequence of vertex weights $\vx = \vx^{\sss(n)}$ and a sequence $q = q^{\sss(n)}$.
	Let $\sigma_k : = \sum_{i \in [n]} x_i^k$ for $k = 2, 3$ and $x_{\max} := \max_{i\in [n]} x_i$. 
	Assume that there exists a constant $\lambda \in \bR$ such that as $n \to \infty$,
	\begin{equation*}
		{\sigma_3}/{(\sigma_2)^3} \to 1, 
		\ \ 
		q - {(\sigma_2)^{-1}} \to \lambda, 
		\ \ \text{ and }\ \ 
		{x_{\max}}/{\sigma_2} \to 0.
	\end{equation*}
\end{ass}
Note that these conditions imply that $\sigma_2 \to 0$. 
Fix $\lambda \in \bR$ and let $B$ be a standard Brownian motion.  
Define the processes $W_\lambda$ and $\tilde{W}_\lambda$ via
\begin{equation}
	\label{eqn:parabolic-bm}
	W_\lambda (t) := B(t) + \lambda t - t^2/2\, ,\ \ \text{ and }\ \  
	\tilde W_\lambda (t) := W_\lambda (t) - \inf_{s \in [0,t]} W_\lambda (s)\, ,\ \ t\geq 0\, .
\end{equation}
An excursion of $\tilde W_\lambda$ is an interval $(l,r) \subset \bR^+$ such that $\tilde W_\lambda(l)=\tilde W_\lambda(r) = 0$ and $\tilde W_\lambda(t)>0$ for all $t \in (l,r)$. Write $r-l$ for the length of such an excursion. 
Aldous \cite{aldous1997brownian} showed that the lengths of the excursions of $\tilde W_\lambda$ can be arranged in decreasing order as
$\gamma_1(\lambda) > \gamma_2(\lambda) > \ldots > 0$.
Conditional on $\tilde W_\lambda$, let $\pois_i(\lambda)$, $i\geq 1$, be independent random variables where $\pois_i(\lambda)$ is Poisson distributed with mean equal to the area underneath the excursion with length $\gamma_i(\lambda)$. 

\begin{thm}[\cite{aldous1997brownian}]\label{thm:aldous-review}
	Under Condition \ref{ass:aldous-basic-assumption}, 	$\left(\mass( \cC_i); i \geq 1 \right) \weakc \mvxi(\lambda):= (\gamma_i(\lambda); i\geq 1)$,  as $n \to \infty$, with respect to the topology on $\ldown$.
	Further, 
	$\big(\big(\mass( \cC_i),\spls(\cC_i) \big);\, i \geq 1 \big) 
	\weakc 
	\mvXi(\lambda):=\big( \big(\gamma_i(\lambda), \pois_i(\lambda)\big);\, i\geq 1)$,  as $n \to \infty$, with respect to the product topology.
\end{thm}

\begin{rem}
	In \cite[Proposition 4]{aldous1997brownian}, Aldous only proves the first assertion in Theorem \ref{thm:aldous-review}, and the second assertion in Theorem \ref{thm:aldous-review} is proved for the special case of the \erdos random graph \cite[Corollary 2]{aldous1997brownian}.
	However, the second convergence in Theorem \ref{thm:aldous-review} follows easily from \cite[Proposition 10]{aldous1997brownian}.
\end{rem}

\subsection{Scaling limits of components in critical \erdos random graph}\label{sec:1}
Now we can define the scaling limit of the maximal components in
$\ERRG(n, 1/n+\lambda/n^{4/3})$ derived in \cite{addario2012continuum}. 
Recall the definitions of tilted Brownian excursions and the metric measure space $\cG(h,g,\cP)$ from Section \ref{sec:cont-limit-descp}. 
Let $\mvxi(\lambda)= (\gamma_i(\lambda), i\geq 1)$ be as in Theorem \ref{thm:aldous-review}. Conditional on $\mvxi(\lambda)$, let $\tilde \ve_{\gamma_i(\lambda)}$, $i\geq 1$, be independent tilted Brownian excursions with $\tilde \ve_{\gamma_i(\lambda)}$ having length $\gamma_i(\lambda)$.  Let $\cP_i$, $i\geq 1$, be independent rate one Poisson processes on $\bR_+^2$ that are also independent of $\big(\tilde \ve_{\gamma_i(\lambda)};\,  i\geq 1\big)$. Now consider the sequence of random metric {measure} spaces
\begin{equation}
	\label{eqn:limit-metric-def}
	\crit_i(\lambda):= \cG(2\tilde \ve_{\gamma_i(\lambda)}, \tilde \ve_{\gamma_i(\lambda)}, \cP_i), \qquad i\geq 1.
\end{equation}
For the rest of the paper we let
\begin{equation}
	\label{eqn:vcrit-def}
	\vCrit(\lambda):= (\crit_i(\lambda), i\geq 1).
\end{equation}

\section{Results: Universality}
\label{sec:res}
Our plan is to extend Theorem \ref{thm:aldous-review} in two stages.
In the first stage, we consider the random graph $\cG(\vx, q)$, and for $i\geq 1$, we view the component $\cC_i$ as a metric measure space endowed with the measure $\mu_{\vx}$.
Then \cite[Theorem 7.3]{SBSSXW14} implies that under Condition \ref{ass:aldous-basic-assumption} and Condition \ref{ass:aldous-gen-1} stated below, the sequence $(\cC_i, i\geq 1)$ properly rescaled converges to $\vCrit(\lambda)$ as in \eqref{eqn:vcrit-def}. 
In the second stage of the extension of Theorem \ref{thm:aldous-review}, we replace each vertex $i \in [n]$ in the graph $\cG(\vx,q)$ with a metric measure space $(M_i, d_i, \mu_i)$.
In Theorem \ref{thm:aldous-gen-2}, we show that the metric measure space $\bar\cC_i$ now associated with $\cC_i$, under Conditions \ref{ass:aldous-basic-assumption} and \ref{ass:aldous-gen-2}, converges to $\crit_i(\lambda)$ after proper rescaling, owing to macroscopic averaging of distances within blobs.

Before moving on to the statements of these results, let us explain how this general theorem can be applied to specific random graph models.
Recall the dynamic formulations of the random graph models in Section \ref{sec:models}, and for this discussion denote these by the generic  $\big(\bG_n(t), t\geq 0\big)$. 
For each model,  we fix an appropriate $\delta \in(0, 1/3)$ and set $t_n = t_c - An^{-\delta}$, where $t_c$ is the critical time for the model, and $A$ is a suitable constant.
We refer to the components of $\bG_n(t_n)$ as `blobs.'
We then study the connections that form between the blobs in the time interval $[t_n, t_c+\lambda/n^{1/3}]$. 
For both CM and IRG, we show that in this interval, the connectivity pattern between blobs can be either exactly or \emph{approximately} described through the graph $\cG(\vx,q)$.  
We emphasize that here, we are viewing each blob as a single vertex. 
Further, $x_i$ is a model dependent functional of the $i$-th blob, and $q$ is a model dependent function of $\lambda$.
Thus, $\cG(\vx,q)$ should be thought of as describing the blob-level superstructure owing to edges created in the interval $[t_n, t_c+\lambda/n^{1/3}]$, where each blob is viewed as a single vertex ignoring the internal structure of the blobs.
Now, once we verify that Condition \ref{ass:aldous-gen-2} is satisfied by the blobs, 
Theorem \ref{thm:aldous-gen-2} applies to the random graph $\overline\bG_n(t_c+\lambda/n^{1/3})$ obtained by replacing the vertices of $\cG(\vx,q)$ by the blobs.
This yields the scaling limit of the maximal components of $\bG_n(t_c+\lambda/n^{1/3})$ if the connectivity pattern between blobs were exactly described through the graph $\cG(\vx,q)$.
Otherwise, an additional argument is needed to show that $\bG_n(t_c+\lambda/n^{1/3})$ is well-approximated by $\overline\bG_n(t_c+\lambda/n^{1/3})$.

\subsection{Blob-level superstructure}
\label{sec:stage-one-blob-level}
Recall the random graph $\cG(\vx,q)$ from Section \ref{sec:smc-def}, and the notation used in Condition \ref{ass:aldous-basic-assumption}.
View $(\cC_i, i \geq 1)$ as metric measure spaces equipped with the graph distance and the weighted measure $\mu_{\vx}$. 

\begin{ass}\label{ass:aldous-gen-1}
	Assume that there exist $\eta_0 \in (0,\infty)$ and $r_0 \in (0,\infty)$ such that ${x_{\max}}/{\sigma_2^{3/2+\eta_0}} \to 0$ and  ${\sigma_2^{r_0}}/{x_{\min}} \to 0$
	as $n \to \infty$, where $x_{\min}:=\min_{i\in [n]}x_i$.
\end{ass}

\begin{thm}\label{thm:aldous-gen-1}
	Under Conditions \ref{ass:aldous-basic-assumption} and \ref{ass:aldous-gen-1}, 
	$\big(\scl(\sigma_2, 1 ) \cC_i^{\sss(n)}, i \geq 1 \big) \weakc \vCrit(\lambda)$ as $n\to\infty$. 
\end{thm}

\begin{rem} 
	Theorem \ref{thm:aldous-gen-1} follows from \cite[Theorem 7.3]{SBSSXW14}.
	Alternately, it can be derived as a special case of Theorem \ref{thm:aldous-gen-2} stated below.
	In the proofs, we will use \cite[Theorem 7.3]{SBSSXW14} directly.
	However, we state Theorem \ref{thm:aldous-gen-1} in the above form as this explains the motivation behind us proceeding in two stages as discussed above.
	Note also that the scaling limit of the critical \erdos random graph $\ERRG(n, n^{-1}+\lambda n^{-4/3})$ derived in \cite{addario2012continuum} can be recovered from Theorem \ref{thm:aldous-gen-1} by taking 
	$q=n^{1/3}+\lambda$ and $x_i = n^{-2/3}$ for $i \in [n]$
	(which results in $\sigma_2 = n^{-1/3}$). 
\end{rem}

\subsection{Incorporating the internal structure of the blobs}
\label{sec:inter-blob-distance}
We need the following ingredients:

\noindent{\upshape (a)} {\bf Blob level superstructure: }  A simple finite graph $\cG$ with vertex set $[n]$ and vertex weight sequence $\vx:=(x_i, i \in [n])$.

\vskip3pt

\noindent{\upshape (b)} {\bf Blobs:} A family of compact metric measure spaces 
$\vM := \big\{(M_i, d_i, \mu_i): i \in [n] \}$, 
where $\mu_i$ is a probability measure for all $i$.

\vskip3pt

\noindent{\upshape (c)} {\bf Blob to blob junction points:} A collection of points $\vX := \{X_{i,j}: i \in [n], j \in [n]\}$ such that $X_{i,j} \in M_i$ for all $i,j$.

\vskip3pt

Using these three ingredients,  we can define a metric measure space $\Gamma(\cG,\vx,\vM,\vX) = (\bar M, \bar d, \bar \mu)$ as follows: Let $\bar M := \bigsqcup_{i \in [n]} M_i$. Define the measure $\bar \mu$ as
\begin{equation}
	\bar \mu( A ) = \sum_{i \in [n]} x_i \mu_i(A \cap M_i), \mbox{ for } A \subseteq \bar M.
	\label{eqn:barmu-on-full-met}
\end{equation}
The metric $\bar d$ is obtained by using the intra-blob
distance functions $\big(d_i;\, i\in [n]\big)$ together with the graph distance on $\cG$ by putting an edge of length one between the pairs of vertices
$\{ \{X_{i,j}, X_{j,i} \}: \{i,j \} \mbox{ is an edge in } \cG \}$.
Thus, for $x, y \in \bar M$ with $x \in M_{j_1}$ and $y \in M_{j_2}$, 
\begin{equation}\label{eqn:44}
	\bar d(x,y) = \inf\Big\{ k +  d_{j_1}(x, X_{j_1,i_1}) + \sum_{\ell=1}^{k-1} d_{i_\ell}(X_{i_\ell, i_{\ell-1}}, X_{i_\ell, i_{\ell+1}}) + d_{j_2}(X_{j_2, i_{k-1}}, y) \Big\},
\end{equation}
where the infimum is taken over $k\geq 1$ and all paths $(i_0, i_1,\ldots,i_{k-1}, i_k)$ in $\cG$ with $i_0 =j_1$ and $i_k = j_2$. 
(Here, the infimum of an empty set is understood to be $+\infty$.
This corresponds to the case where $j_1$ and $j_2$ do not belong to the same compoenent of $\cG$.)

The above is a deterministic procedure for creating a new metric measure space. 
Now assume that we are provided with a parameter sequence $q^{\sss(n)}$, weight sequence $\vx^{\sss(n)} := (x_i^{\sss(n)}, i \in [n])$, and the family of metric measure spaces $\vM^{\sss(n)} := \big((M_i^{\sss(n)}, d_i^{\sss(n)}, \mu_i^{\sss(n)}), i \in [n]\big)$.
As before, we will suppress the dependence on $n$. 
Let $\cG(\vx, q)$ be the random graph defined in Section \ref{sec:smc-def} constructed using the weight sequence $\vx$ and the parameter $q$.   Let $(\cC_i, i \geq 1)$ denote the  connected components of $\cG(\vx,q)$ ranked in decreasing order of their masses. 
Let $X_{i,j}$, $i,j \in [n]$,  be independent random variables (that are also independent of $\cG(\vx,q)$) such that for each fixed $i$, $X_{i,j}$, $j\in [n]$, are i.i.d. $M_i$-valued random variables that are $\mu_i$ distributed. 
Let $\vX = (X_{i,j},\ i,j \in [n])$.

Let $\bar \cG(\vx,q,\vM) := \Gamma(\cG, \vx, \vM, \vX)$ be the (random) compact metric measure space constructed as above. 
Let $\bar \cC_i$ be the component in $\bar \cG(\vx, q, \vM)$ corresponding to the $i$-th maximal component $\cC_i$ in $\cG(\vx,q)$, viewed as a compact metric measure space as explained in \eqref{eqn:barmu-on-full-met} and \eqref{eqn:44}. 
Define the quantities
\begin{gather}
	u_{i} := \int_{M_i \times M_i} d_i(x,y)  \mu_i (dx) \mu_i (dy), \qquad i\in [n],
	\label{eqn:uik-def}\\
	\mbox{$\tau := \sum_{i \in [n]} x_i^2 u_i$ and ${\diam_{\max}} := \max_{i \in [n]} \diam(M_i)$.}
	\label{eqn:tau-dmax-def}
\end{gather}
\chh{Note that $u_i$ is the expectation of \chhh{the distance between} two blob-to-blob junction points in blob $i$ and $\tau/\sigma_2$ is the weighted average of these ``typical'' distances. }
\begin{ass}
	\label{ass:aldous-gen-2}
	There exist $\eta_0 \in (0,\infty)$ and $r_0 \in (0,\infty)$ such that
	{\upshape (a)} Condition \ref{ass:aldous-gen-1} holds, and
	{\upshape (b)}  as $n \to \infty$, $ (\diam_{\max}\cdot \sigma_2^{3/2-\eta_0})/{ (\tau + \sigma_2)} \to 0$ and  $ {\sigma_2 x_{\max} {\diam_{\max}}}/{\tau} \to 0$. 	
\end{ass}
\begin{thm}
	\label{thm:aldous-gen-2}
	Under Conditions \ref{ass:aldous-basic-assumption} and \ref{ass:aldous-gen-2},
	\begin{equation*}
		\Big(\scl\Big(\frac{\sigma_2^2}{\sigma_2 + \tau}, 1 \Big) \bar \cC_i\, ;\, i \geq 1 \Big) \weakc \vCrit(\lambda), \mbox{ as } n \to \infty.
	\end{equation*}
\end{thm}

\begin{rem} \ch{Depending on whether $\lim_{n \to \infty} {\tau}/{\sigma_2}\in [0, \infty)$ or equals $\infty$, the above theorem deals with different scales.
		If each blob is a fixed connected graph $G$, then $\lim_{n\to\infty}\tau/\sigma_2\in[0,\infty)$.
		(The critical \erdos random graph corresponds to the case where $G$ is a single vertex.)
		In \chh{the} applications below, $\sigma_2\sim n^{\delta-1/3}$ while $\tau \sim n^{2\delta-1/3}$ which corresponds to the case where the above limit equals $\infty$.}
	
\end{rem}

\section{Results for associated random graph models}
\label{sec:res-models}

\subsection{Inhomogeneous random graphs}
\label{sec:irg-res}
Recall the IRG model from Section \ref{sec:irg-def}. We work in a slightly more general setup where the kernel could depend on $n$.
Let $\cX=[K]=\set{1,2,...,K}$ be the type space, and $\kappa_n: [K] \times [K] \to \bR^+$ be symmetric kernels for each $n\geq 1$.
Let $\cV^{\sss(n)}:= [n]$ be the vertex \ch{set} where each vertex $i$ has a ``type'' $x_i\in[K]$. 
The IRG $\cG_{\IRG}^{(n)}:=\cG^{(n)}(\kappa_n, \cV^{\sss(n)})$ on the vertex set $\cV^{\sss(n)}$ is constructed by placing an edge with probability 
$
p_{ij}:=1 - \exp\big( - \kappa_n(x_i, x_j)/n\big)
$
between the vertices $i$ and $j$, independently for each $i, j\in \cV^{\sss(n)}$, $i\neq j$.
Denote the empirical distribution of types by $\mu_n$; that is for each $x\in [K]$, $\mu_n(x)$ is the proportion of vertices of type $x$.
Write $\bmu_n$ for the vector $(\mu_n(1), \hdots, \mu_n(K))^t$.
The associated operator $T_{\kappa_n}$ is given by
\begin{equation*}
	(T_{\kappa_n} f ) (x) := \sum_{y\in [K]} \kappa_n(x, y)f(y)\mu_n(y), \;\; x \in [K], f \in \bR^{[K]}.
\end{equation*}
\begin{ass}\label{ass:irg-strong}
	{\upshape (a)} {\bf Convergence of the kernels: } There exists a kernel $\kappa(\cdot, \cdot): [K]\times [K]\to\bR^+$
	and a matrix $A=(a_{xy})_{x, y\in [K]}$ with real valued (not necessarily positive) entries such that
	$\min_{x,y \in [K]} \kappa(x,y) >0$ and for each $x,y\in [K]$, $ \lim_{n \to \infty} n^{1/3}\left(\kappa_{n}(x,y)-\kappa(x,y)\right)= a_{xy}$.
	
	\vskip3pt
	
	\noindent{\upshape (b)} {\bf Convergence of the empirical measures: }There exists a probability measure
	$\mu$ on $[K]$ and a vector $\vb=(b_1,\hdots, b_K)^t$ such that
	$\min_{x\in [K]}\mu(x)>0$, and $\lim_{n \to \infty} n^{1/3}\big(\mu_n(x) - \mu(x)\big)=b_x$ for each $x\in [K]$.
	
	\vskip3pt
	
	\noindent{\upshape (c)} {\bf Criticality}: The operator norm of $T_\kappa$ in $L^2([K], \mu)$ equals one.
\end{ass}

Fix $\delta \in (1/6, 1/5)$ and define the kernel $\kappa_n^-$ by
\begin{equation}\label{eqn:def-kappa^-}
	\kappa_n^-(x,y) := \kappa_n(x,y) - n^{-\delta}, \mbox{ for } x, y \in [K].
\end{equation}
\chh{The IRG $\cG_{\text{IRG}}^{(n),-}:=\cG^{(n)}(\kappa_n^-, \cV^{(n)})$ is barely subcritical. 
	Its components will play the role of the blobs in the proof.}
We will write $\kappa$ (resp. $\kappa_n, \kappa_n^-$) for the $K\times K$ matrix with entries $\kappa(i, j)$
(resp. $\kappa_n(i, j), \kappa_n^-(i, j)$). It will be clear from the context
whether the reference is to the kernel or the matrix.
Define $\mvmu$ to be the vector $(\mu(1),\hdots,\mu(K))^t$, and let
$D:=\mathrm{Diag}(\bmu)$--the $K\times K$ diagonal matrix with diagonal entries $\mu(1), \hdots, \mu(K)$. 
Similarly define $D_n:=\mathrm{Diag}(\bmu_n)$, and $B:=\mathrm{Diag}(\vb)$.
For a square matrix $\om$ with positive entries, write $\rho(\om)$ for its Perron root.
Define $m_{ij}:=\mu(j)\kappa(i,j)$ for $i, j\in [K]$ and let $M=\left(m_{ij}\right)_{i,j \in [K]}$.
{Then} Condition \ref{ass:irg-strong} (c) is equivalent to $\rho(M)=1$. Let $\vu$ and $\vv$
{respectively be the} right and left eigenvectors of $M$ corresponding to {the eigenvalue} $\rho(M)=1$ normalized so that $\vv^t\vu=1$
and $\vu^t\vone=1$, \ie
\begin{align}\label{eqn:defn-u-v}
	M\vu=\vu\, ,\ \ \vv^t M=\vv^t,\ \ \vu^t\vone=1\, , \ \ \text{ and }\ \ \vv^t\vu=1.
\end{align}
Writing $\vu=(u_1,\hdots, u_K)^t$ and $\vv=(v_1,\hdots, v_K)^t$, define
\begin{align}\label{eqn:defn-alpha-beta-zeta}
	\alpha:=\frac{1}{(\vv^t\vone)\cdot(\bmu^t\vu)}\, ,\;\;
	\beta:=\frac{\sum_{x\in[K]}v_x u_x^2}{(\vv^t\vone)\cdot(\bmu^t\vu)^2}\, , \ \text{ and }\
	\zeta:=\alpha\cdot\left[\vv^t(AD+\kappa B)\vu\right].
\end{align}
\begin{thm}
	\label{thm:scaling-limit-irg}
	{For $i\geq 1$}, let $\cC_i(\cG_{\IRG}^{\sss(n)})$ be the $i$-th largest component in $\cG_{\IRG}^{\sss(n)}$, and view it as a metric measure space using the counting measure. Under Condition \ref{ass:irg-strong}, as $n \to \infty$,
	\begin{equation*}
		\Big(\scl\Big( \frac{\beta^{2/3} }{\alpha n^{1/3}}, \frac{\beta^{1/3}}{n^{2/3}} \Big) \cC_i(\cG_{\IRG}^{\sss(n)}),\ i \geq 1 \Big) \weakc \vCrit\Big(\frac{\zeta}{\beta^{2/3}}\Big).
	\end{equation*}
\end{thm}
As a by-product, we obtain the following result
about component sizes.
\begin{thm}
	\label{thm:component-limit-irg}
	Under Condition \ref{ass:irg-strong}, {with respect to the topology on $\ldown$,}
	\begin{equation*}
		\Big(n^{-2/3}\beta^{1/3}\big|\cC_i(\cG_{\IRG}^{\sss(n)})\big|,\ i \geq 1 \Big)
		\weakc
		\mvxi\big(\beta^{-2/3}\zeta\big) \text{ as }n \to \infty.
	\end{equation*}
\end{thm}


\chh{
	The application of Theorem \ref{thm:aldous-gen-2} to
	obtain Theorem \ref{thm:scaling-limit-irg} requires bounds on the barely subcritical graph $\cG_{\text{IRG}}^{(n),-}=\cG^{(n)}(\kappa_n^-, \cV^{(n)})$, whose components will form the blobs. It is technically easier to prove these results for a closely related model, $\cG_{\text{IRG}}^{(n),\star}$, defined as follows. For each $i, j \in [n]$ with $i\neq j$, place an edge independently with probability $p_{i,j}^{\star} := 1 \wedge (\kappa_n^-(x_i, x_j )/n)$. We can transfer the results over to $\cG_{\text{IRG}}^{(n),-}$ via asymptotic equivalence \cite{janson2010asymptotic}}.

For a graph $G = (V(G), E(G))$ and $k\geq 1$, let $\cS_k(G) := \sum_{\cC \text{ component of } G}|\cC|^k$. \chh{Recall that the $k$-th susceptibility function for component sizes is defined as} $s_k(G) := \cS_k(G)/|V(G)|$. Let $\cD(G) := \sum_{i,j \in V(G)} d(i,j) \ind_{\set{d(i,j) < \infty}}$ and $\bar \cD(G) = \cD(G)/|V(G)|$. Let $\cC_i^{\star}$ be the $i$-th largest component in $\cG_{\IRG}^{(n), \star}$, $\cD_{\max}^{\star} = \max_{i \geq 1} \diam(\cC_i^{\star})$,  $\bar\cD^{\star}:= \bar \cD(\cG_{\IRG}^{(n), \star})$ and $\bar s_k^{\star} := s_k(\cG_{\IRG}^{(n), \star})$. Let the analogous objects in $\cG_{\IRG}^{\sss(n), -}$ be  $\cC_i^{-}$,  $\cD_{\max}^{-} = \max_{i \geq 1} \diam(\cC_i^{-})$, $\bar\cD:= \bar \cD(\cG_{\IRG}^{(n), -})$, and $\bar s_k := s_k(\cG_{\IRG}^{(n), -})$ respectively.

\begin{thm}
	\label{prop:irg-barely-subcritical}
	For $l=1,2$, there exist finite constants $A_l=A_l(\kappa, \mu)$ such that
	\begin{gather}\label{eqn:prop-barely-sub-irg-1}
		\frac{\bar s_3^{\star}}{(\bar s_2^{\star})^3} \weakc \beta, \; \;  n^{1/3}\left(\frac{1}{n^{\delta}}-\frac{1}{\bar s_2^{\star}}\right) \weakc \zeta,
		\; \; \pr(|\cC_1^{\star}|\geq A_1 n^{2\delta}\log n)\to 0,\\
		\label{eqn:prop-barely-sub-irg-2}
		{n^{-2\delta}}{\bar \cD^{\star}}\weakc\alpha\, ,\ \text{ and }\ 
		\pr(\cD_{\max}^{\star}\geq A_2 n^{\delta}\log n)\to 0.
	\end{gather}
\end{thm}
A simple consequence of this theorem is the analogous result for $\cG_{\IRG}^{\sss(n), -}$.
\begin{cor}
	\label{cor:irg-barely-subcritical-original}
	The conclusions of Theorem \ref{prop:irg-barely-subcritical} hold if we replace
	$\cC_i^{\star}$, $\cD_{\max}^{\star}$, $\bar\cD^{\star}$, $\bar s_2^{\star}$, and $\bar s_3^{\star}$ by $\cC_i^{-}$, $\cD_{\max}^{-}$, $\bar\cD$, $\bar s_2$, and $\bar s_3$ respectively.
\end{cor}

\subsection{Percolation on supercritical configuration model}
\label{sec:res-cm}
Recall $\CM_n(\vd_n)$ constructed from a degree sequence $\vd_n$ {as} in Section \ref{sec:cm-def}. 
Assume that the degree sequence $\vd_n$ is generated in an i.i.d. fashion using a probability mass function (pmf) $\vp = (p_k, k \geq 1)$; 
if $\sum_{i\in[n]} d_i$ is odd, then we add an extra half-edge to a vertex with the maximum degree.
While this i.i.d. assumption is not essential, it will simplify the statements of the results. 
Let
\begin{equation}
	\label{eqn:def-param}
	\mu := {\sum_{k=1}^\infty k p_k}\, , \ \ 
	\nu := {\sum_{k=1}^\infty k(k-1)p_k}/{\mu}\, , \ \ \text{ and }\ \  
	\beta = \sum_{k=1}^\infty k(k-1)(k-2) p_k\, .
\end{equation}
\begin{ass}
	\label{ass:cm-degree}
	Assume that $\nu> 1$,   $0< \beta <\infty$, and that the degree distribution has exponential moments, i.e., there exists $\Delta_0> 0$ such that $\sum_k \exp(\Delta_0 k) p_k < \infty$.
\end{ass}
Let $\Perc_n(\cdot)$ and $p(\lambda)$ be as in Section \ref{sec:cm-def}. 
Let $\cC_{\sss i, \Perc}(\lambda)$ be the $i$-th maximal component in $\Perc_n(p(\lambda))$, and view it as a metric measure space using the counting measure.

\begin{thm}\label{thm:perc-cm}
	Fix $\lambda \in \bR$ and a sequence $\lambda_n\to\lambda$.
	Then as $n\to\infty$,
	\begin{align}\label{eqn:1}
		\Big(\scl\Big(
		\frac{\beta^{2/3}}{\mu\nu n^{1/3}}\, ,\,
		\frac{\beta^{1/3}}{\mu n^{2/3}}\Big) 
		\cC_{\sss i,\Perc} (\lambda_n)\,;\, i\geq 1
		\Big) 
		\weakc 
		\vCrit\Big(\frac{\mu\nu^2}{\beta^{2/3}}\lambda\Big). 
	\end{align}
\end{thm}

As a {special case},  consider the random $r$-regular multigraph $\cR_r$, where $r\geq 3$.  
Here, $\beta = r(r-1)(r-2)$, $\mu=r$, $\nu=r-1$, and $p(\lambda) = (r-1)^{-1} + \lambda n^{-1/3}$. 
Denote by $\cC_{\sss i,\Perc}^{\sss(\cR_r)}(\lambda)$ the $i$-th largest component of the random graph obtained by performing bond percolation on $\cR_r$ with edge retention probability $p(\lambda)$, and view $\cC_{\sss i,\Perc}^{\sss(\cR_r)}(\lambda)$ as a metric measure space.
Let $\kappa_{1,r} := {(r-2)^{2/3}}/{(r(r-1))^{1/3}} $ and $\kappa_{2,r} := {((r-1)(r-2))^{1/3}}/{r^{2/3}} $.
\begin{cor}\label{cor:random-r-scaling-lim}
	Fix $r\geq 3$, $\lambda \in \bR$, and a sequence $\lambda_n\to\lambda$.
	Then as $n\to\infty$,
	\[
	\left( \scl \left(\frac{\kappa_{1,r}}{n^{1/3}}, \frac{\kappa_{2,r}}{n^{2/3}} \right) 
	\cC_{\sss i,\Perc}^{\sss(\cR_r)}(\lambda)\,;\, i\geq 1\right) 
	\weakc 
	\vCrit\Big(\frac{r^{1/3}(r-1)^{4/3}}{(r-2)^{2/3}}\lambda\Big).  
	\]
\end{cor}

The proof of Theorem \ref{thm:perc-cm} relies on the asymptotics of the dynamic version
$(\CM_n(t)\, ;\, t\geq 0)$ as defined in Section \ref{sec:cm-def}.
As we will see later, the critical time in this process \chh{for the emergence of the giant component} is
\begin{equation}
	\label{eqn:crit-time}
	t_c := 2^{-1}\log{\left({\nu}/{(\nu-1)}\right)}.
\end{equation}
The next theorem gives the critical scaling limit of this process.
\begin{thm}\label{thm:crit-main-res-cm}
	Fix $\lambda \in \bR$ and a sequence $\lambda_n\to\lambda$.
	For $t\geq 0$, let $\cC_i(t)$ be the $i$-th largest component of $\CM_n(t)$, and view it as a metric measure space using the counting measure. 
	Then as $n\to\infty$,
	\[
	\Big(\scl\Big(\frac{\beta^{2/3}}{\mu\nu n^{1/3}},\frac{\beta^{1/3}}{\mu n^{2/3}}\Big)\cC_i\Big(t_c+\frac{\lambda_n}{n^{1/3}}\Big)\, ;\, i\geq 1\Big)
	\weakc
	\vCrit\Big(\frac{2\mu\nu(\nu-1)}{\beta^{2/3}}\lambda\Big).  
	\]
\end{thm}

\ch{The proof of Theorem \ref{thm:crit-main-res-cm} will require bounds in the barely subcritical regime for analogs of the standard susceptibility functions. 
	We first motivate the need for these objects}.
Fix $\delta \in (0, 1/3)$, and define
\begin{equation}
	\label{eqn:tn-def}
	t_n  := t_c -  \nu /{(2(\nu-1)n^{\delta}})\, .
\end{equation}
The constant in front of $n^{-\delta}$ is 
useful in simplifying constants.  
For a component at time $t$, say $\cC_i(t)$, some of the half-edges attached to its vertices may have been used up to form connections in the time interval $[0, t]$. 
Recall that such half-edges are considered dead.
Let $f_i(t)$ be the total number of still alive half-edges in $\cC_i(t)$. 
Note that \chh{two distinct} components $\cC_i(t)$ and $\cC_j(t)$ merge in
$[t,t+dt)$ if one of the free half-edges in $\cC_i(t)$ rings and chooses
one of the free half-edges in $\cC_j(t)$ or vice-versa.
Thus, \chh{the} rate of merger is
\[
f_i(t)\frac{f_j(t)}{n\bars_1(t)-1} + f_j(t) \frac{f_i(t)}{n\bars_1(t)-1} = 2\frac{f_i(t) f_j(t)}{n\bars_1(t)-1}\, ,
\]
where $n\bars_1(t) := \sum_i f_i(t)$. The above is similar to \chh{the} multiplicative coalescent dynamics, but where \chh{the} proxy for \chh{the size of a component} is the number of alive half-edges.
For an alive half-edge $u$ in a component $\cC(t)$,  write $\mvpi(u) \in [n]$ for the vertex that $u$ is attached to, and \chh{define the} distance between two alive half-edges $u$ and $v$ in the same component $\cC(t)$ as $d(u,v) := d_{\cC(t)}(\mvpi(u),\mvpi(v))$. 
\chh{Let $\diam_{\max}(t)$ be the maximal such distance as we range over all pairs of alive half-edges in the same component at time $t$}.  For a connected component $\cC(t)$, let 
$\cD_1(\cC(t)) := \sum_{u,v \text{ alive half-edges, }u,v \in \cC(t)} d(u,v)$.
Define
\begin{equation}
	\label{eqn:suscep-def-cm}
	\bars_l(t) :=  \sum_i \frac{[f_i(t)]^l}{n},\;\;  
	\barg(t):= \sum_i \frac{f_i(t) |\cC_i(t)|}{n},\ \text{ and } \
	\barcd(t):=  \sum_i \frac{\cD_1(\cC_i(t))}{n}.
\end{equation}
\begin{thm}
	\label{thm:config-bare-sibcrit}
	Fix $\delta \in (1/6,1/5)$, and let $t_n$ be as in \eqref{eqn:tn-def}. Then
	\begin{gather}
		\label{eqn:s2-asymp-cm}
		\Big|\frac{n^{1/3}}{\bars_2(t_n)} - \frac{\nu^2 n^{1/3-\delta}}{\mu(\nu-1)^2}\Big| \weakc 0,\\
		\label{eqn:s3-asymp-cm}
		{\bars_3(t_n)}/{[\bars_2(t_n)]^3} \weakc {\beta}/[{\mu^3(\nu-1)^3}],\\
		\label{eqn:d1-g-asymp-cm}
		{n^{-2\delta}}{\barcd(t_n)} \weakc {\mu(\nu-1)^2}/{\nu^3}, \ \text{ and }\
		{n^{-\delta}}{\barg(t_n)} \weakc {(\nu-1)\mu}/{\nu^2}.
	\end{gather}
\end{thm}


\begin{thm}
	\label{thm:config-largest-comp-diam}
	Fix $\delta\in (0, 1/4)$, and let $t_n$ be as in \eqref{eqn:tn-def}. 
	Then there exists $C=C(\delta)>0$ such that
	\[\pr\Big(|\cC_1(t)| \leq \frac{C (\log{n})^3}{(t_c-t)^2}\, ,\, 
	\diam_{\max}(t) \leq \frac{C (\log{n})^3}{(t_c-t)} \mbox{ for all }  0\leq t \leq t_n \Big) \to 1\, , \text{ as }n\to\infty\, .\]
\end{thm}

\section{Discussion}
\label{sec:disc}

\subsection{Critical regime for random graphs}
\label{sec:disc-crit-sizes}

For an overview of network models see
\cite{bollobas-RG-book, janson2011random,
	durrett2007random,van2009random,chung2006complex}. Specific to {the connectivity
	threshold} see
\cite{Boll-Rior-crit-RG-survey} for a nice recent overview.
For the \erdos random graph there is a huge body of work starting with \cite{erd1959random},
and expanded in great detail in
\cite{luczak1990component,janson1993birth}. Of particular relevance to this work is \cite{aldous1997brownian}. 
We now turn to the models discussed in this paper.

{\it IRG model:} This model was introduced in its general form in \cite{BBSJOR07} where a wide array of properties {were} analyzed. 
In the critical regime, the only results we are aware of are for the \emph{rank-one} model which, in the regime of relevance here, is equivalent to the models introduced in \cite{norros2006conditionally,britton2006generating,chung2002connected}.
At criticality, order of magnitude results for the largest component for the rank-one model were derived in \cite{van2013critical}, and this was strenghtened to distributional convergence in \cite{SBVHJL10,SBVHJL12,turova2009diffusion}.

{\it Configuration model:} Necessary and sufficient conditions for existence of a giant component were derived in \cite{MR98}. The continuous time construction used in this paper is similar to the one considered in \cite{janson2009new}. Component sizes in the critical regime have been studied in \cite{nachmias2010critical,joseph2010component,riordan2012phase}.  Further, critical percolation on random regular graphs was studied in \cite{nachmias2008critical} and \chh{the} diameter was shown to scale like $n^{1/3}$.

\subsection{Differential equations method and average distances}
\label{sec:disc-diff-eqn}

One major tool in dealing with dynamic random graph processes is the differential equations technique,
see, e.g., \cite{kurtz-density} for stochastic process applications, \cite{darling2008differential} \chh{for a nice survey of applications relevant to Markov chains}, and \cite{wormald1995differential} for applications  to random graph processes. 
In this paper they play a key role in the analysis of susceptibility functions in Theorem \ref{thm:config-bare-sibcrit}, in particular the average distance scaling within components in the barely subcritical regime of the configuration model; see also Lemma \ref{lem:diff-eq-solution} and Proposition \ref{prop:suscep-close-det}. We are unaware of similar applications for distance scaling {in the literature}. 


\subsection{Barely subcritical regime}
\label{sec:disc-barely-sub}
Some results in this paper (e.g., Theorems \ref{prop:irg-barely-subcritical} and \ref{thm:config-bare-sibcrit}) deal with susceptibility functions $\bars_2$ and $\bars_3$ at times $t_n = t_c-\eps_n$ where $\eps_n = n^{-\delta}$ with $\delta \in (1/6,1/5)$. 
\chhh{For these models, asymptotics for these functionals were studied previously at time} $t_c-\eps$ with $\eps > 0$ fixed; see, e.g., \cite{janson-suscep-CM} for the configuration model, and \cite{janson-suscep-IRG} for the IRG. 
For the \erdos process, more precise estimates of the susceptibility function in the barely subcritical regime are derived in \cite{janson-suscep-subcrit}. See \cite{aldous2000random} for similar estimates for a random graph model with immigration.

\section{Proofs: Universality}
\label{sec:proofs-universality}
In this section we prove Theorem \ref{thm:aldous-gen-2}.
Sections \ref{sec:def-size-bias}--\ref{sec:proof-tilt-p-trees} collect some preliminary results related to the model.
Section \ref{sec:scaling-rank-one-connected} studies the scaling limit of a related model of a {\it connected} random graph. 
Building on this result, Section \ref{sec:proof-aldous-gen-2} completes the proof.

\subsection{Size-biased random order}\label{sec:def-size-bias}
Given a collection of positive weights $\vx=(x_i, i\in [n])$, a \emph{size-biased random order} of $[n]$ is a random permutation $(v(1), v(2), \ldots, v(n))$ such that $\pr(v(1) =k) \propto x_k$ for $k\in [n]$, and for $2\leq i\leq n$ and $k\in [n]\setminus \set{v(1), \ldots, v(i-1)}$,
\begin{align}
	\pr\big(v(i) =k|v(1), v(2), \ldots, v(i-1) \big) \propto x_k. \label{eqn:size-bias}
\end{align}
One construction of such an ordering is as follows:  
Generate independent exponentials $\xi_i \sim \text{Exp}(x_i)$, $i\in [n]$, and arrange them in increasing order as $\xi_{v(1)} < \xi_{v(2)} < \cdots < \xi_{v(n)}$.  
Then $(v(1), \ldots, v(n))$ would be a size-biased permutation of $[n]$ using the weight sequence $\vx$. In  \cite[Section 3.1]{aldous1997brownian} Aldous used this to construct
$\cG(\vx,q)$ simultaneously with a breadth-first exploration of the graph such that the vertices in the exploration appeared in a size-biased random order. 
We describe this construction next. 
Let $\xi_{i,j} \sim \text{Exp}(qx_j)$,
$i\neq j \in [n]$, be independent random variables. 
The exploration process initializes by selecting a vertex $v(1)$ with 
$\pr(v(1) =k) \propto x_k$ for $k\in [n]$.
The neighbors (sometimes referred to as children) of $v(1)$ are the vertices
$\set{i: \xi_{v(1),i} \leq x_{v(1)}}$. Write $c(1)$ for the number of children of
$v(1)$, labeled as $v(2), v(3), \ldots, v(c(1)+1)$ in
increasing order of the $\xi_{v(1), v(i)}$ values. Now move to $v(2)$ and obtain
the children of $v(2)$ as the vertices 
$i\notin\{ v(1),\ldots, v(c(1)+1)\}$ 
such that $ \xi_{v(2), i}\leq x_{v(2)}$. 
Label them as $v(c(1)+2), \ldots, v(c(1)+c(2)+1)$ in increasing order of their $\xi_{v(2),i}$ values. 
Proceed recursively until the component of $v(1)$ has been explored. 
Then select a new vertex amongst the unexplored vertices with probability proportional to their weights and proceed until all vertices have been explored.

\subsection{Partition into connected components}
\label{sec:part-conn-comp}
{We mainly follow \cite{SBSSXW14}.}
Recall that $(\cC_{i}, i \geq 1)$ denoted the ranked components of $\cG(\vx,q)$, and let 
$\cV_i := V(\cC_i)$. 
Thus, $(\cV_i, i\geq 1)$ is a random partition of the vertex set $[n]$. 
For $\cV\subseteq [n]$ write $\bG_{\cV}^{\con}$ {for} the space of simple connected graphs with vertex set $\cV$. 
Given $a > 0$ and a pmf $\vp = (p_v, v \in \cV)$,
define the probability distribution $\pr_{\con}(\; \cdot \; ; \vp, a, \cV)$ on $\bG_{\cV}^{\con}$ via
\begin{equation}
	\label{eqn:pr-con-vp-a-cV-def}
	\pr_{\con}(G; \vp, a, \cV) \propto \prod_{(u,v)\in E(G)} (1-e^{-a p_u p_v}) \prod_{(u,v)\notin E(G)} e^{-a p_u p_v}\, .
\end{equation}
For $i \geq 1$, define
\begin{equation}
	\label{eqn:vp-i-a-i-def}
	\vp_i =(p_{i,v}, v\in \cV_i) := \big( {x_v}/{\sum_{v \in \cV_i}x_v },\ v \in \cV_i\big)
	\ \text{ and }\
	a_i:= q  \cdot \big( \sum_{v \in \cV_i}x_v \big)^2.
\end{equation}

\begin{prop}[{\cite[Proposition 6.1]{SBSSXW14}}]
	\label{prop:generate-nr-given-partition}
	Let $N$ be the number of components in $\cG(\vx,q)$. Then for any $G_i \in  \bG_{\cV_i}^{\con}$, $1 \leq i \leq N$,
	\begin{equation*}
		\pr\left(\cC_i = G_i, \;\forall 1\leq i\leq N \mid  N, (\cV_i)_{1 \leq i \leq N} \right) = \prod_{1\leq i \leq  N} \pr_{\con}( G_i; \vp_i, a_i, \cV_i).
	\end{equation*}
\end{prop}

\subsection{Tilted $\vp$-trees and scaling limits of components }
\label{sec:proof-tilt-p-trees}
We recall here a result from \cite{SBSSXW14} on scaling limits of graphs with distribution $\pr_{\con} = \pr_{\con} (\; \cdot \; ; \vp, a, [m])$ under regularity assumptions on $a = a^{\sss(m)}> 0$ and the pmfs $\vp= \vp^{\sss(m)}=(p_i^{\sss(m)}, 1\leq i \leq m)$. 
We suppress dependence on $m$ to simplify notation. 
Let $\sigma(\vp) := (\sum_{i\in [m]}p_i^2)^{1/2}$, $p_{\max} := \max_{i \in [m]} p_i$, and $p_{\min} := \min_{ i \in [m]} p_i$. Suppose $\cG^{\vp}\sim \pr_{\con}$.  
Recall the notation $\tilde \ve^{\theta}$ from Section \ref{sec:cont-limit-descp}.
\begin{thm}[{\cite[Theorem 7.3]{SBSSXW14}}]
	\label{thm:sbssxw-conn-res}
	Assume that there exist $\bar \gamma, r_0, \eta_0 \in (0,\infty)$ such that
	\begin{equation}
		\label{eqn:ass-connected-1}
		\sigma(\vp) \to 0, \; \; \frac{p_{\max}}{[\sigma(\vp)]^{3/2+\eta_0}} \to 0, \;\;  \frac{[\sigma(\vp)]^{r_0}}{p_{\min}} \to 0, \;\;
		\text{ and }\, \, 	  a \sigma(\vp) \to \bar \gamma,
	\end{equation}
	as $m\to \infty$. 
	View $\cG^{\vp}$ as a metric measure space by endowing it with the graph distance and the probability measure that assigns mass $p_i$ to the vertex $i$.
	Then
	\begin{equation*}
		\sigma(\vp) \cdot \cG^{\vp} = \scl\left(\sigma(\vp), 1\right) \cdot \cG^{\vp} \weakc \cG( 2 \tilde \ve^{\bar \gamma}, \bar \gamma \tilde \ve^{\bar \gamma}, \cP),
	\end{equation*}
	where the metric measure space $\cG( 2 \tilde \ve^{\bar \gamma}, \bar \gamma \tilde \ve^{\bar \gamma}, \cP)$ is as defined in Section \ref{sec:cont-limit-descp}.
\end{thm}
We recall from \cite{SBSSXW14} an ingredient in the proof of this theorem that will be required in this paper. 
Write $\bT_m^{\ord}$ for the space of ordered rooted trees with vertex set $[m]$. 
Here, ``ordered'' means that the children of every vertex is arranged from ``left'' to ``right.''
(In other words, these are rooted labeled trees with a plane embedding.)
Let $\vt\in \bT_m^{\ord}$ with root $\phi$.
For a vertex $v\neq\phi$, let $(\phi=v_0, v_1, \ldots, v_k=v)$ be the path in $\vt$ from $\phi$ to $v$.
Let $R(v; \vt)$ be the set of all vertices $u$ such that for some $0\leq i\leq k-1$, $u$ is a child of $v_i$ and $u$ appears on the right side of $v_{i+1}$ in the plane embedding.
Define the collection of \emph{permitted edges} to be
\[
\sP(\vt):= \set{\set{v, u}: v\in [m]\setminus\phi, u\in R(v; \vt)}\, . 
\]
For $\vt\in
\bT_m^{\ord}$ and $v\in [m]$, write $d_v(\vt)$ for the number of
children of $v$ in $\vt$. 
Define a probability measure $\pr_{\ord}(\cdot)=\pr_{\ord}(\; \cdot \; ; \vp)$ {on $\bT_m^{\ord}$ via}
\begin{equation}
	\label{eqn:ordered-p-tree-def}
	\pr_{\ord}(\vt) = \prod_{v\in [m]} {p_v^{d_v(\vt)}}/{(d_v(\vt)) !}, \qquad \vt \in \bT_m^{\ord}.
\end{equation}
A random tree with distribution $\pr_{\ord}$ is a $\vp$-tree (\cite{aldous2004exploration,pitman2001random}) together with a plane embedding.
Define $L : \bT_m^{\ord} \to \bR_+$ by
\begin{equation}
	\label{eqn:ltpi-def}
	\displaystyle L(\vt):= \Big(\prod_{\set{i,j}\in E(\vt)} \Big[\frac{\exp(a p_i p_j)- 1}{ap_ip_j} \Big]\Big) \exp\big(\sum_{\set{i,j} \in \sP(\vt)} a p_i p_j\big).
\end{equation}
Define the \emph{tilted} $\vp$-tree distribution $\tilde{\pr}_{\ord}(\cdot) = \tilde \pr_{\ord}(\; \cdot \; ; \vp, a)$ via 
\begin{equation}
	\label{eqn:tilt-p-tree-dist}
	\frac{d\tilde{\pr}_{\ord}}{d\pr_{\ord}}(\vt) = \frac{L(\vt)}{\E_{\ord}[L] }\, , \qquad  \vt \in \bT_m^{\ord},
\end{equation}
where $\E_{\ord}$ is the expectation operator with respect to $\pr_{\ord}$.
\begin{prop}[{\cite[Proposition 7.4]{SBSSXW14}}]\label{prop:gp-const-tilt-surp}
	Fix a pmf $\vp$ on $[m]$ and $a> 0$. Then $\cG^{\vp}\sim \pr_{\con}$ can be constructed as follows: Generate $\cT^{\vp}\sim \tilde \pr_{\ord}$. Conditional on $\cT^{\vp}$, add each permitted edge $\set{u,v}\in\sP(\cT^{\vp})$ independently with probability $1-\exp(-ap_up_v)$, and then forget about the root of $\cT^{\vp}$ and the plane embedding of the tree $\cT^{\vp}$.
\end{prop}

\subsection{Connected random graphs with blobs}\label{sec:scaling-rank-one-connected}
For $m \geq 1$, consider $a = a^{\sss(m)} > 0$, and pmfs $\vp = \vp^{\sss(m)} = (p_i, i\in [m])$. 
Suppose we are given the following:

\begin{enumeratea}
	\item 
	{\bf Blob level superstructure:}  Let $\cG^{\vp}\sim \pr_{\con}(\; \cdot \; ;\vp,a,[m]) $.
	\item 
	{\bf Blobs:} Let $\vM = \vM^{\sss(m)} = \big\{(M_i, d_i, \mu_i) :  i\in [m]\big\}$ be a family of compact metric measure spaces where $\mu_i$, $i\in [m]$, are probability measures. 
	Recall $u_i$ from \eqref{eqn:uik-def}.
	\item 
	{\bf Junction points:} $\vX = (X_{i,j},\ i,j \in [m])$ be a family of independent random variables (that is also independent of $\cG^{\vp}$) such that for each $i,j\in [m]$, $X_{i,j}\in M_i$ has law $\mu_i$.
\end{enumeratea}
Recall the operation $\Gamma$ from Section \ref{sec:inter-blob-distance}.
Define
\begin{equation}\label{eqn:g-barp-Am-def}
	\bar \cG^{\vp} := \Gamma(\cG^{\vp}, \vp, \vM, \vX), \;\; 
	A_m := \sum_{i \in [m]} p_i u_{i},\, \,
	\text{ and }\, \,
	{\diam_{\max}} := \max_{i \in [m]} \diam(M_i).
\end{equation}

\begin{thm}\label{thm:augment-metric-space-connected}
	Assume that \eqref{eqn:ass-connected-1} holds and that for some $\eta_0 \in (0,\infty)$,
	\begin{equation}
		\label{eqn:ass-connected-2}
		\lim_{m \to \infty} {[\sigma(\vp)]^{1/2-\eta_0} {\diam_{\max}}}/{(A_m+1)} = 0.
	\end{equation}
	Then
	$\frac{\sigma(\vp)}{A_m + 1} \cdot \bar \cG^{\vp} \weakc \cG(2 \tilde \ve^{\bar \gamma}, \bar \gamma \tilde \ve^{\bar \gamma}, \cP), \mbox{ as } m \to \infty.$
\end{thm}

\noindent \textbf{Proof: }  
We construct $ \bar \cG^{\vp}$ in the following way:
Consider a coupling $(\cG^{\vp}, \cT^{\vp}) \in \bG_{[m]}^{\con}\times \bT_m^{\ord}$ as in Proposition \ref{prop:gp-const-tilt-surp} with $\cG^{\vp} \sim \pr_{\con}$ and $\cT^{\vp} \sim \tilde \pr_{\ord}$.  
Generate the random junction points $\vX$ independent of $(\cG^{\vp}, \cT^{\vp})$.
Set $ \bar \cG^{\vp} := \Gamma(\cG^{\vp}, \vp, \vM, \vX)$.
Further, let $\bar \cT^{\vp} := \Gamma(\cT^{\vp}, \vp, \vM, \vX)$.
By Theorem \ref{thm:sbssxw-conn-res} we have 
$\sigma(\vp)\cdot \cG^{\vp} \weakc \cG(2 \tilde \ve^{\bar \gamma}, \bar \gamma \tilde \ve^{\bar \gamma}, \cP)$. Thus we only need to prove
\begin{equation}
	\label{eqn:1481}
	d_{\GHP}\Big( \sigma(\vp)\cdot \cG^{\vp}, \frac{\sigma(\vp)}{A_m + 1}\cdot \bar \cG^{\vp}  \Big) \weakc 0, \mbox{ as } m \to \infty.
\end{equation}
It would be convenient to work with slight variants of the original spaces $\cG^{\vp}$ and $\bar \cG^{\vp}$. 
Write $\cG_*^{\vp}$ (resp. $\bar \cG_*^{\vp}$) for the metric measure space obtained from $\cT^{\vp}$ (resp. $\bar \cT^{\vp}$) by identifying $i$ and $j$ (resp. $X_{i,j}$ and $X_{j,i}$) for each edge $\{i,j\} \in E(\cG^{\vp}) \setminus E(\cT^{\vp})$, instead of placing an edge of length one between them. 
Recall that $\spls(\cG^{\vp})$ denotes the number of surplus edges in $\cG^{\vp}$. Then one can check that
\begin{equation}
	\label{eqn:1600}
	d_{\GHP}( \cG_*^{\vp}, \cG^{\vp}) \leq \spls(\cG^{\vp}),\ \text{ and }\
	d_{\GHP}( \bar \cG_*^{\vp}, \bar \cG^{\vp}) \leq \spls(\cG^{\vp}).
\end{equation}
View $\cT^{\vp}$ as a metric measure space by endowing it with the tree distance and the probability measure that assigns mass $p_i$ to the vertex $i$.
Let $C_m\subseteq \cT^{\vp}\times \bar \cT^{\vp}$ be given by $C_m := \{ (i,x): i \in [m], x \in M_i\}$.
Let $\nu_m$ be the measure on $\cT^{\vp}\times \bar \cT^{\vp}$ given by
$\nu_m( \{i\} \times A) := p_i \mu_i(A \cap M_i)$ for $i\in [m]$, $A \subseteq \bar \cT^{\vp}$.
Write 
\[
\dsc(\nu_m)= \dsc(\nu_m; \sigma(\vp)\cT^{\vp},  \frac{\sigma(\vp)}{A_m+1}\bar \cT^{\vp})
\ \ \text{ and }\ \ 
\dis(C_m)=\dis\big(C_m; \sigma(\vp)\cT^{\vp}, \frac{\sigma(\vp)}{A_m+1}\bar \cT^{\vp}\big)\, .
\]
By \cite[Lemma 4.2]{addario2013scaling}, we have
\begin{equation}
	\label{eqn:1503}
	d_{\GHP}\Big( \sigma(\vp)\cG_*^{\vp},  \frac{\sigma(\vp)}{A_m+1} \bar \cG_*^{\vp}\Big)
	\leq
	(\spls(\cG^{\vp})+1) \max\Big\{ \frac{\dis(C_m)}{2}, \dsc(\nu_m), \nu_m(C_m^c) \Big\}\, .
\end{equation}
It is easy to check that for all $m$, $\dsc(\nu_m)=0$ and $\nu_m(C_m^c)=0$. Therefore by \eqref{eqn:1503} and \eqref{eqn:1600}, to obtain \eqref{eqn:1481} we only need to show that
\begin{align}\label{eqn:1145-1}
	(\spls(\cG^{\vp}), m\geq 1) \mbox{ is tight, and }
	\dis(C_m) \weakc 0 \mbox{ as } m \to \infty. 
\end{align}
We will prove \eqref{eqn:1145-1} in Lemma \ref{lem:1542} and Lemma \ref{lem:1567} below, which will complete the proof of Theorem \ref{thm:augment-metric-space-connected}. 
\qed

\begin{lem}
	\label{lem:1542}
	The sequence $\big(\spls(\cG^{\vp});\, m\geq 1\big)$ is tight.
\end{lem}
\noindent \textbf{Proof:}
By \cite[Corollary 7.13]{SBSSXW14}, there exists a constant $K_1 > 0$ such that for $m$ large,
\begin{equation}
	\label{eqn:1725}
	\E_{\ord}[L^2] \leq K_1.
\end{equation}
Since $\cG^{\vp}$ is obtained by adding each permitted edge $\set{u,v}$ of $\cT^{\vp}$ independently with probability $1-\exp(-ap_u p_v)$,
$\E [ \spls (\cG^{\vp})] \leq \E\big[ \sum_{\{i,j\} \in \sP(\cT^{\vp})} a p_i p_j\big]$.
Using \eqref{eqn:tilt-p-tree-dist} and \eqref{eqn:1725} together with the relations $L(\vt) \geq 1$ and
$ \sum_{\{i,j\} \in \sP(\vt)} a p_i p_j \leq L(\vt)$,
we have, for all large $m$,
$	\E [ \spls (\cG^{\vp})] \leq \E[L(\cT^{\vp})] =\E_{\ord}[L^2]/\E_{\ord}[L]\leq K_1.$
Thus $\sup_{m \geq 1}\E[\spls(\cG^{\vp})] < \infty$, which implies the desired tightness. \qed

\begin{lem}\label{lem:1567}
	As $m \to \infty$, $\dis(C_m) \weakc 0$.
\end{lem}
\noindent \textbf{Proof:} 
Let $L(\cdot)$ be as in \eqref{eqn:ltpi-def}. 
Suppose $(\Omega,\cF,\bP)$ is a probability space where $\cT^{\vp}$ and $\bar \cT^{\vp}$ are defined.
Let $\bP'$ be the probability measure on $\Omega$ given by
\begin{equation}
	\label{eqn:1537}
	\frac{d \pr}{ d \pr'}(\omega) = \frac{L(\cT^{\vp}(\omega))}{\E_{\ord}[L]} \mbox{ for } \omega \in \Omega.
\end{equation}
By \eqref{eqn:tilt-p-tree-dist}, $\cT^{\vp} \sim \pr_{\ord}$ under $\pr'$.
Suppose we show that for any $\eps>0$,
\begin{align}\label{eqn:9}
	\pr'(\dis(C_m) > \eps)\to 0\ \ \text{ as }\ \ m\to\infty\, .
\end{align}
Then 
$
\big(\pr(\dis(C_m) > \eps)\big)^2
\leq
\pr'(\dis(C_m) > \eps)\cdot\bE_{\ord}[L^2]
\to 0
$
as $m\to\infty$, 
where the first step uses \eqref{eqn:1537}, the fact that $L(\vt)\geq 1$, and the Cauchy-Schwarz inequality, and the second step uses \eqref{eqn:1725} and \eqref{eqn:9}.
So we just have to prove \eqref{eqn:9}.

We introduce some additional random variables: Let $J$ and $Y_1,\ldots, Y_m$ be respectively $[m]$ and $M_1,\ldots, M_m$ valued independent random variables that are also independent of $\cT^{\vp}$ such that $J \sim \vp$, and $Y_i \sim \mu_i$ for $i \in [m]$. We extend the probability space to incorporate these random variables and still write $\pr'$ for the underlying probability measure.

For $ x \in \bar \cT^{\vp}$, write $i(x)$ for the $i \in [m]$ with $x \in M_{i}$. Write $d_{\cT}$ and $d_{\bar \cT}$ for the metrics on $\cT^{\vp}$ and $ \bar \cT^{\vp}$ respectively.   
Then
\begin{equation}
	\label{eqn:1157}
	\dis(C_m) = \sup_{x,y \in \bar \cT^{\vp}}\Big\{ \Big|\sigma(\vp) d_{\cT}(i(x),i(y)) -\frac{\sigma(\vp)}{A_m+1} d_{\bar \cT}(x,y)\Big|\Big\}.
\end{equation}
We first show for any fixed $i_0 \in [m]$ and $x_0 \in M_{i_0}$,
\begin{equation}\label{eqn:1142}
	\dis(C_m) 
	\leq 
	4 \sup_{y \in \bar \cT^{\vp}}
	\Big\{ 
	\Big|\sigma(\vp) d_{\cT}(i_{0},i(y)) - \frac{\sigma(\vp)}{A_m+1} d_{\bar \cT}(x_0,y)\Big|
	\Big\} 
	+ 
	\frac{2\sigma(\vp)}{A_m+1} {\diam_{\max}} \, .
\end{equation}
Indeed, for any two points $x,y \in \bar \cT^{\vp}$, there are two unique paths between $i_0$ and $i(x)$ and between $i_0$ and $i(y)$ in $\cT^{\vp}$. Write $(i_0, \ldots , i_{k})$ for the longest common path shared by these two paths. 
Let $i_* = i_{k}$ and $x_* = X_{i_k, i_{k-1}}$. Since $\cT^{\vp}$ is a tree,
\begin{equation}
	\label{eqn:1146}
	d_{\cT}(i(x), i(y)) = d_{\cT}(i_0, i(x)) + d_{\cT}(i_0, i(y)) - 2d_{\cT}(i_0, i_*).
\end{equation}
By a similar observation for $\bar \cT^{\vp}$, we have
\begin{equation}
	\label{eqn:1151}
	d_{\bar \cT}(x , y) \leq d_{\bar \cT}(x_0, x) + d_{\bar \cT}(x_0, y) - 2d_{\bar \cT}(x_0, x_*) \le	d_{\bar \cT}(x , y) + 2 {\diam_{\max}} .
\end{equation}
Hence, \eqref{eqn:1142} follows using \eqref{eqn:1146} and \eqref{eqn:1151} in \eqref{eqn:1157}.

Next, replace $(i_0, x_0)$ by $(I, Y_I)$, where $I\in [m]$ is the root of $\cT^{\vp}$, and replace every $y \in M_j$ in \eqref{eqn:1142} with $Y_j \in M_j$, which incurs an error of at most ${4\sigma(\vp) {\diam_{\max}}}/{(A_m+1)}$. Thus,
\begin{equation*}
	\dis(C_m) \leq 4 \sup_{j \in [m]} \Big\{ \Big|\sigma(\vp) d_{\cT}(I,j) - \frac{\sigma(\vp)}{A_m+1} d_{\bar \cT}(Y_I,Y_j)\Big|\Big\} + \frac{6\sigma(\vp) {\diam_{\max}}}{A_m+1}.
\end{equation*}
Using \eqref{eqn:ass-connected-2} \ch{and $\sigma(\vp)\to 0$} we can find $m_0$ such that ${6\sigma(\vp) {\diam_{\max}}}/{(A_m+1)} < \eps/5$ for $m > m_0$. Then for $m>m_0$,
\begin{align}
	\pr'( \dis(C_m) > \eps )
	\leq& \pr' \Big( \sup_{j \in [m]} { \Big|\sigma(\vp) d_{\cT}(I,j) - \frac{\sigma(\vp)}{A_m+1} d_{\bar \cT}(Y_I,Y_j)\Big|} > \frac{\eps}{5} \Big) \nonumber \\
	\leq& \frac{1}{p_{\min}}\sum_{j \in [m]} p_j \pr'  \Big( \Big|\sigma(\vp) d_{\cT}(I,j) - \frac{\sigma(\vp)}{A_m+1} d_{\bar \cT}(Y_I,Y_j)\Big| > \frac{\eps}{5} \Big) \nonumber \\
	=& \frac{1}{p_{\min}}\pr' \Big( \Big|\sigma(\vp) d_{\cT}(I,J) - \frac{\sigma(\vp)}{A_m+1} d_{\bar \cT}(Y_I,Y_J)\Big| > \frac{\eps}{5} \Big). \label{eqn:1173}
\end{align}
Write $R^* := d_{\cT}(I, J)+1$,  and let  $(I_0=I, I_1, ..., I_{R^*-1} = J)$ be the path between $I$ and $J$ in $\cT^{\vp}$. Define $\xi_0^* := d_I(Y_I, X_{I,I_1})$, $\xi_{R^*-1}^* := d_J(X_{J,I_{R^*-2}}, Y_J)$, and $\xi_i^* := d_{I_i}(X_{I_i,I_{i-1}}, X_{I_i, I_{i+1}})$ for $1\leq i \leq R^*-2$.
In this notation $d_{\bar \cT}(Y_I,Y_J) = \sum_{i=0}^{R^*-1} \xi_i^* + (R^* - 1)$ and $d_{\cT}(I,J) = R^*-1$. Thus,
\ch{\begin{equation}
		\label{eqn:1590}
		\sigma(\vp)d_{\cT}(I,J) - \frac{\sigma(\vp)}{A_m+1} d_{\bar \cT}(Y_I,Y_J)  = \frac{\sigma(\vp)}{A_m+1} \Big( \Big[\sum_{i=0}^{R^*-1} (A_m - \xi_i^*)\Big] - A_m \Big).
\end{equation}}
The summation above involves the path in a $\vp$-tree connecting the root and a random vertex selected according to the distribution $\vp$. This admits the following alternate construction: 
Let $\vJ:=(J_i, i\geq 0)$ be a sequence of iid $[m]$-valued random variables with law $\vp$.  For each $j \in [m]$, let $\mvxi^{\sss(j)}:=(\xi_i^{\sss(j)}, i\geq 0)$ be an \chh{iid} sequence with $\xi_1^{\sss(j)}\equald d_j(X_{j,1}, X_{j,2})$ such that $\mvxi^{\sss(1)},\ldots, \mvxi^{\sss(m)}$ and $\vJ$ are jointly independent. Let $(\Omega'', \cF'', \pr'')$ be the probability space on which the random variables $\vJ$ and $\big(\mvxi^{\sss(j)}; j\in [m] \big)$ are defined. 
Let $R:= \inf\, \{ k \geq 1: J_k=J_i \text{ for some } i<k \}$ be the first repeat time of the sequence $\vJ$.  By \cite[Corollary 3]{camarri2000limit}, 
\begin{equation*}
	(I_0, ..., I_{R^*-1}; R^*)_{\pr'} \stackrel{d}{=} (J_0, ..., J_{R-1}; R)_{\pr''}.
\end{equation*}
Owing to the independence structure,  \ch{we further} have
\begin{equation*}
	(I_0, ..., I_{R^*-1}; R^*; \xi_0^*, ...,\xi_{R^*-1}^*)_{\pr'} \stackrel{d}{=} (J_0, ..., J_{R-1}; R; \xi_0^{\sss(J_0)}, ..., \xi_{R-1}^{\sss(J_{R-1})} )_{\pr''}.
\end{equation*}
For $i\geq 0$ write $\Delta_i := \xi_i^{\sss(J_i)} - A_m$. Using \eqref{eqn:1590} and the fact that $\sigma(\vp) \to 0$, for large $m$,
\begin{equation}
	\label{eqn:1223}
	\pr'\Big(\Big |\sigma(\vp) d_{\cT}(I,J) - \frac{\sigma(\vp)}{A_m+1} d_{\bar \cT}(Y_I,Y_J)\Big| > \frac{\eps}{5} \Big) \leq \pr''\Big( \frac{\sigma(\vp)}{A_m+1} \Big| \sum_{i=0}^{R-1} \Delta_i \Big| > \frac{\eps}{6} \Big).
\end{equation}
Note that $\big(\xi_i^{\sss (J_i)}; i\geq {0}\big)$ is a collection of iid random variables with mean
$\sum_{k\in [m]} p_k \E[\xi_0^{\sss(k)}] = A_m$.
Thus $(\sum_{i=0}^k \Delta_i,\ k\geq 0)$ is a martingale with respect to \ch{its} natural filtration.
Then, for any $t > 0$,
\begin{equation}
	\label{eqn:1199}
	\pr''\Big(\frac{\sigma(\vp)}{A_m+1} \Big| \sum_{i=0}^{R-1} \Delta_i \Big| > \frac{\eps}{6}\Big) 
	\leq 
	\pr''(R \geq t) + \pr''\Big(\sup_{0\leq k \leq t-1}  \Big|  \sum_{i=0}^{k} \Delta_i \Big| > \frac{\eps (A_m+1)}{6 \sigma(\vp)}\Big).
\end{equation}
For the first term, by \ch{\cite[Lemma 10.7]{SBSSXW14}} (also see \cite[Equations (26), (29)]{camarri2000limit}),
\begin{equation}
	\label{eqn:1794}
	\pr''\left(R \geq t\right)  \leq 2 \exp \big( -t^2 \sigma^2(\vp)/24\big)\ \mbox{ for }\ t \in (0, 1/p_{\max}).
\end{equation}
The second term on the right side of \eqref{eqn:1199} can be bounded by using \chh{Markov's} inequality and the Burkholder-Davis-Gundy (BDG) inequality \cite{burkholder-davis-gundy}.
For fixed $r \geq 1$, we have
\begin{align}
	\pr''\Big(\sup_{0\leq k \leq t-1}  \big|  \sum_{i=0}^{k} \Delta_i \big| > \frac{\eps (A_m+1)}{6 \sigma(\vp)}\Big)
	\leq
	\Big(\frac{6 \sigma(\vp)}{\eps (A_m+1)}\Big)^{2r}  \E''\Big[ \sup_{0\leq k \leq t-1}  \big|  \sum_{i=0}^{k} \Delta_i \big|^{2r} \Big]& \nonumber\\
	\hskip10pt
	\leq \Big(\frac{6 \sigma(\vp)}{\eps (A_m+1)}\Big)^{2r} K_2(r) \E''\Big[ \big(\sum_{i=0}^{t-1} \Delta_i^2\big)^{r} \Big]
	\leq  \Big(\frac{6 \sigma(\vp)}{\eps (A_m+1)}\Big)^{2r} K_2(r)\cdot t^r {\diam_{\max}^{2r}}, \label{eqn:1240}&
\end{align}
where $\bE''$ denotes expectation with respect to $\bP''$, the constant $K_2(r)>0$ comes from the BDG inequality, and the last step uses $|\Delta_i| \leq {\diam_{\max}}$. Define $\alpha_m := -24 r \log \sigma(\vp)$ and $t_m := \alpha_m/\sigma(\vp)$. Combining \eqref{eqn:1199}, \eqref{eqn:1240}, \eqref{eqn:1794} with $t = t_m $, we have \ch{for sufficiently large $m$,}
\begin{align}
	\pr''\Big( \frac{\sigma(\vp)}{A_m+1} \big| \sum_{i=0}^{R-1} \Delta_i \big| > \frac{\eps}{6} \Big)
	\leq& 2 \exp \Big( -\frac{t_m^2 \sigma^2(\vp)}{24}\Big) + K_3(r,\eps)\cdot \frac{\alpha_m^r \sigma^r(\vp) {\diam_{\max}}^{2r}}{(A_m+1)^{2r}} \nonumber \\
	\leq& 2 [\sigma(\vp)]^r + K_3(r,\eps) \alpha_m^r [\sigma(\vp)]^{2r\eta_0}\, , \label{eqn:1244}
\end{align}
where $K_3(r,\eps) = 6^{2r} K_2(r)/\eps^{2r}$, and the last line uses  \eqref{eqn:ass-connected-2}. 
(It is easy to check that $\lim\limits_{m \to\infty} t_m p_{\max} = 0$. 
Thus, we can apply \eqref{eqn:1794} for $m$ large, and \eqref{eqn:1244} is valid for large $m$.) 
By \eqref{eqn:1244}, \eqref{eqn:1223}, and \eqref{eqn:1173}, \chh{taking $r =r_{\ast}:= \max\set{r_0, \lceil (r_0+1)/(2\eta_0) \rceil}$, 
	\begin{equation*}
		\pr'(\dis(C_m) > \eps)
		\leq 
		\Big(2 [\sigma(\vp)]^{r_0} + K_3(r_{\ast},\eps) \alpha_m^{r_{\ast}} [\sigma(\vp)]^{r_0+1}\Big)/p_{\min}.
	\end{equation*}
}
By \eqref{eqn:ass-connected-1}, the above expression goes to zero as $m \to \infty$.
This completes the proof of \eqref{eqn:9}. 
\qed

\subsection{Proof of Theorem \ref{thm:aldous-gen-2}}
\label{sec:proof-aldous-gen-2}
The previous section deals with asymptotics for a single connected component. This section leverages the above results to study maximal components of $\bar\cG(\vx,q,\vM)$ and prove Theorem \ref{thm:aldous-gen-2}.
Recall that $\vx$ is the weight sequence, $\vM$ is the collection of blobs, $\vX$ is the collection of random junction points, and $u_1,\ldots, u_n$ are the average distances.
Also recall that $\sigma_r = \sum_{i \in [n]} x_i^r$, $r \geq  1$, and $\tau = \sum_{i \in [n]} x_i^2 u_i$.  We start with the following auxiliary result.

\begin{prop}\label{prop:moments-convergence-each-component}
	Under Conditions \ref{ass:aldous-basic-assumption} and \ref{ass:aldous-gen-2}, 
	for each $i\geq 1$,
	\begin{equation*}
		\frac{\sum_{v\in \cC_i} x_v^2}{\sum_{v\in \cC_i} x_v } \cdot \frac{\sigma_2}{\sigma_3} 
		\weakc 1,
		\ \ \text{ and }\ \ 
		\frac{\sum_{v\in \cC_i}  x_v  u_{v}}{\sum_{v\in \cC_i}  x_v} \cdot \frac{\sigma_2}{\tau}  
		\weakc 1\, , 
		\ \ \text{ as }\ \  n \to \infty.
	\end{equation*}
\end{prop}
We will make use of the following lemma in the proof.
\begin{lem}[\cite{SBSSXW14}, Lemma 8.2]\label{lem:size-biased-partial-sum}
	For $n\geq 1$, fix $\vx := (x_i > 0,\ i \in [n])$ and $\vu  := (u_i \geq 0,\ i \in [n])$ such that $c_n:= \sum_{i \in [n]} x_i u_i/{\sum_{i\in [n]} x_i} > 0$.  Let $(v(i), i \in [n])$ be a size-biased random \ch{ordering} of $[n]$ using weights $\vx$.
	Let $x_{\max}:= \max_{i \in [n]} x_i$ and $u_{\max} := \max_{i \in [n]}u_i$. Let $\ell= \ell(n)  \in [n]$ be such that
	\begin{equation}
		\label{eqn:ass-size-biased-partial-sum}
		{\ell x_{\max}}/{\sigma_1} \to 0\, ,\ \text{ and }\ \   {u_{\max}}/( \ell c_n) \to 0, \qquad \text{as }n\to\infty.
	\end{equation}
	Then
	\begin{equation*}
		\sup_{k \leq \ell }\Big| \frac{\sum_{i = 1}^k u_{v(i)}}{\ell c_n} - \frac{k}{\ell} \Big| \weakc 0\, ,\  \text{ as }n\to\infty.
	\end{equation*}
\end{lem}
We are now ready to prove Proposition \ref{prop:moments-convergence-each-component}.

\noindent\textbf{Proof of Proposition \ref{prop:moments-convergence-each-component}: }
We only work with $i=1$ to keep notation simple.
Recall the breadth-first construction of $\cG(\vx,q)$ from Section \ref{sec:def-size-bias}, where vertices are explored in a size-biased order $(v(i), 1\leq i\leq n)$ using the weight sequence $\vx$. 
The following properties of this exploration were derived in \cite{aldous1997brownian}:
\begin{inparaenuma}
	\item There exist random variables $m_{L}, m_{R} \in [n]$ such that $V(\cC_1)=\set{ v(i) : m_L +1 \leq i \leq m_R} $.
	\item Under \chh{Condition} \ref{ass:aldous-basic-assumption}, $\sum_{i=1}^{m_{R}} x_{v(i)}$ is tight.
	\item By Theorem \ref{thm:aldous-review},  as $n \to \infty$, $\sum_{i=m_L+1}^{m_R} x_{v(i)} \weakc \chh{\gamma_1(\lambda)}$.
\end{inparaenuma}
In this notation, the claim in Proposition \ref{prop:moments-convergence-each-component}, for $i=1$, is equivalent to
\begin{equation}
	\label{eqn:prop67-equiv}
	\frac{\sum_{i=m_L+1}^{m_R} x_{v(i)}^2}{\sum_{i=m_L+1}^{m_R} x_{v(i)}}\cdot \frac{\sigma_2}{\sigma_3} \weakc 1, 
	\ \ \text{ and }\ \
	\frac{\sum_{i=m_L+1}^{m_R} x_{v(i)}  u_{v(i)}}{\sum_{i=m_L+1}^{m_R} x_{v(i)}}\cdot \frac{\sigma_2}{\tau} \weakc 1.
\end{equation}
Thus, it is enough to prove \eqref{eqn:prop67-equiv}.
We first show that
\begin{equation} \label{eqn:1429}
	\frac{\sum_{i =m_{L}+1}^{m_R}  x_{v(i)} }{m_R - m_L} \cdot \frac{\sigma_1}{\sigma_2} \weakc 1, \mbox{ as } n \to \infty.
\end{equation}
Fix $\eta >0$. Since $\big(\sum_{i=1}^{m_{R}} x_{v(i)};\, n\geq 1\big)$ is a tight sequence, there exists $T >0$ such that for all $n\geq 1 $,
$	\pr\left(\sum_{i=1}^{m_{R}} x_{v(i)} \geq T\right) < \eta$.
Let $m_0 := \sigma_1/\sigma_2$ and apply Lemma \ref{lem:size-biased-partial-sum} with $\ell = 2T m_0$ and $u_i = x_i$ (thus  $c_n = \sigma_2/\sigma_1$). Condition \eqref{eqn:ass-size-biased-partial-sum} is equivalent to $ x_{\max}/\sigma_2 \to 0$, which follows from Assumption \ref{ass:aldous-basic-assumption}. By Lemma \ref{lem:size-biased-partial-sum}, there exists $N_\eta > 0$ such that for all $n > N_\eta$,
\begin{equation}\label{eqn:1438}
	\pr\Big( \sup_{k \leq 2Tm_0}\big|\sum_{i=1}^k x_{v(i)} - \frac{k}{m_0}\big| > \eta\Big) < \eta.
\end{equation}
On the event
$\big\{\sup_{k \leq 2Tm_0}\left|\sum_{i=1}^k x_{v(i)} - \frac{k}{m_0}\right| \leq \eta\big\} 
\cap 
\big\{\sum_{i=1}^{m_{R}} x_{v(i)} \leq T\big\}$,
we must have
$m_L < m_R < 2T m_0$ (assuming $\eta < T$), and hence $|\sum_{i=m_L+1}^{m_R} x_{v(i)} - {(m_R - m_L)}/{m_0}| < 2 \eta$. 
Since $\eta$ can be arbitrarily small, using \eqref{eqn:1438} we have
\begin{equation}\label{eqn:47}
	\Big|\sum_{i=m_L+1}^{m_R} x_{v(i)} - \frac{m_R - m_L}{m_0}  \Big|\weakc 0 \mbox{ as } n \to \infty.
\end{equation}
By property (c) of the exploration and \eqref{eqn:47},  
${(m_R - m_L)}/{m_0} \weakc \gamma_1(\lambda)$,
which coupled with \eqref{eqn:47} yields \eqref{eqn:1429}.

Now, under Conditions \ref{ass:aldous-basic-assumption} and \ref{ass:aldous-gen-2},
\begin{equation}
	\label{eqn:we-want-moment}
	\frac{\sigma_2 x_{\max}^2}{\sigma_3} \to 0\, ,\ \ 
	\text{ and }\ \  \
	\frac{\sigma_2 x_{\max} {\diam_{\max}}}{\tau} \to 0\, ,\ \  \mbox{ as }\ n \to \infty.
\end{equation}
Thus, repeating the above argument with respectively $ x_{v(i)}^2$  and $ x_{v(i)} u_{v(i)}$ in place of $x_{v(i)}$ and using Lemma \ref{lem:size-biased-partial-sum}, we get
\begin{equation}\label{eqn:1447}
	\frac{\sum_{i =m_{L+1}}^{m_R}  x_{v(i)}^2 }{m_R - m_L} \cdot \frac{\sigma_1}{\sigma_3} \weakc 1\, ,\ \ 
	\text{ and }\ \  \
	\frac{\sum_{i =m_{L+1}}^{m_R}  x_{v(i)}  u_{v(i)} }{m_R - m_L} \cdot \frac{\sigma_1}{\tau} \weakc 1.
\end{equation}
Combining \eqref{eqn:1429} and \eqref{eqn:1447} completes the proof of \eqref{eqn:prop67-equiv}, and thus of Proposition \ref{prop:moments-convergence-each-component}. \qed\\

\noindent \textbf{Proof of Theorem \ref{thm:aldous-gen-2}:}
Since we work with the product GHP topology, it is enough to prove the assertion for the maximal component $\bar \cC_1$, i.e.,
\begin{equation}
	\label{eqn:conv-first-component}
	\scl\Big(\frac{\sigma_2^2}{\sigma_2 + \tau}, 1 \Big)\bar \cC_1 \weakc \cG(2 \tilde \ve_{\gamma_1}, \tilde \ve_{\gamma_1}, \cP_1 ), \mbox{ as } n \to \infty.
\end{equation}
In \eqref{eqn:conv-first-component} and the proof below, to simplify notation, we suppress dependence on $\lambda$ and write $\gamma_1$ for $\gamma_1(\lambda)$. 
By Theorem \ref{thm:aldous-review} and Proposition \ref{prop:moments-convergence-each-component},  without loss of generality, we \ch{may} consider \chh{a} probability space on which the following \chh{convergences} hold almost surely: 
\begin{equation}
	\mass(\cC_1) \convas \gamma_1,   \ \
	\frac{\sum_{v\in \cC_1} x_v^2}{\mass(\cC_1)} \cdot \frac{\sigma_2}{\sigma_3} \convas 1, \ \ 
	\text{ and }\ \ 
	\frac{\sum_{v\in \cC_1}  x_v  u_{v}}{\mass(\cC_1)} \cdot \frac{\sigma_2}{\tau}  \convas 1. \label{eqn:1616-1109}
\end{equation}
Consider the following construction of $\bar \cC_1$: 
Define $\vp := \big( {x_v}/{\mass(\cC_1)}, v \in \cC_1\big)$ and 
$a := q \cdot [\mass(\cC_1)]^2$.
Conditioned on $\cV_1 := V(\cC_1)$, let $\cC_1'$ be a $\bG^{\con}_{\cV_1}$-valued random variable with distribution {$\bP_{\con}(\cdot; \vp, a, \cV_1)$}.
By Proposition \ref{prop:generate-nr-given-partition}, 
\begin{align}\label{eqn:94}
	\bar \cC_1 {\equald} \scl(1, \mass(\cC_1)) \Gamma(\cC_1', \vp, \vM, \vX)\, .
\end{align} 
In order to apply Theorem \ref{thm:augment-metric-space-connected}, we will verify the statements in \eqref{eqn:ass-connected-1} and \eqref{eqn:ass-connected-2} with $\bar \gamma = \gamma_1^{3/2}$. Note that $p_{\max} \leq {x_{\max}}/{\mass(\cC_1)}$, $p_{\min} \geq {x_{\min}}/{\mass(\cC_1)}$, and by \chh{\eqref{eqn:1616-1109}},
\begin{equation}
	\label{eqn:1866}
	\sigma(\vp) = \frac{\sqrt{\sum_{v \in \cC_1} x_v^2}}{\mass(\cC_1)} \sim \frac{\sigma_2}{\gamma_1^{1/2}}
	\, ,
	\ \ \text{ and }\ \
	A_m = \frac{\sum_{v\in \cC_1}  x_v  u_{v}}{\mass(\cC_1)} \sim \frac{\tau}{\sigma_2}\, ,
\end{equation}
where the simplification uses the relation $\sigma_3 \sim \sigma_2^3$ (by Condition \ref{ass:aldous-basic-assumption}). 
In addition, {using} $q \sim 1/\sigma_2$ from Condition \ref{ass:aldous-basic-assumption}, $a \cdot \sigma(\vp) \sim q\gamma_1^2 \cdot \frac{\sigma_2}{\gamma_1^{1/2}} \sim {\gamma_1^{3/2}}$. The remaining conditions in \eqref{eqn:ass-connected-1} and \eqref{eqn:ass-connected-2} follow directly from Conditions \ref{ass:aldous-basic-assumption} and \ref{ass:aldous-gen-2}, and the details are omitted.

Applying Theorem \ref{thm:augment-metric-space-connected} and using \eqref{eqn:94}, we see that
\begin{equation*}
	\scl\Big( \frac{\sigma_2^2}{(\sigma_2 + \tau)\gamma_1^{1/2}}, \frac{1}{\gamma_1} \Big) \bar \cC_1
	\weakc \cG(2 \tilde \ve^{\gamma_1^{3/2}}, \gamma_1^{3/2} \tilde \ve^{\gamma_1^{3/2}}, \cP_1).
\end{equation*}
Note that for any excursions $h$ and $g$ and Poisson point process $\cP$, for $\alpha, \beta >0$, we have
$	\scl(\alpha, \beta) \cG(h,g,\cP) \stackrel{d}{=} \cG\big( \alpha h( \cdot/\beta), \frac{1}{\beta}g(\cdot/\beta), \cP\big)$.
Thus,
\begin{equation*}
	\scl(\gamma_1^{1/2}, \gamma_1) \scl\Big( \frac{\sigma_2^2}{(\sigma_2 + \tau)\gamma_1^{1/2}}, \frac{1}{\gamma_1} \Big) \bar \cC_1  \weakc \cG\Big(2 \gamma_1^{1/2}\tilde \ve^{\gamma_1^{3/2}}(\cdot/ \gamma_1), \gamma_1^{1/2} \tilde \ve^{\gamma_1^{3/2}}(\cdot/\gamma_1), \cP_1\Big).
\end{equation*}
By \cite[Display (4.9)]{SBSSXW14}, we have $\gamma_1^{1/2} \tilde \ve^{\gamma_1^{3/2}}(\cdot/\gamma_1) \stackrel{d}{=} \tilde \ve_{\gamma_1}(\cdot)$. 
Combining this observation with the last display completes the proof of \eqref{eqn:conv-first-component}. 
\qed


\section{Proofs: Scaling limits of IRG}
\label{sec:proof-irg}
In Section \ref{sec:irg-res} we introduced  $\cG_{\IRG}^{(n),\star}$. 
This graph is closely related to $\cG_{\IRG}^{\sss(n), -}$, but it is technically easier to prove that \chh{its} components satisfy \chh{the conditions} (Theorem \ref{prop:irg-barely-subcritical}) required to apply Theorem \ref{thm:aldous-gen-2}. In Section \ref{sec:scaling-irg}, assuming Theorem \ref{prop:irg-barely-subcritical}, we complete the proof of all the other results.
Section \ref{sec:irg-branching-process} builds technical  machinery connected to multitype branching \chh{processes}. \chh{This is then} used in Section \ref{sec:irg-final-section} to prove Theorem \ref{prop:irg-barely-subcritical}.

\subsection{Scaling limit for the IRG model}\label{sec:scaling-irg}

In this section we prove Corollary \ref{cor:irg-barely-subcritical-original}, Theorem \ref{thm:scaling-limit-irg}, and Theorem \ref{thm:component-limit-irg} assuming Theorem \ref{prop:irg-barely-subcritical}.
Recall the kernel $\kappa_n^{-}$ from \eqref{eqn:def-kappa^-}.

\vskip3pt

\noindent{\bf Proof of Corollary \ref{cor:irg-barely-subcritical-original}:} Define $p_{ij}^-=p_{ij}^{\sss(n), -}=1-\exp(-\kappa_n^-(x_i, x_j)/n)${, and recall from Section \ref{sec:irg-res} that $p_{ij}^{\star}:= 1\wedge \big(\kappa_n^{-}(x_i, x_j)/n\big)$.}
Note that $\cG_{\IRG}^{\sss(n),-}$ and $\cG_{\IRG}^{\sss(n), \star}$ are both models of random graphs where we place edges independently between different vertices $i,j\in [n]$, using {connection probabilities} $p_{ij}^{-}$ for $\cG_{\IRG}^{\sss(n),-}$ and $p_{ij}^{\star}$ for $\cG_{\IRG}^{\sss(n), \star}$.
Now,
\begin{equation}
	\label{eqn:pij-pij-star-small}
	\sum_{1\leq i<j\leq n}{\big(p_{ij}^{-}-p_{ij}^{\star}\big)^2}/{p_{ij}^\star}\leq
	\sum_{1\leq i<j\leq n}{p_{ij}^\star}^3=O\big({n}^{-1}\big).
\end{equation}
The claim follows from Theorem \ref{prop:irg-barely-subcritical} and the asymptotic equivalence between  $\cG_{\IRG}^{\sss(n),-}$ and $\cG_{\IRG}^{\sss(n), \star}$ \ch{under \eqref{eqn:pij-pij-star-small}} using \cite[Corollary 2.12]{janson2010asymptotic}.
\qed

\vskip3pt

Given $\cG_{\IRG}^{\sss(n), -}$, construct $\cG^{\sss(n)}_{\IRG}$ by placing edges independently with probability $1-\exp(-n^{-1-\delta})$ between distinct vertices $i, j\in [n]$ such that $\{i,j\}$ is not an edge in $\cG_{\IRG}^{\sss(n), -}$. 
Thus, for any two distinct components $\cC_i^-$ and $\cC_j^-$ in $\cG_{\IRG}^{\sss(n), -}$, the number of edges added between them in $\cG_{\IRG}^{\sss(n)}$, say $N_{ij}$, follows a $\mathrm{Binomial}\big(|\cC_i^-| |\cC_j^-|, 1-\exp(-1/n^{1+\delta})\big)$ distribution, 
and given $N_{ij}$, the edges that link the two components can be generated by drawing a sample of size $N_{ij}$ without replacement from the set of $|\cC_i^-| |\cC_j^-|$ many possible edges.
Further, for any component $\cC_i^-$ in $\cG_{\IRG}^{\sss(n), -}$, 
$N_{ii}\sim
\mathrm{Binomial}\Big(\Big(\begin{array}{c}|\cC_i^-|\\ 2\end{array}\Big), 1-\exp(-1/n^{1+\delta})\Big)$ 
many edges are added between the vertices of $\cC_i^-$ in $\cG_{\IRG}^{\sss(n)}$

Our plan is to apply Theorem \ref{thm:aldous-gen-2} with
\begin{equation}\label{eqn:def-x-i}
	x_i = n^{-2/3}\beta^{1/3}|\cC_i^{-}|,\ 
	q =\beta^{-2/3}n^{1/3-\delta},
	\ \text{ and }\
	M_i =  \scl(1, 1/|\cC_i^{-}|) \cC_i^{-}.
\end{equation}
Following the notation introduced around \eqref{eqn:uik-def}, let $\bar \cC_i$ be the $i$-th maximal component (in terms of the number of vertices) of $\bar \cG(\vx, q, \vM)$.
We will write $\bar \cC_i$ to denote both the graph and the metric measure space (where the metric is given by \eqref{eqn:44} and coincides with the graph distance, and the measure is as described in \eqref{eqn:barmu-on-full-met} and it coincides with $n^{-2/3}\beta^{1/3}$ times the counting measure on the vertices).
Now, in the setup of Theorem \ref{thm:aldous-gen-2}, we only ever place {\bf one} edge between distinct components $\cC_i^-$ and $\cC_j^-$ of $\cG_{\IRG}^{\sss(n), -}$ with probability $1-\exp(-|\cC_i^-||\cC_j^-|/n^{1+\delta})$;
if such an edge is added, it is distributed uniformly on the set of $|\cC_i^-||\cC_j^-|$ many possible edges.
There is a natural coupling in which 
$V\big(\cC_i(\cG_{\IRG}^{\sss(n)})\big)=V(\bar\cC_i)$, 
and 
$E\big(\cC_i(\cG_{\IRG}^{\sss(n)})\big)\supseteq E(\bar\cC_i)$ for all $i\geq 1$.
In the next lemma we show that after proper rescaling, $\cC_i(\cG_{\IRG}^{\sss(n)})$ and $\bar\cC_i$ are close with respect to Gromov-Hausdorff-Prohorov distance in the above coupling. 
This will enable us to apply Theorem \ref{thm:aldous-gen-2}.
\begin{lem}\label{lem:multiple-edges}
	Fix $k\geq 1$. Then, in the above coupling,
	\[
	d_{\GHP}\big(
	\scl\big(n^{-1/3},\, n^{-2/3}\beta^{1/3}\big)\cC_k(\cG_{\IRG}^{\sss(n)}), \,
	n^{-1/3}\bar\cC_k\big)\weakc 0\, .
	\]
\end{lem}
We first prove Theorems \ref{thm:component-limit-irg} and \ref{thm:scaling-limit-irg} assuming Lemma \ref{lem:multiple-edges}.

\noindent\textbf{Proof of Theorem \ref{thm:component-limit-irg} and Theorem \ref{thm:scaling-limit-irg}:} 
Recall the parameters from \eqref{eqn:def-x-i}.
As before, write $\sigma_k=\sum_i x_i^k$. Here we have
\begin{equation}\label{eqn:99}
	\sigma_2 = \frac{\beta^{2/3} \bar s_2}{n^{1/3}},\
	\sigma_3 = \frac{\beta\bar s_3}{n},\
	x_{\max} =\frac{\beta^{1/3}|\cC_1^{-}|}{n^{2/3}},\text{ and }
	x_{\min} \geq \frac{\beta^{1/3}}{n^{2/3}}.
\end{equation}
By Corollary \ref{cor:irg-barely-subcritical-original} and \eqref{eqn:prop-barely-sub-irg-1},
Condition \ref{ass:aldous-basic-assumption} holds with $\lambda=\zeta\beta^{-2/3}$. Theorem \ref{thm:component-limit-irg} now follows from Theorem \ref{thm:aldous-review}.

In view of Lemma \ref{lem:multiple-edges}, it is enough to show that Condition \ref{ass:aldous-gen-2} holds.
\ch{Now, by Corollary \ref{cor:irg-barely-subcritical-original} and \eqref{eqn:prop-barely-sub-irg-1},
	\begin{equation}\label{eqn:100}
		\bar s_2 \sim n^{\delta},\
		\bar s_3 \sim \beta n^{3\delta},\text{ and }
		|\cC_1^{-}| =o_P(n^{2\delta}\log n).
	\end{equation}
	From \eqref{eqn:99} and \eqref{eqn:100}, it follows that for any $\eta_0>0$,
	\[{x_{\max}}/{\sigma_2^{3/2+\eta_0}}
	=o_P\Big({n^{-\frac{1}{6}+\frac{\delta}{2}+\eta_0\big(\frac{1}{3}-\delta\big)}}\cdot{\log n}\Big).\]
	Since $\delta<1/5$, $\eta_0$ can be chosen small enough so that the last expression goes to zero.
	Similarly, it follows from \eqref{eqn:99} and \eqref{eqn:100} that $\sigma_2^{r_0}/x_{\min}\probc 0$ for some large $r_0$.
	This shows that Condition \ref{ass:aldous-gen-1} holds.}

By definition, $u_{\ell} =\sum_{i,j \in \cC_\ell^{-}} d^-(i,j)/|\cC_\ell^{-}|^{2}$ where $d^-$ denotes the graph distance in $\cG_{\IRG}^{\sss(n), -}$. Therefore, \ch{by \eqref{eqn:prop-barely-sub-irg-2} and Corollary \ref{cor:irg-barely-subcritical-original},}
\begin{equation*}
	\ch{\tau=}\sum_{\ell \geq 1} x_\ell^2 u_{\ell}
	= n^{-4/3}\beta^{2/3}\sum_{\ell \geq 1} \sum_{i,j \in \cC_\ell^{-}} d^{-}(i,j) = n^{-1/3}\beta^{2/3}\bar \cD
	\ch{\sim\alpha\beta^{2/3}n^{2\delta-1/3},}
\end{equation*}
\ch{and $\cD_{\max}^-=o_P(n^{\delta}\log n)$. {The above asymptotics then directly imply},
	\[
	\frac{\cD_{\max}^-\sigma_2^{\frac{3}{2}-\eta_0}}{\tau+\sigma_2}
	=o_P\Big(n^{-\big(\frac{1}{3}-\delta\big)\big(\frac{1}{2}-\eta_0\big)}\cdot \log n\Big), \text{ and }
	\frac{\sigma_2 x_{\max}\cD_{\max}^-}{\tau}=o_P\Big(\frac{(\log n)^2}{n^{\frac{2}{3}-2\delta}}\Big).
	\]
	Since $\delta<1/5$, it follows that (b) of Condition \ref{ass:aldous-gen-2} holds for any choice of $\eta_0<1/2$.}
This completes the proof of Theorem \ref{thm:scaling-limit-irg}.\qed

\noindent{\bf Proof of Lemma \ref{lem:multiple-edges}:} 
Recall that $N_{ii}$ denotes the number of edges added between the vertices of $\cC_i^-$ while going from $\cG_{\IRG}^{\sss(n), -}$ to $\cG_{\IRG}^{\sss(n)}$. 
Similarly, for $i\neq j$, $N_{ij}$ denotes the number of edges between $\cC_i^-$ and $\cC_j^-$ in $\cG_{\IRG}^{\sss(n)}$. 
Let $\cF_{-}$ denote the $\sigma$-field generated by $\cG_{\IRG}^{\sss(n), -}$. 
Define $X_n:=\sum_{i\neq j}\ind\set{N_{ij}\geq 2}+\sum_{i}\ind\set{N_{ii}\geq 1}$.
Then, for any $k\geq 1$ and $x, y\in V(\bar\cC_k)$,
$|d(x, y)-\bar d(x,y)|\leq 2X_n\cD_{\max}^-$, 
where $d$ and $\bar d$ denote the metrics in $\cC_k(\cG_{\IRG}^{\sss(n)})$ and $\bar\cC_k$ respectively.
Thus,
$d_{\GH}\big(\cC_k(\cG_{\IRG}^{\sss(n)}), {\bar\cC}_k\big)\leq X_n\cD_{\max}^-$.
From Corollary \ref{cor:irg-barely-subcritical-original}, $n^{-1/3}\cD_{\max}^-\weakc 0$. So it is enough to show that $X_n$ is tight. To this end, note that
\[\pr(N_{ij}\geq 2|\cF_-)\leq |\cC_i^-|^2|\cC_j^-|^2/n^{2+2\delta}\text{ and }\pr(N_{ii}\geq 1|\cF_-)\leq |\cC_i^-|^2/n^{1+\delta}.\]
Hence, $\E\left[X_n|\cF_-\right]\leq \bar s_2^2/n^{2\delta}+\bar s_2/n^{\delta}$. 
Now, an application of Corollary \ref{cor:irg-barely-subcritical-original} will show that $\bar s_2/n^{\delta}\weakc 1$. This proves tightness of $X_n$. Hence, we have shown that $d_{\GH}\big(\cC_k(\cG_{\IRG}^{\sss(n)}), {\bar\cC}_k\big)=o_P(n^{1/3})$ for any $k\geq 1$. 
Now the corresponding statement for $d_{\GHP}$ follows trivially.
\qed

\subsection{Branching process approximation}\label{sec:irg-branching-process}
As has been observed in \cite{BBSJOR07}, one key tool in {the study of} the $\IRG$ model is a closely related multitype branching process.
The aim of this section is to introduce this object and study its properties in the barely subcritical regime. 

\chh{First,} for any graph $G$ and a vertex $w \in V(G)$, define $\cC(w ; G)$ to be the component in $G$ that contains the vertex $w$. Denote by $d$ the graph distance on $G$. Define
$\cD(w; G) := \sum_{z \in \cC(w;G)} d(w,z)$.
\chh{To ease notation we suppress dependence on $n$ and write}  $\cG^{\star}_{\IRG}:=\cG_{\IRG}^{(n), \star}$. Let
$	\cC(i) = \cC(i;\cG_{\IRG}^{\star})$ and $\cD(i) = \cD(i; \cG_{\IRG}^{\star})$.
Let $v$ and $u$ be two uniformly chosen vertices from $[n]$, independent of each other and of $\cG_{\IRG}^{\star}$. 
Let $\cF_{\star}$ be the $\sigma$-field generated by $\cG_{\IRG}^{\star}$. 
Then
\begin{gather}
	\label{eqn:1162}
	\bar s_{k+1}^{\star} = \E\big[ |\cC(v)|^k \big| \cF_{\star}\big]
	\ \text{ for }\ k\geq 1\, , \ \
	\bar \cD^{\star}= \E[ \cD(v)| \cF_{\star}]\, ,\\
	\label{eqn:1166}
	(\bar s_{k+1}^{\star})^2 = \E\big[ |\cC(v)|^k |\cC(u)|^k \big| \cF_{\star} \big]
	\ \text{ for }\ k\geq 1,\ \text{ and }\ 
	(\bar \cD^{\star})^2 = \E\left[ \cD(v)\cD(u)\big| \cF_{\star}\right]\, .
\end{gather}

\begin{prop}\label{prop:irg-moments-true}
	We have,
	\begin{equation}\label{eqn:jh}
		\lim_n\frac{\E[\bar s_{2}^{\star}]-n^{\delta}}{n^{2\delta-1/3}}=\zeta,\
		\lim_n\frac{\E[\bar s_{3}^{\star}]}{n^{3\delta}}=\beta,
		\text{ and }\lim_n\frac{\E[\bar \cD^{\star}]}{n^{2\delta}}=\alpha.
	\end{equation}
	In addition, there exist positive constants $C_1=C_1(k, \kappa, \mu), C_2=C_2(\kappa, \mu)$, and some positive integer $n_0$ such that for all $n\geq n_0$,
	\begin{align}\label{eqn:var-s-bound}
		&\var(\bar s_{k+1}^{\star}) \leq C_1 n^{(4k+1)\delta-1},\text{ and }\\
		\label{eqn:var-D-bound}
		&\hskip25pt \var(\bar \cD^{\star}) \leq C_2 n^{(8\delta-1)}.
	\end{align}
	The cut-off $n_0$ depends only on the sequences $\{\mu_n\}$ and $\{\kappa_n\}$.
\end{prop}

It turns out that bounding $\var(\bar \cD^\star)$ is the most difficult part. We will prove the other asymptotics  first and leave this for the end of this section.  We start with the following lemma.
\begin{lem}\label{lem:bound-variance-s}
	For all $k \geq 1$ and all $n$, we have
	\begin{align}
		\var(\bar s_{k+1}^{\star})\leq \frac{1}{n}\E[ |\cC(v)|^{2k+1}]\text{ and}\label{eqn:var-s-bound-auxilliary}\\
		\cov(\bar s_{k+1}^{\star}, \bar \cD^{\star}) \leq \frac{1}{n}\E[|\cC(v)|^{k+1}\cD(v)].\label{eqn:cov-s-D}
	\end{align}
\end{lem}
\noindent\textbf{Proof:} By \eqref{eqn:1166}, we have $\E[(\bar{s}^{\star}_{k+1})^2] = \E[|\cC(v)|^k |\cC(u)|^k]$.
Write $\cG'$ for the graph induced by $\cG_{\IRG}^{\star}$ on the vertex set $[n]\setminus V(\cC(v))$. 
Then
\begin{align*}
	\E\Big[ |\cC(u)|^k \big| \set{\cC(v), v} \Big]
	&= \frac{|\cC(v)|}{n} \cdot |\cC(v)|^k + \frac{1}{n} \E\Big[\sum_{\cC \subseteq \cG'} |\cC|^{k+1} \big|\set{\cC(v), v}\Big],\\
	&\leq \frac{1}{n}|\cC(v)|^{k+1} + \E[\bar s_{k+1}^{\star}] = \frac{1}{n}|\cC(v)|^{k+1} + \E[ |\cC(v)|^{k} ],
\end{align*}
where $\sum_{\cC \subseteq \cG'}$ denotes sum over all components in $\cG'$. Therefore,
\begin{align*}
	&\E\left[(\bar s_{k+1}^\star)^2\right]= \E\left[|\cC(u)|^k |\cC(v)|^k\right]
	= 
	\E\left[ |\cC(v)|^k \E\left(|\cC(u)|^k \big| \set{\cC(v), v}\right)\right]\\
	&\hskip20pt \leq 
	\frac{1}{n}\E\big[|\cC(v)|^{2k+1}\big] + \big(\E \big[ |\cC(v)|^{k} \big]\big)^2
	= 
	\frac{1}{n}\E\big[|\cC(v)|^{2k+1}\big]  {+ \big(\bE[\bar s_{k+1}^{\star}]\big)^2},
\end{align*}
which gives the {desired} bound on $\var(\bar s_{k+1}^\star)$. Similarly,
\begin{align*}
	&\E\big[\bar s_{k+1}^{\star} \bar \cD^{\star}\big]=\E\big[\big|\cC(u)\big|^k\cD(v)\big]
	=\E\Big[\cD(v) \E\big(|\cC(u)|^k\big| \set{\cC(v), v}\big)\Big] \\
	&\le\frac{1}{n}\E\big[\cD(v)|\cC(v)|^{k+1}\big] + \E\big[ |\cC(v)|^{k} \big] \E[ \cD(v) ]
	=
	\frac{1}{n}\E\big[\cD(v)|\cC(v)|^{k+1}\big] +\bE\big[\bar s_{k+1}^{\star} \big]\bE\big[\bar \cD^{\star}\big].
\end{align*}
This completes the proof of Lemma \ref{lem:bound-variance-s}. \qed

Recall the definition of $\kappa_n^{-}$ from \eqref{eqn:def-kappa^-}.
We now consider a $K$-type branching process in which each particle of type $y \in [K]$ in any generation has a Binomial$\left(n\mu_n(z),\ \kappa_n^{-}(z, y)/n\right)$ number of type $z$ children in the next generation for $z \in [K]$, and the numbers of children of different types are independent.
Suppose in the $0$-th generation, there is only one particle and its type is $x \in [K]$.
For $k\geq 0$, define $G_k(x)=G_k(x; n, \mu_n, \kappa_n^-)$ to be the total number of particles in the $k$-th generation of such a branching process. (Thus, $G_0(x) \equiv 1$.) For $k \geq 0$, define
$T_k(x)=T_k(x; n, \mu_n, \kappa_n^-):= \sum_{\ell = 0}^\infty \ell^k G_\ell(x)$,
and $H_k(x)=\sum_{s=1}^k G_s(x)$.
Denote by \ch{$G_k(\mu_n)$} and $T_k(\mu_n)$ the corresponding quantities for the branching process when the type of the first particle follows the distribution $\mu_n$.
\ch{Note that $T_0(x)$ (resp. $T_0(\mu_n)$) represents the total progeny of such a branching process started from only one particle of type $x$ (resp. whose type has distribution $\mu_n$).}

Given a random vector $\mvw=(w_1,\hdots, w_K)^t$ with $w_y\geq 0$,
\ch{consider $\sum_{y\in [K]}w_y$ independent branching processes
	among which $w_y$ start from one particle of type $y$, $y\in[K]$.}
For $k\geq 0$, define
\vskip-.2in
\begin{align}\label{eqn:def-G-k-mvw}
	G_k(\mvw)=\sum_{y\in[K]}\sum_{i=1}^{w_y}G_k^{(i)}(y), \text{ and } 
	H_k(\mvw)=\sum_{s=1}^k G_s(\mvw),
\end{align}
\vskip-.1in
\noindent where $(G_k^{(i)}(y), k\geq 0)\equald (G_k(y), k\geq 0)$
for $1\leq i\leq w_y$ and  conditional on $\mvw$, the sequences
$(G_k^{(i)}(y), k\geq 0)$ are independent across $y\in[K]$ and $1\leq i\leq w_y$.
\ch{In words, $G_k(\mvw)$ (resp. $H_k(\mvw)$) represents the total number of particles in generation $k$
	(resp. total number of particles between generation $1$ and generation $k$) of such a system \chh{of} branching processes.}
Analogously define
\vskip-.3in
\begin{align}\label{eqn:def-T-0-mvw}
	T_0(\mvw)=\sum_{y\in[K]}\sum_{i=1}^{w_y}T_0^{(i)}(y),
\end{align}
where $T_0^{(i)}(y)\equald T_0(y)$ for $1\leq i\leq w_y$,
and the random variables $T_0^{(i)}(y),\ y\in[K], 1\leq i\leq w_y$, are independent conditional on $\mvw$.
\ch{Thus, $T_0(\mvw)$ represents the total progeny in such a system of branching processes.}
\chh{To prove Proposition \ref{prop:irg-moments-true}, we will couple the IRG with a $K$-type branching process. The asymptotics given in the next lemma will be useful in the proof.}

\begin{lem}\label{lem:branching-moments}
	\noindent{\upshape(a)} {\bf Growth rates for $T_k(\mu_n)$:} For any $r, k\geq 0$, there exists a constant $C_1=C_1(r, k; \kappa, \mu)>0$ such that 
	\begin{equation}
		\label{eqn:branching-moments-bound-2}
		\E[T_0^r(\mu_n)T_{k}(\mu_n)]\leq C_1 n ^{(2r + k + 1)\delta}
	\end{equation}
	for all $n\geq n_0$, where $n_0$ depends only on $k, r$ and the sequences $\{\mu_n\}$ and $\{\kappa_n\}$. In particular, for any $r\geq 1$ and $x\in[K]$, there exists a constant $C_2=C_2(r; \kappa, \mu)>0$ such that
	\begin{equation}\label{eqn:branching-moments-bound-1}
		\sup_{n\geq n_1}n^{-(2r-1)\delta}\E[ T_0^r(x)]\leq C_2
	\end{equation}
	for some $n_1$ depending only on $r$ and the sequences $\{\mu_n\}$ and $\{\kappa_n\}$. Further, for any $J>0$, $\delta'>\delta$, and integers $r, k\geq 0$, there exists a constant $C_3=C_3(J,\delta', r,k; \kappa, \mu)>0$ and an integer $n_2$ depending only on $J, \delta', r, k$, and the sequences $\{\mu_n\}$ and $\{\kappa_n\}$ such that for all $n\geq n_2$,
	\begin{align}\label{eqn:T-0-T-k-mu-n'}
		\E\left[T_0^r(x; n, \mu_n', \kappa_n^-)\times T_{k}(x; n, \mu_n', \kappa_n^-)\right]\leq C_3 n ^{(2r + k + 1)\delta}
	\end{align}
	for any $x\in[K]$ and any sequence of measures $\{\mu_n'\}$ on $[K]$ satisfying 
	\begin{align}\label{eqn:29}
		\sum_{x\in[K]}|\mu_n'(x)-\mu_n(x)|\leq J n^{-\delta'}
	\end{align} 
	for all $n$. 
	
	\noindent{\upshape(b)} {\bf Exact asymptotics for $T_0(\mu_n)$ and $T_0^2(\mu_n)$:} We have,
	\begin{equation}\label{eqn:branching-moments-limit}
		\lim_n n^{1/3-2\delta}\big(\E[T_0(\mu_n)]-n^{\delta}\big)=\zeta,\text{ and }
		\lim_n n^{-3\delta}\E[T_0^2(\mu_n)]=\beta.
	\end{equation}
	
	\noindent{\upshape(c)} {\bf Exact asymptotics for $T_1(\mu_n)$:} We have,
	\begin{equation}\label{eqn:T-1-limit}
		\lim_n n^{-2\delta}\E[T_1(\mu_n)]=\alpha.
	\end{equation}
	
	\noindent{\upshape(d)} {\bf Tail bound on height and component size:} For $x\in[K]$, let $\mathrm{ht}(x)$ denote the the height of the $K$-type branching process started from one initial particle of type $x$. Then there exist constants $C_4, C_5, C_6>0$ depending only on $\kappa$ \chh{and} $\mu$ such that for all $x\in[K]$ and $m\geq 1$, 
	\begin{align}
		&\pr(\mathrm{ht}(x)\geq m)\leq C_4\exp\big( - C_5 m/n^{\delta}\big)\text{ and }\label{eqn:ht-bound}\\
		&\hskip10pt\pr\left( T_0(x)\geq m \right) \leq 2 \exp\big( - C_6 m/n^{2\delta}\big)\label{eqn:component-tail-bound}
	\end{align}
	for $n\geq n_3$ where $n_3$ depends only on the sequences $\set{\mu_n}$ and $\set{\kappa_n}$.
\end{lem}
While proving $\eqref{eqn:var-D-bound}$, we will need an analogue of \eqref{eqn:branching-moments-bound-2} in the setting where the empirical distribution of types on $[K]$ may be different from $\mu_n$ but is sufficiently concentrated around $\mu_n$. This is the only part where we will use \eqref{eqn:T-0-T-k-mu-n'}. However, to avoid introducing additional notation, we will only prove \eqref{eqn:branching-moments-bound-2}. The proof of \eqref{eqn:T-0-T-k-mu-n'} for a sequence $\{\mu_n'\}$ satisfying \eqref{eqn:29} follows the exact same steps. We will continue to write $G_k(x), T_k(x)$ etc. without any ambiguity as the underlying empirical measure will always be $\mu_n$.

We will make use of the following result in the proof of \eqref{eqn:branching-moments-limit}.
\begin{lem}\label{lem:derivative-of-perron-root}
	Let $A_1, A_2$, and  $A_3$ be square matrices of order $K$. 
	Assume that the entries of $A_1$ are positive. Let $\vw_{\ell}$ and $\vw_r$ respectively be the left and right eigenvectors of $A_1$ corresponding to $\rho(A_1)$ subject to $\vw_{\ell}^t\vw_{r}=1$. 
	Then
	\[
	\lim_{(x, y)\to (0, 0)}\
	\frac{1}{y}\Big(\rho(A_1+xA_2+yA_3)-\rho(A_1+xA_2)\Big)=\vw_{\ell}^t A_3\vw_{r}.
	\]
\end{lem}

\noindent\textbf{Proof:} Since the entries of $A_1$ are positive, $\rho(A_1)$ is a simple eigenvalue of $A_1$.
An application of the implicit function theorem shows that $\rho(x, y):=\rho(A_1+xA_2+yA_3)$ is a $C^{\infty}$ function of $x, y$ in a small neighborhood of $(0, 0)$. So the required limit is simply $\partial\rho(0, 0)/\partial y$.

For some small $\eps>0$, let $\vw_{\ell}(y)\ \left(\text{resp. }\vw_{r}(y)\right):[-\eps, \eps]\to\bR^K$
be a $C^{\infty}$ function such that $\vw_{\ell}(0)=\vw_{\ell}$ (resp. $\vw_{r}(y)=\vw_{r}$) and for each
$y\in[-\eps, \eps]$, $\vw_{\ell}(y)$ (resp. $\vw_{r}(y)$) is a left (resp. right)
eigenvector of $A_1+yA_3$ corresponding to $\rho(0, y)$.
We further assume that $\vw_{\ell}(y)\cdot\vw_{r}(y)=1$ for $y\in[-\eps, \eps]$. Hence, we have
\begin{align}\label{eqn:2727}
	\Big(\frac{\partial}{\partial y}\vw_{\ell}(y)\Big)^t\vw_{r}(y)
	+\vw_{\ell}(y)^t\Big(\frac{\partial}{\partial y}\vw_{r}(y)\Big)=0\text{ for }y\in[-\eps, \eps].
\end{align}
Note that, $\vw_{\ell}(y)^t(A_1+yA_3)\vw_{r}(y)=\rho(0, y)$. Hence,
\begin{align*}
	&y\vw_{\ell}(y)^t A_3\vw_{r}(y)
	=\rho(0, y)-\vw_{\ell}(y)^t A_1\vw_{r}(y)\\
	&\hskip10pt=\rho(0, y)-\rho(0, 0)
	+\rho(A_1)\big(\vw_{\ell}-\vw_{\ell}(y)\big)^t\vw_{r}
	+\vw_{\ell}(y)^t A_1\big(\vw_r-\vw_{r}(y)\big).\notag
\end{align*}
The result follows upon dividing by $y$ and taking limits in the last equation and using \eqref{eqn:2727}.\qed

\noindent\textbf{Proof of Lemma \ref{lem:branching-moments}:}
For $x, y\in [K]$, define
$m_{xy}^{(n)}=\mu_n(y)\kappa_n^{-}(x, y)$, and let
$M_n = \big(m_{xy}^{(n)}\big)_{K \times K}$.
Note that for large $n$, $M_n$ is a matrix with positive entries.
Let $\rho_n=\rho(M_n)$, and let $\vu_n$ and $\vv_n$ be the associated right and left eigenvectors of $M_n$
respectively subject to $\vv_n^t\vu_n=\vu_n^t\vone=1$. Now,
\begin{align}\label{eqn:expectation-T-0-I}
	\E T_0(x)=\sum_{\ell=0}^{\infty}\E G_{\ell}(x)=\sum_{\ell=0}^{\infty}\ve_x^t M_n^{\ell}\vone
\end{align}
where $\ve_x$ denotes the unit vector with one at the $x$-th coordinate.
\ch{Recall the definitions of $M, \vu, \vv, A, B$, and $D$ from around \eqref{eqn:defn-u-v} and Condition \ref{ass:irg-strong}.}
From the Perron-Frobenius theorem
for positive matrices (see, e.g., \cite{horn-johnson}), it follows that there exist $c>0$ and $0<r<1$ such that
$M=\vu\vv^t+R$, where $R^t\vv=R\vu=\vzero$,
and $\max_{i, j}\big|R^{\ell}(i, j)\big|\leq cr^{\ell}$
for every $\ell\geq 1$. A similar decomposition holds for $M_n$:
\begin{align}\label{eqn:M-n-decomposition-I}
	M_n=\rho_n\vu_n\vv_n^t+R_n\text{ where }R_n^t\vv_n=R_n\vu_n=\vzero.
\end{align}
Since \ch{under Condition \ref{ass:irg-strong}}, $\max_{i,j}|m_{ij}^{(n)}-m_{ij}|=O(n^{-\delta})$ and similar statements are true for
$\|\vu_n-\vu\|$, $\|\vv_n-\vv\|$ and $(1-\rho_n)$, it follows that $\max_{i,j}|(R_n-R)(i, j)|=O(n^{-\delta})$.
Hence, there exist positive constants $c_1, c_2$ such that
\begin{align}\label{eqn:bound-on-R-n}
	\max_{i, j}\big|R_n^{\ell}(i, j)\big|\leq c_1(r+c_2 n^{-\delta})^{\ell}\text{ for }\ell\geq 1.
\end{align}
Using this decomposition, \eqref{eqn:expectation-T-0-I} yields
\begin{align}\label{eqn:expectation-T-0-II}
	\lim_n\frac{\E T_0(\mu_n)-n^{\delta}}{n^{2\delta-1/3}}
	=\lim_n\frac{\left((\bmu_n^t\vu_n)(\vv_n^t\vone)/(1-\rho_n)\right)-n^{\delta}}{n^{2\delta-1/3}}.
\end{align}
We can write
\begin{align}\label{eqn:M-n-decomposition}
	M_n&=\kappa_n^{-}D_n=\kappa D+\kappa(D_n-D)+(\kappa_n-\kappa)D_n+(\kappa_n^{-}-\kappa_n)D_n\notag\\
	&=M+\frac{\kappa B}{n^{1/3}}+\frac{AD}{n^{1/3}}-\frac{\vone\bmu^t}{n^{\delta}}+o(n^{-1/3}).
\end{align}
Note that
\begin{align}\label{eqn:n-1/3-rho-n}
	\rho_n=
	\rho\Big(M+\frac{\kappa B}{n^{1/3}}+\frac{AD}{n^{1/3}}-\frac{\vone\bmu^t}{n^{\delta}}\Big)
	+o(n^{-1/3}).
\end{align}
Lemma \ref{lem:derivative-of-perron-root} coupled with \eqref{eqn:n-1/3-rho-n} and \eqref{eqn:M-n-decomposition} gives
\begin{align}\label{eqn:derivative-rho-n}
	\lim_n n^{\delta}(1-\rho_n)=\lim_n n^{\delta}\left(\rho(M)-\rho(M_n)\right)=(\bmu^t\vu)(\vv^t\vone).
\end{align}

Using \eqref{eqn:derivative-rho-n} together with the facts $\|\bmu_n-\bmu\|=O(n^{-1/3})$,
$\|\vu_n-\vu\|=O(n^{-\delta})$, $\|\vv_n-\vv\|=O(n^{-\delta})$ and $\delta>1/6$,
we conclude from \eqref{eqn:expectation-T-0-II} and \eqref{eqn:derivative-rho-n} that
\begin{align}
	&
	\lim_n\frac{\E T_0(\mu_n)-n^{\delta}}{n^{2\delta-1/3}}
	=
	\lim_n\frac{\left((\bmu^t\vu)(\vv^t\vone)/(1-\rho_n)\right)-n^{\delta}}{n^{2\delta-1/3}}
	=
	\lim_n\frac{(\bmu^t\vu)(\vv^t\vone)-n^{\delta}(1-\rho_n)}{n^{2\delta-1/3}(1-\rho_n)}
	\notag\\
	&
	=\lim_n\Big[n^{\frac{1}{3}-\delta}-\frac{n^{\frac{1}{3}}(1-\rho_n)}{(\bmu^t\vu)(\vv^t\vone)}\Big]
	=\lim_n\Big[n^{\frac{1}{3}-\delta}-\frac{n^{\frac{1}{3}}\Big(1-\rho\left(M+\frac{\kappa B}{n^{1/3}}
		+\frac{AD}{n^{1/3}}-\frac{\vone\bmu^t}{n^{\delta}}\right)\Big)}{(\bmu^t\vu)(\vv^t\vone)}\Big],
	\label{eqn:expectation-T-0-III}
\end{align}
the last \ch{step} being a consequence of \eqref{eqn:n-1/3-rho-n}.

Since $f(x):=x^{-1}\left(1-\rho(M-x\vone\bmu^t)\right)$ is $C^{\infty}$ on a compact
interval around zero and $f(0)=(\bmu^t\vu)(\vv^t\vone)$, $|f(0)-f(x)|=O(|x|)$ on an interval
around zero. We thus have
\begin{align*}
	n^{1/3-\delta}&=\frac{n^{1/3-\delta}f(n^{-\delta})}{f(0)}
	+\frac{n^{1/3-\delta}\left[f(0)-f(n^{-\delta})\right]}{f(0)}\\
	&=\frac{n^{1/3}\left(1-\rho(M-n^{-\delta}\vone\bmu^t)\right)}{(\bmu^t\vu)(\vv^t\vone)}+O(n^{1/3-2\delta}).
\end{align*}
Plugging this in \eqref{eqn:expectation-T-0-III} and using Lemma \ref{lem:derivative-of-perron-root}, we get
\[\lim_n\frac{\E T_0(\mu_n)-n^{\delta}}{n^{2\delta-1/3}}
=\left(\vv^t(\kappa B+AD)\vu\right)/\left((\bmu^t\vu)(\vv^t\vone)\right).\]

This proves the first part of \eqref{eqn:branching-moments-limit}.
Here, we make note of the following fact:
\begin{align}\label{eqn:lim-T-0-x}
	\lim_n n^{-\delta}\E T_0(x)=u_x/\bmu^t\vu
\end{align}
which is a direct consequence of \eqref{eqn:expectation-T-0-I},
\eqref{eqn:M-n-decomposition-I}, \eqref{eqn:bound-on-R-n}, \eqref{eqn:derivative-rho-n}
and the facts that $\vu_n\to\vu$ and $\vv_n\to\vv$. We will need this result later.

To get the other part of \eqref{eqn:branching-moments-limit}, recall that for a
random vector $\vY=(Y_1,\hdots, Y_{r})^t$, the $p$-th order cumulants are given by
\begin{align}\label{eqn:def-cumulant}
	\cumu_{\vY}(\ell_1,\hdots, \ell_p)&=\cumu(Y_{\ell_1}, \hdots, Y_{\ell_p})\notag\\
	&:=\sum_{q=1}^p\sum\nolimits_1 (-1)^{q-1}(q-1)!\prod_{i=1}^q\E\Big(\prod_{j\in I_i}Y_{\ell_j}\Big)\, ,
\end{align}
where $1\leq\ell_i\leq r$ and $\sum_1$ denotes the sum over all partitions of
$I=\{1,\hdots, p\}$ into $q$ subsets $I_1, \hdots, I_q$.
Moments of $\vY$ can be expressed in terms of the cumulants as follows:
\begin{align}\label{eqn:moments-from-cumulant}
	\E\Big[\prod_{i=1}^p Y_{\ell_i}\Big]
	=\sum_{q=1}^p\sum\nolimits_1\prod_{i=1}^q\cumu_{\vY}\Big(\set{\ell_j}_{j\in I_i}\Big),
\end{align}
where $\sum_1$ has the same meaning as in \eqref{eqn:def-cumulant}.

For $x\in[K]$, let $Z(x, y)$, $y\in[K]$, be independent random variables with 
$Z(x, y)$ having $\mathrm{Binomial}(n\mu_n(y), \kappa_n^{-}(x, y)/n)$ distribution, and let $a_x(y_1, \hdots, y_q):=\cumu(Z(x, y_1), \hdots, Z(x, y_q))$ for $y_1, \hdots, y_q\in[K]$.
\ch{Let $T_0(x, y)$ be the total number of type-$y$ particles in the branching process starting from a particle of type $x$.}
Then it follows from \cite[Equation (13)]{lange1981moment} that
\begin{align}\label{eqn:lange}
	\cumu(T_0(x, y_1), \hdots, T_0(x, y_p))
	=\sum_{y\in[K]}m_{x, y}^{(n)}\cumu(T_0(y, y_1), \hdots, T_0(y, y_p))&\notag\\
	\hskip30pt
	+\sum_{q=2}^p\sum\nolimits_1\sum_{k_1,\hdots, k_q}a_x(k_1,\hdots, k_q)
	\prod_{m=1}^q\cumu\big(\set{T_0(k_m, y_j)}_{j\in I_m}\big),&
\end{align}
where $\sum_1$ is sum over all partitions of $I=\set{1,\hdots, p}$ into $q$ subsets $I_1,\hdots, I_q$.
For $p=2$, \eqref{eqn:lange} reduces to
\begin{align*}
	\cov\left(T_0(x, y_1), T_0(x, y_2)\right)=&\sum_{u\in[K]}m_{xu}^{(n)}\cov(T_0(u, y_1), T_0(u, y_2))\\
	&+\sum_{u\in[K]}\var(Z(x, u))(\E T_0(u, y_1))(\E T_0(u, y_2)).
\end{align*}
Summing both sides over all $y_1, y_2\in[K]$ and using the relations $\var(Z(x, u))=m_{xu}^{(n)}+O(n^{-1})$ and
$\max_{y\in[K]}\E T_0(y)=O(n^{\delta})$ \ch{from \eqref{eqn:lim-T-0-x}}, we get
\begin{align*}
	\var(T_0(x))&=\sum_{u\in[K]}m_{xu}^{(n)}\var(T_0(u))+\sum_{u\in[K]}m_{xu}^{(n)}\left[\E T_0(u)\right]^2+O(n^{2\delta}/n)\\
	&=\sum_{u\in[K]}m_{xu}^{(n)}\left[\E(T_0(u)^2)\right]+O(n^{2\delta}/n).
\end{align*}
Letting {$\va_1(n)=\left[\E T_0^2(x)\right]_{x\in[K]}$} and {$\va_2(n)=\left[\left(\E T_0(x)\right)^2\right]_{x\in[K]}$}, we have
\begin{align}\label{eqn:3465}
	(I-M_n)\va_1(n)=\va_2(n)+O(n^{2\delta}/n),
\end{align}
where the second term represents a vector with each coordinate $O(n^{2\delta}/n)$.
Since $\rho_n=\rho(M_n)<1$ for large $n$,
\begin{align}\label{eqn:inverse-M-n}
	(I-M_n)^{-1}=I+\sum_{k=1}^{\infty}M_n^k
	\chh{=\big[I+\big(\sum_{k\geq 1}\rho_n^k\big)\vu_n\vv_n^t+\sum_{k\geq 1}R_n^k\big],}
\end{align}
where the last step uses \eqref{eqn:M-n-decomposition-I}.
	It follows respectively from \eqref{eqn:bound-on-R-n}, \eqref{eqn:derivative-rho-n}, and \eqref{eqn:lim-T-0-x} that 
	\[\max_{i,j}\sum_{k\geq 1}|R_n^k(i,j)|=O(1),
	\sum_{k\geq 1}\rho_n^k=\frac{(1+o(1))n^{\delta}}{\bmu^t\vu\vv^t\vone},
	\  \text{ and }\ 
	\va_2(n)=(1+o(1))\cdot n^{2\delta}\cdot(\bmu^t\vu)^{-2}\cdot\big[u_x^2\big]_{x\in[K]}.
	\]
	These observations together with \eqref{eqn:3465} and \eqref{eqn:inverse-M-n} imply that
	\begin{align*}
		n^{-3\delta}\E T_0^2(\mu_n)
		=n^{-3\delta}\bmu_n^t\va_1(n)
		=n^{-3\delta}\bmu_n^t(I-M_n)^{-1}\big[\va_2(n)+O(n^{2\delta-1})\big]
		\to\beta.
	\end{align*}
	This completes the proof of \eqref{eqn:branching-moments-limit}.

Suppose we have proven that all cumulants (and hence all moments via \eqref{eqn:moments-from-cumulant}) of order $r$ are $O(n^{(2r-1)\delta})$ for $r\leq p-1$. To prove the same for $r=p$, note that the second term on the right side of \eqref{eqn:lange} is $O(\prod_{m=1}^q n^{(2|I_m|-1)\delta})=O(n^{(2p-q)\delta})=O(n^{(2p-2)\delta})$. 
As observed right below \eqref{eqn:inverse-M-n}, every entry of $(I-M_n)^{-1}$ is $O(n^{\delta})$. These two observations combined yield \eqref{eqn:branching-moments-bound-1}, \ch{by induction}.

Next, note that
\[\E T_1(\mu_n)
={\sum_{\ell=1}^{\infty}\sum_{x\in [K]}\ell\mu_n(x)\E G_{\ell}(x)}
=\sum_{\ell=1}^{\infty}\ell\bmu_n^t\Big(\rho_n^{\ell}\vu_n\vv_n^t+R_n^{\chh{\ell}}\Big)\vone.
\]
From Assumption \ref{ass:irg-strong} (b), \eqref{eqn:bound-on-R-n} and the facts $\|\vu_n-\vu\|+\|\vv_n-\vv\|=O(n^{-\delta})$, it follows that
\[\lim_n\frac{\E T_1(\mu_n)}{n^{2\delta}}
=(\bmu^t\vu\vv^t\vone)\lim_n\frac{1}{n^{2\delta}}\Big(\sum_{\ell=1}^{\infty}\ell\rho_n^{\ell}\Big)
=\lim_n\frac{(\bmu^t\vu\vv^t\vone)}{n^{2\delta}(1-\rho_n)^2}=\alpha,\]
where the last equality is a consequence of \eqref{eqn:derivative-rho-n}. This proves \eqref{eqn:T-1-limit}.

To prove \eqref{eqn:ht-bound}, notice that \eqref{eqn:derivative-rho-n} ensures the existence of $n_3\geq 1$ such that for $n\geq n_3$, we have $\rho_n\leq 1-C_5/n^{\delta}\leq \exp(-C_5/n^{\delta})$, where $C_5=(\bmu^t\vu\vv^t\vone)/2$. 
Now \eqref{eqn:M-n-decomposition-I} yields, for all $m\geq 1$,
\begin{align*}
	\pr(\mathrm{ht}(x)\geq m)&=\pr(G_m(x)\geq 1)\leq\E G_m(x)
	\leq C_4\rho_n^m\leq C_4\exp\Big(-C_5 m/n^{\delta}\Big)
\end{align*}
for $n\geq n_3$ and some constant $C_4>0$. 
This establishes \eqref{eqn:ht-bound}.

Now, \eqref{eqn:component-tail-bound} can be proved by imitating the proof of \cite[Lemma 6.13]{bhamidi2015aggregation}, and using \eqref{eqn:derivative-rho-n} and the fact \chh{that} $\|\kappa_n - n^{-\delta} \|_{L^2(\mu_n)}=\rho(M_n)$. Since no new idea is involved, we omit the proof. 

Finally, we prove \eqref{eqn:branching-moments-bound-2}. Note that
it is enough to prove the same bound on $\E[T_0^r(x)T_k(x)]$ for all $x \in [K]$. Consider one initial particle of type $x$. Let $N=G_1(x)$ and for $i = 1, 2, ..., N$, denote by $x_i$ the type of the $i$-th particle in generation one. 
Let $\Theta^{\sss (\ell)}(y)=\big(G_j^{\sss(\ell)}(y),\ j\geq 0\big)$, $y \in [K]$, $\ell\geq 1$, be independent random sequences that are also independent of $N, x_1,\ldots, x_N$ such that for all $\ell\geq 1$ and $y\in [K]$, 
$\Theta^{\sss (\ell)}(y)\equald\big(G_j(y),\ j\geq 0\big)$. 
For $k\geq 0$, let $T_k^{\sss(\ell)}(y)=\sum_{j\geq 1} j^kG_j^{\sss(\ell)}(y)$.
By the branching structure, we have
\begin{align*}
	T_0(x) \stackrel{d}{=}& 1 + \sum_{i = 1}^N T_0^{\sss(i)}(x_i),\ \text{ and}\\
	T_k(x) \stackrel{d}{=}& \sum_{i=1}^N \sum_{j=0}^\infty (j+1)^k G_j^{\sss(i)}(x_i) = \sum_{i=1}^N \sum_{j=0}^\infty \sum_{\ell = 0}^k {k \choose \ell} j^\ell G_j^{\sss(i)}(x_i) \\
	=& \sum_{i=1}^N \sum_{\ell = 0}^k {k \choose \ell} \Big[\sum_{j=0}^\infty  j^\ell G_j^{\sss(i)}(x_i)\Big] = \sum_{i=1}^N \sum_{\ell = 0}^k {k \choose \ell} T_{\ell}^{\sss(i)}(x_i), \mbox{ for } k = 1,2,...
\end{align*}
In the above sums, our convention is to take $0^0 = 1$.
The above distributional equalities also hold jointly. Observe that by \eqref{eqn:M-n-decomposition-I} and \eqref{eqn:derivative-rho-n}, for $k\geq 0$,
\[\E T_k(x)=\sum_{\ell\geq 0}\ell^k\ve_x^t\chh{M_n^{\ell}}\vone=O(\sum_{\ell\geq 0}\ell^k\chh{\rho_n^{\ell}})=O(n^{(k+1)\delta}).\]
So it is enough to prove \eqref{eqn:branching-moments-bound-2} for $r\geq 1, k\geq 1$. We will prove this by induction on $r+k$. First, we show the inductive step as follows. Assume \eqref{eqn:branching-moments-bound-1} and that for all $\set{(r',k'): r' < r \mbox{ or } k'<k}$,
\begin{equation}
	\label{eqn:induction-assumption}
	\E[T_0^{r'}(x)T_{k'}(x)] = O(n^{(2r'+k'+1)\delta}).
\end{equation}
Then for $(r,k)$, observe that
\begin{align}\label{eqn:6667}
	& (T_0(x)-1)^r T_k(x)
	\stackrel{d}{=}\Big[\sum_{i = 1}^N T_0^{\sss(i)}(x_i)\Big]^r \Big[  \sum_{i=1}^N \sum_{\ell = 0}^k {k \choose \ell} T_{\ell}^{\sss(i)}(x_i)\Big]\\
	&\hskip25pt= \Big[\sum_{r_1,...,r_N} {r \choose r_1,...,r_N} \prod_{j=1}^N (T_0^{\sss(j)}(\chh{x_j}))^{r_j} \Big] \Big[\sum_{i=1}^N \sum_{\ell = 0}^k {k \choose \ell} T_{\ell}^{\sss(i)}(x_i) \Big], \nonumber
\end{align}
\chh{where the summation $\sum_{r_1,...,r_N}$ is over $\vr:=(r_1, r_2, ...,r_N) \in \bN_0^N$, with $\sum_{i=1}^N r_i = r$.}
\chh{The right side of \eqref{eqn:6667} reduces to the following expression:}
\begin{equation}
	\label{eqn:1261}
	\sum_{r_1,...,r_N} \sum_{i=1}^N \sum_{\ell = 0}^k \Big[ {r \choose r_1,...,r_N}  {k \choose \ell} [T_0^{\sss(i)}(x_i)]^{r_i}T_{\ell}^{\sss(i)}(x_i) \prod_{j\neq i} (T_0^{\sss(j)}(x_j))^{r_j}\Big].
\end{equation}
By independence, writing $\cF_1$ for the $\sigma$-field generated by ${x_1,x_2,...,x_N}$,
\begin{align*}
	\E \Big[(T_0^{\sss(i)}(x_i))^{r_i}T_{\ell}^{\sss(i)}(x_i) \prod_{j\neq i} (T_0^{\sss(j)}(x_j))^{r_j} \Big| \cF_1 \Big]
	= \E[(T_0^{\sss(i)}(x_i))^{r_i}T_{\ell}^{\sss(i)}(x_i) | \cF_1]
	\prod_{j\neq i} \E[(T_0^{\sss(j)}(x_j))^{r_j} | \cF_1].
\end{align*}
Then whenever $\ell < k$ or $r_i < r$, we can apply the assumptions \eqref{eqn:induction-assumption}. Therefore,
\begin{equation}
	\E \Big[ (\chh{T_0^{\sss(i)}}(x_i))^{r_i}T_{\ell}^{\sss(i)}(x_i) \prod_{j\neq i} (T_0^{\sss(j)}(x_j))^{r_j}\Big| \cF_1 \Big] = O( n^{\phi(\vr,\ell,i) \delta}),
\end{equation}
where $\phi(\vr,\ell, i) = 2r + \ell+1 - |\set{j:r_j >0}| + \ind_{\set{r_i>0}}$. 
One can check that when $r_i < r$ or $\ell < k$, we have $\phi(\vr,\ell,i) \leq 2r + k$. Further, using \eqref{eqn:induction-assumption} again, we get
\begin{equation*}
	\E[(T_0(x)-1)^r T_k(x)] = \E[(T_0(x))^r T_k(x)]+O(n^{(2r+k-1)\delta}).
\end{equation*}
Therefore, from \eqref{eqn:1261}, we have
\begin{align*}
	\E[T_0(x)^r T_k(x)] =& \E\Big[\sum_{i=1}^N (T_0^{\sss(i)}(x_i))^{r}T_{k}^{\sss(i)}(x_i)\Big] + O (n^{(2r+k)\delta})\\
	=& \sum_{y \in [K]} m_{xy}^{(n)} \E[T_0(y)^r T_k(y)] + O (n^{(2r+k)\delta}).
\end{align*}
\chh{Thus $\Big[\E[T_0(x)^r T_k(x)]\Big]_{x\in[K]}=(I-M_n)^{-1}O(n^{(2r+k)\delta})$, where the last term represents a vector with norm $O(n^{(2r+k)\delta})$.
	Further, $\max_{i, j}|(I-M_n)^{-1}(i,j)|=O(n^{\delta})$ as observed right below \eqref{eqn:inverse-M-n}.
	This completes the induction step.}
Now we only need to bound $E[ T_0(x)T_1(x)]$. \chh{Using \eqref{eqn:6667} with $r=1$ and $k=1$}, we have
\begin{align*}
	T_0(x)T_1(x) - T_1(x) \stackrel{d}{=} \Big[\sum_{i = 1}^N T_0^{\sss(i)}(x_i)\Big] \Big[  \sum_{j=1}^N (T_0^{\sss(j)}(x_j) + T_1^{\sss(j)}(x_j))\Big].
\end{align*}
Now we can use the facts $\E[(T_0(x))^2] = O(n^{3\delta})$, $\E T_0(x)=O(n^{\delta})$ and $\E[T_1(x)] = O(n^{2\delta})$ to conclude that
$\E T_0(x)T_1(x)=O(n^{4\delta})$.
This proves the starting point of the induction and thus finishes the proof of \eqref{eqn:branching-moments-bound-2}. We have now completed the proof of Lemma \ref{lem:branching-moments}.\qed

\vskip3pt

The following lemma shows how closely we can approximate $\cG_{\IRG}^{\sss(n), \star}$ by the branching process. \chh{Fix $i\in[n]$. Recall that $x_i$ is the type of $i$ and $\cC(i)$ is the component of $i$ in $\cG_{\IRG}^{\sss (n), \star}$}.

\begin{lem}\label{lem:graph-to-branching}
	As before, let $v$ be a vertex chosen uniformly from $[n]$. 
	Then
	\begin{gather}
		\left|\E T_0(\mu_n)-\E|\cC(v)|\right|=O(n^{4\delta-1}),\
		\left|\E T_0(\mu_n)^2-\E|\cC(v)|^2\right|=O(\sqrt{n^{9\delta-1}}),
		\label{eqn:graph-to-branching-size}\\
		\text{ and }\ \ \
		\left|\E T_1(\mu_n)-\E\cD(v)\right|=O(n^{4\delta-1}). \label{eqn:graph-to-branching-distance}
	\end{gather}
	Further, for $r\geq 0$, $n\geq 1$, and {$i\in [n]$}, we have
	\begin{gather}
		\E\left[|\cC(v)|^r\cD(v)\right]\leq\E\left[T_0^r(\mu_n)T_1(\mu_n)\right],
		\text{ and} \label{eqn:component-times-distance-expectation}\\
		|\cC(i)|\leq_{st}T_0(x_i),\
		{\cD(i)\leq_{st} T_1(x_i),}\
		\diam(\cC(i))\leq_{st}2\times\mathrm{ht}(T_0(x_i)),\label{eqn:stoch-dom}
	\end{gather}
	where $X\leq_{st}Y$ means $Y$ dominates $X$ stochastically.
\end{lem}

We now set some notation which we will follow throughout the rest of this section.
For sequences $\set{a_m}$ and $\set{b_m}$, we will write ``$a_m\preceq_m b_m$" if there exists a positive constant $c$ depending only on $\kappa$ and $\mu$ and an integer $m_0$ depending only on the sequences $\set{\mu_n}$ and $\set{\kappa_n}$ such that $a_m\leq cb_m$ for $m\geq m_0$. If we have two sequences $\set{a_m(k)}_{m\geq 1}$ and $\set{b_m(k)}_{m\geq 1}$ for each $k\geq 1$, we will write ``$a_m(k)\preceq_m b_m(k)$ for $k\geq 1$" if $a_m(k)\leq cb_m(k)$ for $m\geq m_0$ and all $k\geq 1$ where $c$ and $m_0$ are as before. We emphasize that the same $c$ and $m_0$ work for all $k$.

We will use the following lemma in the proof of Lemma \ref{lem:graph-to-branching}.
Recall the definitions of $G_{\ell}(\mvw)$ and $H_{\ell}(\mvw)$ from \eqref{eqn:def-G-k-mvw}.
\begin{lem}\label{lem:square-of-generation}
	We have,
	\begin{equation}\label{eqn:G-ell-2}
		\E[G_{\ell}(x)^2]\preceq_n\rho_n^{\ell}/(1-\rho_n)\text{ for }x\in[K]\text{ and }\ell\geq 1,
	\end{equation}
	and for any non-random vector $\mvw=(w_y,\ y\in[K])$ with $w_y\geq 0$ for each $y\in[K]$,
	\begin{equation}\label{eqn:G-ell-H-ell}
		\E [G_{\ell}(\mvw)H_{\ell}(\mvw)]
		\preceq_n\frac{1}{1-\rho_n}
		\Big[\big(\vone\cdot\mvw\big)\ell\rho_n^{\ell}+\big(\vone\cdot\mvw\big)^2\rho_n^{\ell}\Big]\text{ for }\ell\geq 1.
	\end{equation}
\end{lem}
\noindent\textbf{Proof:}
Let $G_k(x, y)$ denote the number of type-$y$ particles in the $k$-th generation
of the multitype branching process started from a single particle of type $x$. 
Define
$\fN_k(x)^t=\left(G_k(x, y),\ y\in[K]\right)$.
{For any vector $\mvw$}, define the vector $\chh{\fN_k}(\mvw)$ in a similar fashion.

Let $\cF_{s}=\sigma\set{\chh{\fN_{k}}(x):\ 0\leq k\leq s}$ for $s\geq 0$. For any vector $\mvw=(w_1,\hdots, w_K)^t$, let
$\mvw^{(2)}=(w_1^2,\hdots, w_K^2)^t$ and $\|\mvw\|_{\infty}=\max_j w_j$. Also define $\mvw_k=M_n^k\mvw$ for $k\geq 0$.
From \eqref{eqn:M-n-decomposition-I} and \eqref{eqn:bound-on-R-n}, it follows that
\begin{equation}\label{eqn:M-n-^k-bound}
	\vone^t M_n^k\vone\preceq_n \rho_n^k\text{ and }
	\|\mvw_k\|_{\infty}\preceq_n\rho_n^k\|\mvw\|_{\infty}
	\text{ for }k\geq 1.
\end{equation}
Now,
\begin{align*}
	&\E\big[\big(\chh{\fN_{\ell}(x)}\cdot\mvw\big)^2\big]
	=\E\Big[\E\Big(\big(\chh{\fN_{\ell}(x)}\cdot\mvw\big)^2\Big|\cF_{\ell-1}\Big)\Big]\\
	&=\E\Big[\var\Big(\chh{\fN_{\ell}(x)}\cdot\mvw\Big|\cF_{\ell-1}\Big)\Big]
	+\E\big[\left(\chh{\fN_{\ell-1}(x)^t} M_n\mvw\right)^2\big]\\
	&=\E\bigg(\sum_{y\in[K]}w_y^2\var\left(G_{\ell}(x, y)\big|\cF_{\ell-1}\right)\bigg)+\E\big[\left(\chh{\fN_{\ell-1}(x)}\cdot\mvw_1\right)^2\big]\\
	&\leq\E\left[\chh{\fN_{\ell-1}(x)^t} M_n\mvw^{(2)}\right]
	+\E\big[\left(\chh{\fN_{\ell-1}(x)}\cdot\mvw_1\right)^2\big]
	=\ve_x^t M_n^{\ell}\mvw^{(2)}
	+\E\big[\left(\chh{\fN_{\ell-1}(x)}\cdot\mvw_1\right)^2\big].
\end{align*}
Using this recursion and making use of \eqref{eqn:M-n-^k-bound}, we get
\begin{align*}
	\E\big[\left(\chh{\fN_{\ell}(x)}\cdot\mvw\right)^2\big]
	&\leq\ve_x^t M_n^{\ell}\mvw^{(2)}+\ve_x^t M_n^{\ell-1}\mvw_1^{(2)}
	+\hdots+\ve_x^t M_n\mvw_{\ell-1}^{(2)}+(\ve_x\cdot\mvw_{\ell})^2\\
	&\preceq_n\|\mvw\|_{\infty}^2\left(\rho_n^{\ell}+\rho_n^{\ell+1}+\hdots+\rho_n^{2\ell-1}+\rho_n^{2\ell}\right)
	\leq\|\mvw\|_{\infty}^2\rho_n^{\ell}/(1-\rho_n).
\end{align*}
We get \eqref{eqn:G-ell-2} by taking $\mvw=\vone$. Next, note that
\begin{align*}
	&\E[G_{\ell}(\mvw)H_{\ell}(\mvw)]=\E[G_{\ell}(\mvw)^2]+\E[G_{\ell}(\mvw)H_{\ell-1}(\mvw)]\\
	&=\E[G_{\ell}(\mvw)^2]+\E\left[\left(\chh{\fN_{\ell-1}(\mvw)^t} M_n\vone\right)\times H_{\ell-1}(\mvw)\right]\\
	&\leq\E[G_{\ell}(\mvw)^2]
	+\E\left[\left(\vone^t M_n\vone\right)\times G_{\ell-1}(\mvw)^2\right]+\E\left[\left(\chh{\fN_{\ell-1}(\mvw)^t} M_n\vone\right)\times H_{\ell-2}(\mvw)\right],
\end{align*}
where the last step uses the inequality 
$\fN_{\ell-1}(\mvw)^t M_n\vone\leq G_{\ell-1}(\mvw)\cdot(\vone^t M_n\vone)$.
Proceeding in this way, we get
$
\E[G_{\ell}(\mvw)H_{\ell}(\mvw)]
\leq
\E[G_{\ell}(\mvw)^2]
+\sum_{k=1}^{\ell-1}(\vone^t M_n^k\vone)\times\E[G_{\ell-k}(\mvw)^2].
$
Since
\begin{align}\label{eqn:last-line}
	\E[G_{k}(\mvw)^2]&=\var[G_{k}(\mvw)]+[\E G_{k}(\mvw)]^2\\
	&\leq(\mvw\cdot\vone)\max_{y\in[K]}\ \E[G_k(y)^2]+(\mvw\cdot\vone)^2(\vone^t M_n^k\vone)^2,\notag
\end{align}
\eqref{eqn:G-ell-H-ell} follows by an application of \eqref{eqn:G-ell-2} and \eqref{eqn:M-n-^k-bound}.\qed

\vskip3pt

\noindent\textbf{Proof of Lemma \ref{lem:graph-to-branching}:}
It is enough to get bounds on
$\big|\E T_0(x_i)-\E|\cC(i)|\big|$, $\big|\E T_0(x_i)^2-\E|\cC(i)|^2\big|$, and
$\big|\E T_1(x_i)-\E\cD(i)\big|$ that are uniform in $i$. The proof proceeds via a coupling between
a $K$-type branching process and the breadth-first tree of $\cC(i)$. A similar coupling
in the \erdos case is standard and can be found in, for example, \cite{durrett2007random}.

Let $\eps_{jk}^{(\ell)}$ be independent Bernoulli$(\kappa_n^-(x_j, x_k)/n)$ random variables
for $1\leq j, k\leq n$ and $\ell\geq 1$. Let $\fI_{0}=\set{i}$ and $\fS_0=[n]\setminus\set{i}$.
Assume that we have defined $\fI_t$ and $\fS_t$ for $1\leq t\leq\ell-1$. For each $j\in\fI_{\ell-1}$
and $k\in\fS_{\ell-1}$, place an edge between $j$ and $k$ iff $\eps_{jk}^{(\ell)}=1$
\chh{and $j=\min\{u\in\fI_{\ell-1} : \eps_{uk}^{(\ell)}=1\}$.}
Let
\begin{align*}
	\fI_{\ell}=\big\{k\in\fS_{\ell-1}:\ \eps_{jk}^{(\ell)}=1\text{ for some }j\in\fI_{\ell-1}\big\}\text{ and }
	\fS_{\ell}=\fS_{\ell-1}\setminus\fI_{\ell}.
\end{align*}
\chh{Stop when $\fI_{\ell}=\emptyset$.
	This process returns a finite tree. Root this tree at $i$ and call it $\bbf(i)$.}
\ch{In words, $\fI_{\ell}$ corresponds to the set of new vertices found in the $\ell$-th step of a breadth-first exploration of $\cC(i)$ started from $i$, and $\fS_{\ell}$ corresponds to the set of vertices that have not been discovered up to step $\ell$. Note that the vertices in generation $\ell$ of $\bbf(i)$ is given by $\fI_{\ell}$, and that the set of vertices in $\cC(i)$ has the same distribution as $\bigcup_{\ell\geq 0}\fI_{\ell}$.}

Set $\fZ_{0}=\set{i}$.
\chh{We will now construct a $K$-type branching process tree $\bbp(i)$ rooted at $i$ in which $\fZ_{\ell}$, $\ell\geq 1$, will be the set of vertices in the $\ell$-th generation.
	To do so, we proceed inductively starting from $\ell=1$, and add three kinds of excess vertices to $\bbf(i)$:}

\noindent{\upshape(a)} \chhh{Assume that we have defined $\fZ_{t}$ for $1 \leq t \leq \ell-1$.} For each $u\in\fZ_{\ell-1}\setminus\fI_{\ell-1}$ and $y\in[K]$, create a collection $\fE_{uy}^{(\ell)}$
of type-$y$ children of $u$ independently, where $|\fE_{uy}^{(\ell)}|$ is distributed as a Binomial$(n\mu_n(x_u), \kappa_n^-(x_u, y)/n)$ random variable.
As usual, $x_u\in[K]$ denotes the type of $u$.
\ch{(Thus, $\fE_{uy}^{(\ell)}$ denotes the set of type-$y$ vertices spawned by a vertex $u$ present in the $(\ell-1)$-th generation of the branching process that is not in $\fI_{\ell-1}$.)}

\noindent{\upshape(b)} For each $j\in\fI_{\ell-1}$ and $k\in\fS_{\ell-1}^c$, create a child of $j$ of type $x_k$
iff $\eps_{jk}^{(\ell)}=1$. Call this collection of newly created vertices $\fB_{\ell}$.
\ch{(If there is an edge between a vertex $j\in\fI_{\ell-1}$ and another vertex, say $k$, discovered before or in step $\ell-1$ of the breadth-first exploration, then $k$ does not appear in $\fI_{\ell}$, but every such connection accounts for an extra vertex in the $\ell$-th generation of the coupled branching process;
	the set $\fB_{\ell}$ keeps track of that.)}

\noindent{\upshape(c)} For each $k\in\fS_{\ell-1}$ \chh{and $j\in\fI_{\ell-1}$, create a child of $j$ of type $x_k$ if $\eps_{jk}^{(\ell)}=1$ and $j\neq\min\{u\in\fI_{\ell-1} : \eps_{uk}^{(\ell)}=1\}$.}
Call this collection of newly created vertices $\fC_{\ell}$.
\ch{(If $s\geq 2$ vertices in $\fI_{\ell-1}$ have an edge common with some $k\in\fS_{\ell-1}$, then $k$ is only counted once in $\fI_{\ell}$, but this should account for an additional $s-1$ vertices in the coupled branching process.
	The set	$\fC_{\ell}$ keeps track of that.)}

Set
\[\fZ_{\ell}:=
(\fI_{\ell}\ \ \cup\ \ \fB_{\ell}\cup\fC_{\ell})\cup
\Big(\bigcup_{u\in\fZ_{\ell-1}\setminus\fI_{\ell-1},\ \ y\in[K]}\fE_{uy}^{(\ell)}\Big).\]

\chh{The resulting tree, which we call $\bbp(i)$, is a $K$-type branching process tree that starts from one particle of type $x_i$ as described right after Lemma \ref{lem:bound-variance-s}.
	The set of vertices in the $\ell$-th generation of $\bbp(i)$ is given by $\fZ_{\ell}$.
	For $\ell\geq 1$, define $\fA_{\ell}:=\fB_{\ell}\cup\fC_{\ell}$.
	Note that $\bbp(i)$ can be constructed by first adding the vertices $\cup_{\ell\geq 1}\fA_{\ell}$ to $\bbf(i)$, and then running a $K$-type branching process starting from each vertex in $\cup_{\ell\geq 1}\fA_{\ell}$ independently.}

\chh{For any vertex $v$ of $\bbp(i)$, let $T_0^{\circ}(v)$ denote the total number of descendants of $v$ (including $v$) in $\bbp(i)$.
	For $k\geq 0$, let $G_k^{\circ}(v)$ denote the number of descendants of $v$ in $\bbp(i)$ that are at distance $k$ from $v$, and let $H_k^{\circ}(v)=\sum_{s=1}^k G_s^{\circ}(v)$. 
	Then these random variables have the same distributions as $T_0(x_v)$, $G_k(x_v)$, and $H_k(x_v)$ respectively.
	For a subset $\cU$ of vertices of $\bbp(i)$, let $T_0^{\circ}(\cU)=\sum_{v\in\cU}T_0^{\circ}(v)$.
	Define $G_k^{\circ}(\cU)$ and $H_k^{\circ}(\cU)$ similarly.}

Let $\cF_0$ be the trivial $\sigma$-field and for $\ell\geq 1$, define
\begin{align*}
	\cF_{\ell}=\sigma\Big\{\fI_s, \fB_s, \fC_s, \fE_{uy}^{(s)}\, \Big|\, u\in\fZ_{s-1}\setminus\fI_{s-1},
	y\in[K], 1\leq s\leq\ell\Big\}.
\end{align*}
Write
\begin{align}
	&Z_{\ell}=|\fZ_{\ell}|, I_{\ell}=|\fI_{\ell}|, A_{\ell}=|\fA_{\ell}|,
	B_{\ell}=|\fB_{\ell}|,
	C_{\ell}=|\fC_{\ell}|, S_{\ell}=|\fS_{\ell}|,
	\text{ and }R_{\ell}=\sum_{j=0}^{\ell}I_j.\notag
\end{align}
Then,
$\E(B_{\ell}|\cF_{\ell-1})\preceq_n I_{\ell-1}R_{\ell-1}/n\leq Z_{\ell-1}(\sum_{s=0}^{\ell-1}Z_s)/n$. Hence,
\begin{align}\label{eqn:50}
	\sum_{\ell=1}^{\infty}\E B_{\ell}
	&\preceq_n\frac{1}{n}\sum_{\ell=0}^{\infty}\E\Big[Z_{\ell}(1+\chh{H_{\ell}^{\circ}(i)})\Big]
	=\frac{1}{n}\E\Big[T_0(x_i)+\sum_{\ell=0}^{\infty}G_{\ell}(x_i)H_{\ell}(x_i)\Big]
	\preceq_n \frac{n^{3\delta}}{n},
\end{align}
the last inequality being a consequence of \eqref{eqn:branching-moments-bound-1} and \eqref{eqn:G-ell-H-ell}.
We also have
$\E(C_{\ell}|\cF_{\ell-1})\preceq_n I_{\ell-1}^2\times S_{\ell-1}/n^2\leq Z_{\ell-1}^2/n.$
Hence,
\begin{align}\label{eqn:51}
	\sum_{\ell=1}^{\infty}\E C_{\ell}
	\preceq_n n^{-1}\sum_{\ell=0}^{\infty}\E(Z_{\ell}^2)
	\preceq_n n^{2\delta-1},
\end{align}
\chh{where the last step follows upon noting that $Z_{\ell}=G_{\ell}^{\circ}(i)\equald G_{\ell}(x_i)$ and then using \eqref{eqn:G-ell-2} and \eqref{eqn:derivative-rho-n}.}
Now
\begin{align*}
	\big|\E T_0(x_i)-\E|\cC(i)|\big|
	=\sum_{\ell=1}^{\infty}\E\big[\chh{T_0^{\circ}(\fA_{\ell})}\big]
	\preceq_n n^{\delta}\sum_{\ell=1}^{\infty}\E(A_{\ell})
	=n^{\delta}\sum_{\ell=1}^{\infty}\E(B_{\ell}+C_{\ell}),
\end{align*}
the second step being a consequence of \eqref{eqn:lim-T-0-x}. 
Combining this with \eqref{eqn:50} and \eqref{eqn:51}, we get the first half of \eqref{eqn:graph-to-branching-size}.
To get the other inequality, note that
\begin{align}\label{eqn:2746}
	\big|\E T_0(x_i)^2-\E|\cC(i)|^2\big|
	=\E\Big(\sum_{\ell=1}^{\infty}Z_{\ell}\Big)^2
	-\E\Big(\sum_{\ell=1}^{\infty}I_{\ell}\Big)^2\hskip75pt&\\
	\leq\Big[\E\Big(\sum_{\ell=1}^{\infty}(Z_{\ell}-I_{\ell})\Big)^2\Big]^{\frac{1}{2}}\times
	\Big[\E\Big(\sum_{\ell=1}^{\infty}(Z_{\ell}+I_{\ell})\Big)^2\Big]^{\frac{1}{2}}
	\preceq_n n^{\frac{3\delta}{2}}\Big[\E\Big(\sum_{\ell=1}^{\infty}(Z_{\ell}-I_{\ell})\Big)^2\Big]^{\frac{1}{2}},&\notag
\end{align}
where the last step follows from \eqref{eqn:branching-moments-bound-1}. Now,
\begin{align}\label{eqn:2747}
	\E\Big[\sum_{\ell=1}^{\infty}(Z_{\ell}-I_{\ell})\Big]^2
	=\E\Big[\sum_{\ell=1}^{\infty}\chh{T_0^{\circ}(\fA_{\ell})}\Big]^2
	\leq 
	2\E\Big[\sum_{\ell=1}^{\infty}\chh{T_0^{\circ}(\fA_{\ell})^2}
	+\sum_{1\leq\ell\leq s}\chh{T_0^{\circ}(\fA_{\ell})T_0^{\circ}(\fA_{s+1})}\Big].
\end{align}
By an argument similar to the one used in \eqref{eqn:last-line} and
the estimate from \eqref{eqn:branching-moments-bound-1},
\begin{align}\label{eqn:2748}
	\E\big[\chh{T_0^{\circ}(\fA_{\ell})^2}\big|\cF_{\ell}\big]\preceq_n A_{\ell}n^{3\delta}+A_{\ell}^2 n^{2\delta}.
\end{align}
Also,
\begin{align}\label{eqn:2749}
	\E\big[\chh{T_0^{\circ}(\fA_{\ell})T_0^{\circ}(\fA_{s+1})}\big]\preceq_n n^{2\delta}\E(A_{\ell}A_{s+1}).
\end{align}
From \eqref{eqn:50} and \eqref{eqn:51}, we have $\sum_{\ell\geq 1}\E A_{\ell}\preceq_n n^{3\delta}/n$. Further,
$A_{\ell}^2\leq 2(B_{\ell}^2+C_{\ell}^2)$ and
\begin{align*}
	\E\Big(B_{\ell}^2\big|\cF_{\ell-1}\Big)
	=\var(B_{\ell}\big|\cF_{\ell-1})+(\E(B_{\ell}\big|\cF_{\ell-1}))^2
	\preceq_n \frac{Z_{\ell-1}}{n}\Big(\sum_{s=0}^{\ell-1}Z_s\Big)
	+\Big[\frac{Z_{\ell-1}}{n}\Big(\sum_{s=0}^{\ell-1}Z_s\Big)\Big]^2.
\end{align*}
Hence,
\begin{align*}
	\sum_{\ell=1}^{\infty}\E(B_{\ell}^2)
	&\preceq_n n^{-1}\E\Big[\sum_{\ell=0}^{\infty}Z_{\ell}\Big(\sum_{s=0}^{\infty}Z_s\Big)\Big]
	+n^{-2}\E\Big[\sum_{\ell=0}^{\infty}Z_{\ell}^2\Big(\sum_{s=0}^{\infty}Z_s\Big)^2\Big]\\
	&\leq n^{-1}\E\Big[T_0(x_i)^2\Big]
	+n^{-2}\E\Big[T_0(x_i)^4\Big]\preceq_n n^{3\delta-1},
\end{align*}
where the final inequality uses \eqref{eqn:branching-moments-bound-1}
and the relation $\delta<1/5$. Similarly,
\begin{align*}
	\E(C_{\ell}^2\big|\cF_{\ell-1})
	&=\var(C_{\ell}\big|\cF_{\ell-1})+[\E(C_{\ell}\big|\cF_{\ell-1})]^2\\
	&\preceq_n I_{\ell-1}^2S_{\ell-1}/n^2+\left[I_{\ell-1}^2S_{\ell-1}/n^2\right]^2
	\leq Z_{\ell-1}^2/n+Z_{\ell-1}^4/n^2.
\end{align*}
Using \eqref{eqn:G-ell-2} and \eqref{eqn:branching-moments-bound-1}, we conclude that
\begin{align*}
	\sum_{\ell=1}^{\infty}\E(C_{\ell}^2)
	\preceq_n n^{-1}\sum_{\ell=0}^{\infty}\E Z_{\ell}^2
	+n^{-2}\E T_0(x_i)^4
	\preceq_n n^{2\delta-1}+n^{7\delta-2}.
\end{align*}
Combining these observations with \eqref{eqn:2747}, \eqref{eqn:2748}, and \eqref{eqn:2749}, we get
\begin{align}\label{eqn:2750}
	\E\Big(\sum_{\ell=1}^{\infty}(Z_{\ell}-I_{\ell})\Big)^2
	\preceq_n n^{6\delta-1}+n^{2\delta}\sum_{\ell=1}^{\infty}\sum_{s=\ell}^{\infty}
	\E\Big[A_{\ell}B_{s+1}+A_{\ell}C_{s+1}\Big].
\end{align}
Now,
\begin{align*}
	\E[A_{\ell}B_{s+1}]&\preceq_n\E\Big[A_{\ell}\Big(\frac{I_sR_s}{n}\Big)\Big]
	\leq\frac{1}{n}\E\Big[A_{\ell}\chh{G_{s-\ell}^{\circ}(\fI_{\ell})} \Big(R_{\ell}+\chh{H_{s-\ell}^{\circ}(\fI_{\ell})}\Big)\Big]\\
	&\preceq_n\frac{1}{n}\cdot\rho_n^{s-\ell}\E\left[A_{\ell}I_{\ell}R_{\ell}\right]
	+\frac{1}{n(1-\rho_n)}\E\left[A_{\ell}\left((s-\ell)\rho_n^{s-\ell}I_{\ell}+\rho_n^{s-\ell}I_{\ell}^2\right)\right],
\end{align*}
the last step being a consequence of \eqref{eqn:G-ell-H-ell} and the relation
\[
\E\big[G_{s-\ell}^{\circ}(\fI_{\ell})\big|\cF_{\ell}\big]\preceq_n I_{\ell}\times\max_{x\in[K]}\E G_{s-\ell}(x).
\]
Using \eqref{eqn:derivative-rho-n}, a simple computation yields
\begin{align}\label{eqn:2751}
	\sum_{\ell=1}^{\infty}\sum_{s=\ell}^{\infty}\E A_{\ell}B_{s+1}
	\preceq_n\sum_{\ell=1}^{\infty}
	\Big[{n^{\delta-1}}\E(A_{\ell}I_{\ell}R_{\ell})
	+{n^{3\delta-1}}\E(A_{\ell}I_{\ell})
	+{n^{2\delta-1}}\E(A_{\ell}I_{\ell}^2)\Big].
\end{align}
We can write
\begin{align}\label{eqn:2752}
	\sum_{\ell=1}^{\infty}\E(A_{\ell}I_{\ell}R_{\ell})
	=\sum_{\ell=1}^{\infty}\Big[\E(B_{\ell}I_{\ell}^2)+\E(B_{\ell}I_{\ell}R_{\ell-1})
	+\E(C_{\ell}I_{\ell}^2)+\E(C_{\ell}I_{\ell}R_{\ell-1})\Big].
\end{align}
To bound the first term on the right side, note that
\begin{align*}
	&\E(B_{\ell}I_{\ell}^2)=\E\left[\E\left(B_{\ell}\big|\cF_{\ell-1}\right)\E\left(I_{\ell}^2\big|\cF_{\ell-1}\right)\right]\\
	&\hskip20pt\preceq_n\E\Big[n^{-1}I_{\ell-1}R_{\ell-1}
	\cdot\big[\var\big(I_{\ell}\big|\cF_{\ell-1}\big)
	+\left(\E\left(I_{\ell}\big|\cF_{\ell-1}\right)\right)^2\big]\Big],
\end{align*}
where the first equality holds because of independence between $B_{\ell}$ and $I_{\ell}$ conditional on $\cF_{\ell-1}$.
Thus,
\begin{align}\label{eqn:2753}
	\sum_{\ell=1}^{\infty}\E(B_{\ell}I_{\ell}^2)
	&\preceq_n\sum_{\ell=1}^{\infty}\E\Big[\frac{I_{\ell-1}R_{\ell-1}}{n}
	\cdot(I_{\ell-1}+I_{\ell-1}^2)\Big]\leq\frac{2}{n}\sum_{\ell=1}^{\infty}\E\Big[I_{\ell-1}^3 R_{\ell-1}\Big]\notag\\
	&\leq\frac{2}{n}\sum_{\ell=0}^{\infty}\E\Big[Z_{\ell}^3\Big(\sum_{\ell=0}^{\infty}Z_{\ell}\Big)\Big]
	\leq\frac{2}{n}\E\Big[T_0(x_i)^4\Big]\preceq_n\frac{n^{7\delta}}{n},
\end{align}
by an application of \eqref{eqn:branching-moments-bound-1}. By a similar argument,
\begin{align}\label{eqn:2754}
	\sum_{\ell=1}^{\infty}\E\left[B_{\ell}I_{\ell}R_{\ell-1}\right]
	\preceq_n n^{7\delta-1}.
\end{align}
Next,
$
\E\big[C_{\ell}I_{\ell}^2\big|\cF_{\ell-1}\big]
\leq\sum\nolimits_1\E\big[\ind\big\{\eps_{j_1, k}^{(\ell)}=\eps_{j_2, k}^{(\ell)}=1\big\}
\cdot I_{\ell}^2\big|\cF_{\ell-1}\big]
$,
where $\sum_1$ stands for sum over all $j_1, j_2\in\fI_{\ell-1}$ and $k\in\fS_{\ell-1}$ such that $j_1\neq j_2$.
For any such $j_1, j_2, k$, {let
	\[
	V_{j_1,j_2;k}^{\sss (\ell)}=\big|\fI_{\ell}\setminus\{k\}\big|+\ind\{\eps_{j, k}^{(\ell)}=1\text{ for some }j\in\fI_{\ell-1}\setminus\{j_1, j_2\}\}.
	\]
	Then $V_{j_1,j_2;k}^{\sss (\ell)}$ is independent of $(\eps_{j_1, k}^{(\ell)}, \eps_{j_2, k}^{(\ell)})$
	conditional on $\cF_{\ell-1}$. 
	Further, $I_{\ell}-1\leq V_{j_1,j_2;k}^{\sss (\ell)}\leq I_{\ell}$.}
Hence,
\begin{align*}
	\E\big[C_{\ell}I_{\ell}^2\big|\cF_{\ell-1}\big]
	&\leq\sum\nolimits_1\E\big[\ind\big\{\eps_{j_1, k}^{(\ell)}=\eps_{j_2, k}^{(\ell)}=1\big\}\big|\cF_{\ell-1}\big]
	\cdot\E\big[\big(1+V_{j_1,j_2;k}^{(\ell)}\big)^2\big|\cF_{\ell-1}\big]\\
	&\preceq_n n^{-2}\times I_{\ell-1}^2 S_{\ell-1}\times\big[1+\E\big(I_{\ell}^2\big|\cF_{\ell-1}\big)\big]
	\preceq_n n^{-1}I_{\ell-1}^4.
\end{align*}
We thus have
\begin{align}\label{eqn:2755}
	\sum_{\ell=1}^{\infty}\E\left[C_{\ell}I_{\ell}^2\right]
	\preceq_n n^{-1}\sum_{\ell=0}^{\infty}\E\left[I_{\ell}^4\right]\leq n^{-1}\E\big[T_0(x_i)^4\big]\preceq_n n^{7\delta-1}.
\end{align}
We can similarly argue that
$
\sum_{\ell=1}^{\infty}\E\left[C_{\ell}I_{\ell}R_{\ell-1}\right]\preceq_n n^{7\delta-1}.
$
Combining this last observation with \eqref{eqn:2752}, \eqref{eqn:2753}, \eqref{eqn:2754}, and \eqref{eqn:2755}, we have
\begin{align}\label{eqn:387}
	\sum_{\ell=1}^{\infty}\E\left[A_{\ell}I_{\ell}R_{\ell}\right]\preceq_n n^{7\delta-1}.
\end{align}
We can use similar reasoning to bound the second and third terms on the right side of \eqref{eqn:2751}
and the term $\E(A_{\ell}C_{s+1})$ appearing on the right side of \eqref{eqn:2750};
we omit the details. The final estimates will be
\begin{align}\label{eqn:388}
	\sum_{\ell=1}^{\infty}\E\left[A_{\ell}I_{\ell}\right]\preceq_n n^{5\delta-1}
	,\ \sum_{\ell=1}^{\infty}\E\left[A_{\ell}I_{\ell}^2\right]\preceq_n n^{7\delta-1}
	\text{ and }\sum_{\ell=1}^{\infty}\sum_{s=\ell}^{\infty}\E\left[A_{\ell}C_{s+1}\right]\preceq_n n^{8\delta-2}.
\end{align}
We get the second inequality in \eqref{eqn:graph-to-branching-size} by combining
\eqref{eqn:2746}, \eqref{eqn:2750}, \eqref{eqn:2751}, \eqref{eqn:387} and \eqref{eqn:388}.

Next, note that $|\E T_1(x_i)-\E\cD(i)|=\sum_{\ell=1}^{\infty}\E\left[\ell Z_{\ell}-\ell I_{\ell}\right]
=\sum_{\ell=1}^{\infty}\ell\E B_{\ell}+\sum_{\ell=1}^{\infty}\ell\E C_{\ell},
$
so we can again argue similarly to get the estimate \eqref{eqn:graph-to-branching-distance}.
Finally, \eqref{eqn:component-times-distance-expectation} and \eqref{eqn:stoch-dom} are immediate from the coupling between $\bbf(i)$ and $\bbp(i)$. This completes the proof of Lemma \ref{lem:graph-to-branching}.\qed


We will need the following lemma to prove \eqref{eqn:var-D-bound}.
\begin{lem}\label{lem:irg-monotonicity-d}
	Fix $n\geq 1$ and as before, let $\cV=\cV^{\sss(n)} = [n]$. Let $\cV^- = [n]\setminus \set{1}$. Recall that for each $i \in [n]$, $x_i\in [K]$ denotes the type of the vertex $i$. Let $\bar\kappa$ be a kernel on $[K]\times[K]$. Let $\cG_1$ be the IRG model on the vertex set $\cV$ where we place an edge between $i, j\in[n], i\neq j,$ independently with probability $(\bar\kappa(x_i, x_j)/n)\wedge 1$. Let $\cG_0$ be the graph on the vertex set $\cV^-$ induced by $\cG_1$. Define $A := \max_{x, y\in[K]}\bar\kappa(x, y)$. Then, we have,
	\begin{equation*}
		\E[\cD(\cG_0)]  \leq  \E[\cD(\cG_1)] + {A^2} \E[ \cD(\cG_0) \cS_2(\cG_0)]/(2n^2).
	\end{equation*}
\end{lem}
\noindent\textbf{Proof:} Let $E$ denote the event that
there exists a component $\cC$ of $\cG_0$ such that there are at least two edges between vertex $1$ and the component $\cC$.
One important observation is that, on the event $E^c$, $\cD(\cG_0)\le\cD(\cG_1)$. Note also that the connection probability of each pair of vertices is bounded by $A/n$. Hence,
\begin{equation*}
	\pr(E|\cG_0) \leq \sum_{\cC \subseteq \cG_0} {|\cC| \choose 2} \left(\frac{A}{n}\right)^2 \leq \frac{A^2\cS_2(\cG_0)}{2n^2}\, .
\end{equation*}
Thus, we have
\begin{align*}
	\E[\cD(\cG_0)]
	=\E\left[\cD(\cG_0) \ind_{E^c}\right]+\E\left[\cD(\cG_0)\E(\ind_{E} \mid \cG_0)\right]
	\le\E[\cD(\cG_1)]+\frac{A^2}{2n^2} \E[\cD(\cG_0)\cS_2(\cG_0)].
\end{align*}
This completes the proof. \qed

\vskip3pt

\noindent\textbf{Proof of Proposition \ref{prop:irg-moments-true}:} Most of our work is already done. \eqref{eqn:jh} follows from \eqref{eqn:1162}, \eqref{eqn:graph-to-branching-size}, \eqref{eqn:graph-to-branching-distance}, \eqref{eqn:branching-moments-limit}, and \eqref{eqn:T-1-limit}. \eqref{eqn:var-s-bound} is a consequence of \eqref{eqn:var-s-bound-auxilliary}, \eqref{eqn:stoch-dom}, and \eqref{eqn:branching-moments-bound-1}. 
So, we turn directly to the proof of \eqref{eqn:var-D-bound}. 
Our goal is to bound $\bE[(\bar \cD^{\star})^2]=\E[\cD(u)\cD(v)]$.
Write $N=|V(\cC(V))|$, and let $V(\cC(v))=\{v_1, v_2, ...,v_N\}$.  
Define 
$\cV_0:= [n] \setminus V(\cC(v))$, and consider a sequence of sets $\cV_0 \subseteq \cV_1 \subseteq ... \subseteq \cV_{N} = [n]$ such that $\cV_i = \cV_{i-1} \cup \{ v_i\}$.
Write $\bE_1[\cdot]=\bE[\ \cdot\mid \{v, \cC(v) \} ]$. 
Define $\cov_1(\cdot,\cdot)$ similarly.
Let $\cG^{\star\star}\equald \cG^{\star}_{\IRG}$ be independent of $(u, v, \cG^{\star}_{\IRG})$, and for $0\leq i\leq N$, let $\cG_i$ be the restriction of $\cG^{\star\star}$ to the vertex set $\cV_i$ 
(thus, $\cG_N=\cG^{\star\star}$).
Note that
\begin{align}\label{eqn:1372}
	\bE_1\big[ \cD(u)\big]
	= 
	\frac{1}{n} \sum_{i,j \in \cC(v)}d_{\cC(v)}(i,j) 
	+ 
	\frac{1}{n} \E_1\big[\cD(\cG_0)\big].
\end{align}
Applying Lemma \ref{lem:irg-monotonicity-d} repeatedly, we have
\begin{align}\label{eqn:1383}
	\bE_1\big[ \cD(\cG_0) \big]
	\leq 
	\E_1[\cD(\cG_1) ] + \frac{A^2}{2} \E_1[ \bar \cD(\cG_0) s_2(\cG_0)]
	\leq
	\hdots
	\hskip110pt
	\\
	\leq 
	\E_1[\cD(\cG_N) ] + \frac{A^2}{2} \sum_{i=0}^{N-1} \E_1[ \bar \cD(\cG_i) s_2(\cG_i)]
	= 
	n\E[\cD(v) ] + \frac{A^2}{2} \sum_{i=0}^{N-1} \E_1[ \bar \cD(\cG_i) s_2(\cG_i)].
	\notag
\end{align}
Here, $A=2\max_{x, y\in[K]}\kappa(x, y)$, and each inequality holds for $n\geq n_0$, where $n_0$ depends only on the sequence $\{\kappa_n\}$.

Define the event $F_n=\set{ |\cC(v)|\leq Bn^{2\delta}\log n}$ where $B$ is a positive constant such that $n^5\pr(F_n^c)\to 0$. (This can be done because of \eqref{eqn:component-tail-bound} and \eqref{eqn:stoch-dom}.) 
\chh{For $x\in[K]$ and $0\leq i\leq N$, let $\mu_n^{\sss (i)}(x)$ denote the proportion of type-$x$ vertices in $\cV_i$. 
	Then
	$\mu_n(x)-N/n\leq \mu_n^{\sss (i)}(x)\leq n\mu_n(x)/(n-N)$.
	Using $\delta<1/5$, we see that on $F_n$,  \eqref{eqn:29} is satisfied by $\mu_n'=\mu_n^{\sss (i)}$ with $\delta'=3/5$ and $J=2BK$ for all $0\leq i\leq N$ and $n\geq n_{\star}$ where $n_{\star}$ depends only on $B$ and $\delta$. 
}
Now, note that $\cov_1(\bar \cD(\cG_i), s_2(\cG_i))$ can be bounded by following the proof techniques of \eqref{eqn:cov-s-D} and \eqref{eqn:component-times-distance-expectation}. Then an application of \eqref{eqn:T-0-T-k-mu-n'} {yields that on $F_n$}, $\cov_1(\bar \cD(\cG_i), s_2(\cG_i))\preceq_n n^{6\delta-1}$ for $0\leq i\leq N$. 
Similarly, using \eqref{eqn:T-0-T-k-mu-n'} and the first two relations in \eqref{eqn:stoch-dom}, on the event $F_n$,  $\E_1[s_2(\cG_i)]\preceq_n n^{\delta}$, and $\E_1 [\bar \cD(\cG_i)]\preceq_n n^{2\delta}$ for $0\leq i\leq N$. 
Thus, on $F_n$,
\begin{align*}
	\E_1\left[\bar\cD(\cG_i) s_2(\cG_i)\right]
	=\cov_1\left(\bar \cD(\cG_i), s_2(\cG_i)\right) + \E_1\left[\bar\cD(\cG_i)\right] \E_1\left[ s_2(\cG_i)\right]\preceq_n(n^{6\delta-1} + n^{3\delta})\preceq_n n^{3\delta}
\end{align*}
for $0\leq i\leq N$. 
This together with \eqref{eqn:1383} yields
\begin{align}\label{eqn:98}
	\ind_{F_n}\cdot\bE_1\big[ \cD(\cG_0) \big]
	\leq n\E\cD(v)+\eps_n,
\end{align}
where $\eps_n\preceq_n n^{5\delta}\log n$.
Hence,
\begin{align}\label{eqn:1395}
	\E[\cD(v)\cD(u)]
	&
	=
	\E[\cD(v)\cD(u) \ind_{F_n^c}] + \E[\cD(v)\cD(u)\ind_{F_n}]
	\notag\\
	&
	\leq 
	n^4 \pr(F_n^c) + \bE\big[\cD(v) \ind_{F_n} \E_1[ \cD(u)]\big] 
	\nonumber \\
	&
	\preceq_n 
	n^4 \pr(F_n^c)
	+
	\frac{1}{n}\E\Big[\cD(v)\sum_{i,j \in \cC(v)}d_{\cC(v)}(i,j)\Big]+[\E\cD(v)]^2
	+
	\frac{\eps_n}{n}\E\cD(v)
	\nonumber\\
	&
	=:
	n^4 \pr(F_n^c) +Q_n+[\E\cD(v)]^2+\eps_n\E\cD(v)/n \, ,
\end{align}
the third step being a consequence of \eqref{eqn:1372} and \eqref{eqn:98}. Since $n^5\pr(F_n^c)\to 0$ and $\E\cD(v)\preceq_n n^{2\delta}$ (by the last convergence in \eqref{eqn:jh}), we just need to show $Q_n\preceq_n n^{8\delta-1}$. To this end, note that
\begin{align*}
	\sum_{i, j\in\cC(v)}d_{\cC(v)}(i, j)
	\leq 
	2\sum_{i, j\in\cC(v)}d_{\cC(v)}(i, v)
	=
	2|\cC(v)|\sum_{i\in\cC(v)}d_{\cC(v)}(i, v)
	=
	2|\cC(v)|\cD(v).
\end{align*}
Further, we trivially have $\cD(v)\leq|\cC(v)|^2$. Thus,
\[Q_n\leq 2n^{-1}\E\left[|\cC(v)|^3\cD(v)\right]\preceq_n n^{8\delta-1},\]
by an application of \eqref{eqn:component-times-distance-expectation} and \eqref{eqn:branching-moments-bound-2}. This completes the proof of \eqref{eqn:var-D-bound}.\qed

\subsection{Proof of Theorem \ref{prop:irg-barely-subcritical}}\label{sec:irg-final-section}
First note that the claim $n^{-2\delta}\bar\cD^{\star}\weakc\alpha$ follows from \eqref{eqn:var-D-bound} and the last convergence in \eqref{eqn:jh}.

Next, by a simple union bound, $\pr(|\cC_1^{\star}|\geq m)\leq\sum_{i\in[n]}\pr(|\cC(i)|\geq m)$ for $m\geq 1$. Hence, the tail bound on $|\cC_1^{\star}|$ is immediate from  \eqref{eqn:component-tail-bound} and \eqref{eqn:stoch-dom}.
Similarly, the tail bound on $\cD_{\max}^{\star}$ follows from \eqref{eqn:ht-bound} and \eqref{eqn:stoch-dom}.

Since $2\delta>1/3$, the first convergence in \eqref{eqn:jh} shows that
\begin{align}\label{eqn:7826}
	\lim_n n^{-\delta}\E\bar s_2^\star=1,
\end{align}
which together with \eqref{eqn:jh} implies
\begin{align}\label{eqn:7827}
	\lim_n\ n^{1/3}\left(\frac{1}{n^{\delta}}-\frac{1}{\E\bar s_2^\star}\right)=\zeta.
\end{align}
Further, for each $\eps>0$, $\pr(|\bar s_2^\star-\E\bar s_2^\star|>\eps n^{\delta})\leq\eps^{-2}n^{-2\delta}\var(\bar s_2^\star)\to 0$ by \eqref{eqn:var-s-bound}. Hence,
\begin{align}\label{eqn:7828}
	n^{-\delta}\bar s_2^\star\weakc 1.
\end{align}
Now, for each $\eps>0$,
\begin{align}\label{eqn:7829}
	\pr\Big(\Big| n^{1/3}\Big(\frac{1}{\bar s_2^\star}-\frac{1}{\E\bar s_2^\star}\Big)\Big|>\eps\Big)
	&\leq\pr\Big(2n^{1/3}\frac{\left|\bar s_2^\star-\E\bar s_2^\star\right|}{n^{\delta}\E\bar s_2^\star}>\eps \Big)
	+\pr\Big(\bar s_2^\star\leq n^{\delta}/2\Big)\notag\\
	&\leq\frac{4n^{2/3}\var(\bar s_2^\star)}{\eps^2 n^{2\delta}(\E\bar s_2^\star)^2}
	+\pr\Big(\bar s_2^\star\leq n^{\delta}/2\Big)\to 0,
\end{align}
where the last convergence follows from \eqref{eqn:var-s-bound}, \eqref{eqn:7826}, and \eqref{eqn:7828}. 
The second convergence in \eqref{eqn:prop-barely-sub-irg-1} now follows from \eqref{eqn:7827} and \eqref{eqn:7829}.
Finally, the claim $\bar s_3^\star/(\bar s_2^\star)^3\weakc\beta$ is a simple consequence of \eqref{eqn:7828}, \eqref{eqn:var-s-bound}, and the second convergence in \eqref{eqn:jh}.
\qed

\section{Proofs: Scaling limits of CM}
\label{sec:proof-cm}
Section \ref{sec:proof-cm-max-diam} proves Theorem \ref{thm:config-largest-comp-diam}.
Section \ref{sec:cm-prop-entrance} uses Theorem \ref{thm:config-largest-comp-diam} to prove asymptotics for the susceptibility functions, namely, Theorem \ref{thm:config-bare-sibcrit}. 
The dynamic version of CM does not have the exact same merger dynamics as the multiplicative coalescent. In Section  \ref{sec:cm-mod-process-def} we construct a modification $\cG_n^{\modi}$ of the dynamic CM that does evolve as the multiplicative coalescent, and derive properties of this modified process including an application of the main universality result (Theorem \ref{thm:aldous-gen-2}). 
Section \ref{sec:proof-dynamic-cm-scaling-limit} derives properties of the dynamic CM by coupling it with both the percolated graph $\Perc_n$ and the modified process $\cG_n^{\modi}$, and proves Theorem \ref{thm:crit-main-res-cm} building on these results.   
Finally, Section \ref{sec:proof-perc-cm} proves Theorem \ref{thm:perc-cm} on critical percolation for CM.

\subsection{Preliminaries}
\label{sec:cm-proof-prelim}
Recall the construction of $\big(\CM_n(t);\, t\geq 0\big)$ from Section~\ref{sec:cm-def}, the parameters from \eqref{eqn:def-param}, and Condition \ref{ass:cm-degree} on the degree distribution $\vp$. 
For $d\sim \vp$ and $r\geq 1$, let $\sigma_r = \E(d^r) $. 
For a degree sequence $\vd_n=(d_1, \ldots, d_n)$, let 
$\sigma_r(\vd_n) := n^{-1} \sum_{i=1}^n d_i^r$, and $\mu^{\sss(n)} := \sigma_1(\vd_n)$.
We will work with \chh{deterministic degree sequences satisfying the following} regularity conditions, as this will be convenient for us.
\begin{ass}\label{ass:cm-degree-det}
	We have $\nu>1$, 
	$n\mu^{(n)}$ is even for all $n$,
	and $n^{-1}\sum_{i\in [n]}\delta_{\{d_i\}}$ converges, as $n\to\infty$, to the law associated to the pmf $\vp$.
	Further, there exists $N <\infty$ such that the sequence $\big(\vd_n;\, n\geq N\big)$ of degree sequences satisfies the following: 
	\begin{enumeratea}
		\item  There exists $B >0$ such that the \chh{maximum} degree  $\deg_{\max} := \max_{i\in [n]} d_i < B \log{n}$.
		\item 
		There exists $q > 0$ such that
		$|\sigma_r(\vd_n) - \sigma_r| \leq {(\log{n})^q}/{\sqrt{n}}$
		for $r=1, 2, 3$.
	\end{enumeratea}
\end{ass}
Under Condition \ref{ass:cm-degree} on the pmf $\vp$, degrees generated in an i.i.d. fashion from $\vp$ satisfy Condition \ref{ass:cm-degree-det} (in this case, $N$ is a finite random variable).
(Recall from Section \ref{sec:res-cm} our convention of adding an extra half-edge to a vertex with the maximum degree when the total degree is odd.)
As before, $\cC_i(t)$ denotes the $i$-th largest component in $\CM_n(t)$, and $f_i(t) \leq \sum_{v\in \cC_i(t)} d_v$ denotes the number of alive half-edges in $\cC_i(t)$. 
Let $s_1(t) := \mu\exp(-2t)$, $t\geq 0$, and recall the definition of $\bars_l(t)$ from \eqref{eqn:suscep-def-cm}.

\begin{lem}\label{lem:s1-conc}
	Under Condition \ref{ass:cm-degree-det}, for any $T> 0$,  there exists $C:=C(T)>0$ and $n_0(T)\geq 1$ such that for all $n > n_0(T)$ and $\eta \in [0,1/2)$,
	$\pr(\sup_{0\leq t\leq T}|\bars_1(t) - s_1(t)| > n^{-\eta})\leq \exp(-C n^{1-2\eta}).$
\end{lem}
\noindent {\bf Proof:} By construction, $\bars_1(\cdot)$ is a pure death process with jumps of size $-2/n $ occurring at rate $n\bars_1(t)$.
Since $\bars_1(0) = \mu^{\sss(n)}$, {by \cite{ethier-kurtz,kurtz-density}}, $\bars_1$ can be constructed as the unique solution of the equation
\begin{equation}
	\label{eqn:856}
	\bars_1(t):= \mu^{\sss(n)} - \frac{2}{n}\cdot Y\Big(n\int_0^t \bars_1(u) du\Big), \qquad t\geq 0,
\end{equation}
where $\big(Y(t);\, t\geq 0\big)$ is a rate one Poisson process. 
Let $\bar{Y}(\cdot) := Y(n\,\cdot)/n$. 
Analogous to $s_1$, define $s_1^*(t) := \mu^{\sss(n)} \exp(-2t)$. By Condition \ref{ass:cm-degree-det}{,} for large $n$,
$\sup_{t\geq 0}|s_1^*(t) - s_1(t)| \leq {\log^q{n}}/{\sqrt{n}}.$
Thus, it is enough to prove the lemma with $s_1^*(t)$ in place of \chhh{$s_1(t)$}. Note {that} $s_1^*(\cdot)$ satisfies
\begin{equation}
	\label{eqn:904}
	s_1^*(t):= \mu^{\sss(n)} - 2\int_0^t s_1^*(u) du, {\qquad t\geq 0.}
\end{equation}
Using \eqref{eqn:856}, \eqref{eqn:904}, and $\bars_1(\cdot)\leq \mu^{\sss(n)}$, we have, for any $t\in [0, T]$,
\begin{align*}
	|\bars_1(t) - s_1^*(t)|
	&\leq \sup_{0\leq \css{v}\leq \mu^{\sss(n)} T} \left|2 \bar{Y}(v) - 2v\right| + 2 \int_0^t |\bars_1(u) - s_1^*(u)| du.
\end{align*}
Now Gronwall's lemma \cite[page 498]{ethier-kurtz} yields
\[
\sup_{0\leq t\leq T} |\bars_1(t) - s_1^*(t)| 
\leq 
e^{2T} \sup_{0\leq t\leq \mu^{\sss(n)} T } \left|2 \bar{Y}(t) - 2t\right|. 
\]
Standard \chh{large deviation estimates for a Poisson process complete the proof.} \qed

The next two results describe an equivalence between the dynamic construction of $\CM$ and percolation on $\CM_n(\vd_n)$.  
The first result is obvious by construction; we omit its proof.

\begin{lem}\label{lem:cmn-t-k-2k-equiv}
	For any $0<t_1<t_2$ and $0\leq k_1\leq k_2\leq \sum_{i\in [n]}d_i/2$, conditional on the event
	$\big\{|E(\CM_n(t_i))|=k_i\, ,\, i=1, 2\big\}$, 
	$\big(\CM_n(t_1),\,\CM_n(t_2) \big)$ has the same distribution as $\big(\cQ_1,\, \cQ_2\big)$ constructed as follows:
	\begin{inparaenuma}
		\item 
		Choose a uniform subset of $2k_1$ half-edges from the set of the $\sum_{i\in [n]} d_i$ half-edges.
		Perform a uniform matching of these $2k_1$ half-edges to construct $\cQ_1$.
		\item 
		Conditional on the previous step, choose a uniform subset of $2k_2-2k_1$ half-edges from the set of the remaining $\sum_{i\in [n]} d_i-2k_1$ half-edges.
		Perform a uniform matching of these $2k_2-2k_1$ half-edges; this will result in a graph that contains $\cQ_1$ as a subgraph. 
		Declare this graph to be $\cQ_2$.
	\end{inparaenuma}
\end{lem}
\begin{prop}[\cite{dhara-hofstad-leeuwaarden-sen-cm-finite-third-moment}, Lemmas 8.1 and 8.2]\label{prop:perc-cm}
	For any $0<p_1<p_2<1$ and $0\leq k_1\leq k_2\leq \sum_{i\in [n]}d_i/2$, conditional on the event 
	$\big\{ |E(\Perc_n(p_i))|=k_i\, ,\, i=1, 2\big\}$, 
	$\big(\Perc_n(p_1),\,\Perc_n(p_2) \big)$ has the same distribution as $\big(\cQ_1,\, \cQ_2\big)$ as in Lemma \ref{lem:cmn-t-k-2k-equiv}.
\end{prop}

\subsection{Bounds in the subcritical regime}
\label{sec:proof-cm-max-diam}
In this section we prove Theorem \ref{thm:config-largest-comp-diam}. 
Let $\big(\CM_n^*(k);$ $\,k\in\bZ_{\geq 0}\big)$ be the embedded discrete time Markov chain in $\big(\CM_n(t);\, t\geq 0\big)$. 
Thus, $\CM_n^*(k)=\CM_n(\fT_k^{\sss(n)})$, where $\fT_k^{\sss(n)}$ denotes the time when the $k$-th (full) edge gets added in $\CM_n(\cdot)$.
Recall that the giant component emerges when $k\approx n\mu/(2\nu)$. We start by describing a result for the discrete step model in the barely subcritical regime. 
\begin{thm}\label{thm:CM-fixed-p}
	Fix $A_1, A_2, \kappa > 0$, and $\delta \in(0, 1/4)$. Let
	\begin{equation}
		\label{eqn:kn-inter-def}
		k_n^- := \frac{n\mu}{2\nu} - \frac{A_1 n}{{\log{n}}},\ \ 
		\text{ and }\ \
		k_n^+ := \frac{n\mu}{2\nu} - A_2n^{1-\delta} \, .
	\end{equation}
	Then there exist $\theta^*>0$ and $n^*\geq 1$ depending only on $A_1, A_2, \kappa, \delta$, and $(\vd_i,\, i\geq 1)$ such that for all $n\geq n^*$, we have, for all $k\in [k_n^-, k_n^+]$,
	\[
	\pr\big(\big\{ |\cC_1(\CM_n^*(k))|   \geq \theta^* (\log{n})/ a_n(k)^2\big\}
	\cup
	\big\{\diam(\CM_n^*(k)) \geq {\theta^* (\log{n})}/{a_n(k)}\big\} \big) 
	\leq n^{-\kappa} ,
	\]
	where $a_n(k) = \mu/(2\nu) - k/n$, and $\cC_1(\CM_n^*(k))$ denotes the largest component of $\CM_n^*(k)$.
\end{thm}

Let us first prove Theorem \ref{thm:config-largest-comp-diam} assuming Theorem \ref{thm:CM-fixed-p}. The proof of Theorem \ref{thm:CM-fixed-p} is given in the next section.

\noindent{\bf Proof of Theorem \ref{thm:config-largest-comp-diam}:} Recall the definition of $t_n$ from \eqref{eqn:tn-def}. 
We will show how Theorem \ref{thm:CM-fixed-p} implies the existence of $\theta_0 > 0$ such that
\begin{equation}
	\label{eqn:max-cm}
	\eps_n(\theta_0):=	
	\pr\left(\exists~ t \in [0,t_n] \mbox{ such that } |\cC_{1}(t)| \geq \theta_0 \log{^3 n}/{(t_c-t)^2}\right) \to 0\, .
\end{equation}
The proof for the diameter is similar.
As mentioned before, we will work with a deterministic degree sequence satisfying Condition \ref{ass:cm-degree-det}. 
Recall that $\sum_{i\in [n]} d_i = n\mu^{\sss(n)}$.
Fix $\eta\in (\delta, 1/2)$.
Using Lemma \ref{lem:s1-conc}, with probability at least $1-\exp(-C_1n^{1-2\eta})$, 
$| E(\CM_n(t_n)) |=(\sum_i d_i - n\bar{s}_1(t_n))/2 \in [\vE_{t_n}^-, \vE_{t_n}^+]$, 
where 
$\vE_{t_n}^{\pm}: = 2^{-1} (n\mu -ns_1(t_n)\pm 2n^{1-\eta} )$.
Note that
\begin{equation}
	\label{eqn:etn-pm-def}
	\vE_{t_n}^+\leq n\mu/(2\nu) - C_2n^{1-\delta}
\end{equation}
for an appropriate constant $C_2$. 
Recall that $\fT_i^{\sss(n)}$ denotes the time when the $i$-th edge is added in $\CM_n(\cdot)$. 
By definition, $n\bar{s}_1(\fT_i^{\sss(n)}) = \sum_{j\in [n]} d_j - 2i$. By Lemma \ref{lem:s1-conc} and for an appropriate choice of constant $C_3$,
\begin{equation}
	\label{eqn:fti-conc}
	\pr\big(
	\fT_i^{\sss(n)} \in \big[\, 
	s_1^{-1}(\mu-2i/n) \mp C_3 n^{-\eta}
	\,\big]
	\big) 
	\geq 
	1-\exp(-C_1 n^{1-2\eta}) \quad \forall 1\leq i\leq \vE_{t_n}^+. 
\end{equation}
Write $\vT_i^{\sss (n),\pm} = s_1^{-1}(\mu-2i/n) \pm C_3 n^{-\eta}$ for the end points of the above intervals. We now finish the proof of \eqref{eqn:max-cm}. 
Let $A_1=1$ and $A_2=C_2$, where $C_2$ is as in \eqref{eqn:etn-pm-def}. 
Let $k_n^{\pm}$ be as in \eqref{eqn:kn-inter-def} with this choice of $A_1$ and $A_2$.
Then
\begin{align}
	\eps_n(\theta) 
	\leq 
	\exp(-C_1 n^{1-2\eta})
	+
	\sum_{i=1}^{\vE_{t_n}^+} \pr(\fT_i^{\sss(n)} \notin [\vT_i^{\sss(n),-},\vT_i^{\sss(n),+}]) 
	+ 
	\sum_{i=1}^{\vE_{t_n}^+} 
	\pr\Big(|\cC_1(\CM_n^*(i))|\geq \frac{\theta \log^3{n}}{(t_c - \vT_i^{\sss(n),-})^2}\Big) 
	\notag\\
	\leq 
	C_4 n \exp(-C_1 n^{1-2\eta}) 
	+ 
	\sum_{i=1}^{k_n^-} \epsilon_n^{\sss(i)}(\theta) 
	+ 
	\sum_{i=k_n^- + 1}^{k_n^+}\epsilon_n^{\sss(i)}(\theta)\, ,
	\hskip70pt
	\label{eqn:65}
\end{align}
where 
$\epsilon_n^{\sss(i)}(\theta) 
:= 
\pr\big(|\cC_1(\CM_n^*(i))|\geq \theta \log^3{n}/(t_c - \vT_i^{\sss(n),-})^2\big)$. 
Now, for all large $n$, for any $i\in[1, k_n^-]$, 
\begin{align}\label{eqn:344}
	\epsilon_n^{\sss(i)}(\theta) 
	\leq
	\pr\big(|\cC_1(\CM_n^*(k_n^{-}))|\geq \theta \log{n}/(t_c a_n(k_n^-))^2\big)\, .
\end{align}
Further, a direct computation shows that there exist constants $C_5$ and $C_6$ such that for $n$ large enough, for any $i\in [k_n^{-}, k_n^+]$,  
$
\big| 
\big(t_c-\vT_i^{\sss(n),-}\big)-C_5 a_n(i)\big| 
\leq 
C_6 \big(a_n^2(i) + n^{-\eta}\big)$,
which in particular implies that for all large $n$, for any $i\in [k_n^-, k_n^+]$,
\begin{align}\label{eqn:767}
	\epsilon_n^{\sss(i)}(\theta) 
	\leq
	\pr\big(|\cC_1(\CM_n^*(i))|\geq \theta \log^3{n}/(2C_5 a_n(i))^2\big)\, .
\end{align}
Thus, \eqref{eqn:max-cm} follows upon combining 
\eqref{eqn:65}, \eqref{eqn:344}, and \eqref{eqn:767}, 
and then applying Theorem \ref{thm:CM-fixed-p} with $\kappa =2$.
\qed

\subsubsection{Proof of Theorem \ref{thm:CM-fixed-p}}
All constants in this section will depend only on 
$A_1, A_2, \kappa, \delta$, and $(\vd_i ;\, i\geq 1)$.
Simililarly, statements that hold for all large $n$ will be true for all $n$ bigger than a threshold that depends only on 
$A_1, A_2, \kappa, \delta$, and $\big(\vd_i;\, i\geq 1\big)$.
Inequalities involving $k$ will be true uniformly for $k\in [k_n^-,k_n^+]$. 
We will make use of the following result:
\begin{thm}[\cite{dubhashi1996negative}, Theorem 10]\label{thm:negative-association}
	Let $y_1,\ldots, y_n\in\bR$, and let $\pi$ be a uniform permutation on $[n]$.
	Then $y_{\pi(i)}$, $i=1,\ldots, n$, are negatively associated.
	Consequently, if a sample (of size $r$, say) is drawn without replacement from a collection of $N$ objects, the random variables $\ind_{\{i\text{-th object is selected in the sample}\}}$, $1\leq i\leq N$, are negatively associated.
\end{thm}

The random graph $\CM_n^*(k)$ can be constructed in two steps: 
(a) sample $2k$ half edges from the set of $\sum_{i\in[n]}d_i$ many half-edges without replacement; 
(b) perform a uniform matching of the selected half edges.  
Let $\td_i$ denote the number of half-edges selected in step (a) above that are attached to vertex $i$, and let $\vd^*_n := (\td_i,\, i\in [n])$. 
(Thus, $\sum_{i\in [n]}d^*_i=2k$.)
Let
\[
\fh = \fh(k) :=\frac{2k}{n\mu} = \frac{1}{\nu} - \frac{2}{\mu}a_n(k)\, .\]
We suppress dependence on $k$ unless required and write, e.g., $a$ for $a_n(k)$.

\begin{lem}\label{lem:per-degree-prop}
	There exists $C_1>0$ such that for all $n\geq 2$,
	{\bf \upshape (a)} 
	$|2k - \fh \sum_{i=1}^n d_i | \leq C_1\log^{q}{n}\cdot\sqrt{n}$, and 
	{\bf \upshape (b)} 
	with probability at least $1-n^{-\kappa-1}$,  
	$|\sum_{i=1}^n (\td_i)_r  - \fh^r \sum_{i=1}^n (d_i)_r | \leq C_1\log^{q+3}{n}\cdot\sqrt{n}$
	for $r=2, 3$, where $(m)_r=m(m-1)\ldots(m-r+1)$ for $m\in\bN$.
\end{lem}
\noindent {\bf Proof:} 
The claim in {\bf \upshape (a)} follows from Condition~\ref{ass:cm-degree-det}.
Next, by Theorem \ref{thm:negative-association}, for every $r\geq 1$, the random variables $(\td_i)_r\,, i\in[n]$, are negatively associated.
Using the expressions for factorial moments of a hypergeometric distribution and Condition \ref{ass:cm-degree-det}, one directly checks that
$\E[\sum_{i=1}^n (\td_i)_r] = \sum_{i=1}^n (d_i)_r \cdot \fh^r + O(\log^q n\cdot\sqrt{n})$ for $r=2, 3$. 
By \cite[Proposition 5]{dubhashi1998balls}, the bounded difference inequality applies to a sum of negatively associated random variables. 
Now using the bounds $d_i^*\leq d_i$ together with Condition~\ref{ass:cm-degree-det} completes the proof of {\bf \upshape (b)}.
\qed  

\vskip3pt

We will henceforth work conditional on $\vd^*_n$, and further, on the event 
$E_n:=\{|\sum_{i=1}^n (\td_i)_r  - \fh^r\sum_{i=1}^n (d_i)_r | \leq C_1\log^{q+3}{n}\cdot\sqrt{n}\text{ for }r=1,2,3\}$. 
By Lemma \ref{lem:per-degree-prop}, $\pr(E_n)\geq 1-n^{-\kappa-1}$.
Now, Condition \ref{ass:cm-degree-det} implies that there exist $C_2,C_3,C_4>0$ such that for all large $n$, on $E_n$,
\begin{gather}
	\Big|
	\sum_{i=1}^n \td_i(\td_i-1)/(2k) - (1-C_2a)
	\Big|
	\leq
	C_4(\log{n})^{q+3}/{\sqrt{n}} ,
	\ \ \text{ and }\label{eqn:vp-asymp}\\
	2C_3 < \sum_{i=1}^n (\td_i-2)^2 \td_i/(2k)  < C_4 
	\label{eqn:c1-c2-boundd-2}\, .
\end{gather}

\noindent{\bf Analysis of the largest component:}
Write 
$\pr_1(\cdot)=\pr(\cdot\, | \vd^*_n)$, 
$\bE_1(\cdot)=\bE(\cdot\, | \vd^*_n)$,
and let $\cC^*(i)$ be the component of the vertex $i\in [n]$ {in $\CM_n^*(k)$}.
\begin{lem}\label{lem:comp-size-bias}
	Sample $v(1)\in [n]$ such that $\pr_1(v(1)=i) \propto \td_i\, ,\, i\in [n]$. 
	Then there exists $\theta>0$ such that on the event $E_n$,
	$\pr_1\left(|\cC^*({v(1)})| \geq {\theta (\log{n})}/{a^2}\right) \leq n^{-\kappa-2}$.
\end{lem}
This lemma implies the bound on $|\cC_1(\CM_n^*(k))|$ in Theorem \ref{thm:CM-fixed-p}, since 
\begin{align}
	&\pr_1\left(|\cC_1(\CM_n^*(k))| \geq {\theta (\log{n})}/{a^2}\right) 
	\leq 
	\Big(\sum_{i=1}^n \td_i\Big) \sum_{i=1}^n \frac{\td_i}{\sum_{j=1}^n \td_j} 
	\pr_1\left(|\cC^*(i)|\geq {\theta (\log{n})}/{a^2}\right) \notag\\  
	&\hskip30pt
	= 2k\cdot \pr_1\left(|\cC^*({v(1)})| \geq {\theta (\log{n})}/{a^2}\right) 
	\leq
	(n\mu/\nu)\cdot \pr_1\left(|\cC^*({v(1)})| \geq {\theta (\log{n})}/{a^2}\right)
	.\label{eqn:999}
\end{align}

\noindent {\bf Proof of Lemma \ref{lem:comp-size-bias}:}
\chhh{Throughout this proof, we will work conditional on the degree sequence $\vd_n^*$.} 
Declare all $2k$ half-edges to be {\it available}.
Sample $v(1)$ as in the statement of Lemma \ref{lem:comp-size-bias}, and declare $v(1)$ and its $\td_{v(1)}$ many half-edges to be \emph{current}. 
(A current half-edge is still available.)
Continue exploring the component $\cC^*(v(1))$ in a breadth-first manner, where at each step, a current half-edge is paired with another available half-edge chosen uniformly, and both these half-edges are declared {\it unavailable}.
In this exploration, two possibilities arise at each step:
\begin{inparaenuma}
	\item A current half-edge may select an available half-edge that is not current.
	In this case, declare the newly found vertex and all its \emph{remaining} half-edges to be \emph{current}.
	\item A current half-edge may select one of the other current half-edges. In this case, no new vertex is added to the component in this step.
\end{inparaenuma}

We stop when there are no current half-edges (at which point we have completed exploring $\cC^*(v(1))$).
Let $v(j)$ denote the $j$-th new vertex found in the above exploration.
By construction, $v(1), v(2),\ldots$ appear in a size-biased order using $\vd^*_n$ as vertex weights. 
Let $R_n(i)$ be the number of current half-edges at step $i$, with $R_n(1)= \td_{v(1)}$. 
Let $\current(i)\leq i$ be the number of current vertices found by the $i$-th step. Then
\begin{equation}
	\label{eqn:sn-i-transition}
	R_n(i+1)=\left\{\begin{array}{ll}
		R_n(i)+ \td_{v(\current(i)+1)}-2\, , & \mbox{ under the event in (a), } \\
		R_n(i) -2\, , & \mbox{ under the event in (b).}
	\end{array}  \right.
\end{equation}
If $H_n$ is the first time the above walk hits zero, then $|\cC^*({v(1)})|\leq H_n$. 
Write $(\tR_n(i),\, i\geq 1)$ {for} the walk that ignores transitions modulated by the event in (b), i.e., $\tR_n(i) = 2+ \sum_{j=1}^i (\td_{v(j)}-2)$. Let $\tH_n$ be the corresponding hitting time of zero. Since $H_n\leq \tH_n$, {it is enough to show} that we can choose $\theta_1$ so that for all large $n$, on the event $E_n$, 
\begin{equation}
	\label{eqn:th-hit}
	\pr_1\big(\tH_n > {\theta_1 (\log{n})}/{a^2} \big) \leq n^{-\kappa-2}.
\end{equation}

Let $m_n = m_n(\theta_1):= \theta_1 \log{n}/a^2 $, and $m_n^\prime := m_n/2$. 
Following \cite{aldous1997brownian}, we use  a reformulation of the problem using an artificial time parameter. Let $(\xi_i,\, i\in [n])$ be an independent sequence of random variables with \chhh{$\xi_i\sim\mbox{Exponential}(\td_i/(2k))$}. 
If $\xi_{\tilde v(1)}< \xi_{\tilde v(2)}< \ldots< \xi_{\tilde v(n)}$, then $\big(\tilde v(i) ,\, i\in [n]\big)$ is a size-biased permutation of $[n]$ using $\vd^*_n$ as weights.
For $t\geq 0$, let $Q_n(t) :=2+ \sum_{i=1}^n (\td_i-2) \ind\set{\xi_i \leq t}$, and
$\Upsilon_n(t):=\sum_{i=1}^n \ind\set{\xi_i \leq t}$.
Further, define $T_m :=\inf\big\{t>0 :  \Upsilon_n(t) =m\big\}$, $m\geq 1$.
Then
\begin{align}
	\pr_1(\tH_n > m_n) 
	&\leq 
	\pr_1(Q_n(t) > 0 \mbox{ for all } 0< t \leq m_n^\prime, T_{m_n} > m_n^\prime) 
	+ 
	\pr_1(T_{m_n} < m_n^\prime)\notag\\
	&\leq 
	\pr_1(Q_n(m_n^\prime) > 0) 
	+ 
	\pr_1(T_{m_n} < m_n^\prime).   \label{eqn:sn-tmn}
\end{align}
To deal with the first term, note that 
$\E_1\big[Q_n(m_n^\prime)\big] 
= 
2+\sum_{i=1}^n (\td_i-2) (1-\exp(-m_n^\prime \td_i/(2k)))$. 
Expanding and using \eqref{eqn:vp-asymp}  gives, for all large $n$, on the event $E_n$,
\begin{align}
	\E_1\big[Q_n(m_n^\prime)\big] = -(1+o(1))C_2 m_n^\prime a. \label{eqn:sn-expec-form}
\end{align}
Now, 
$Q_n(m_n^\prime)-2= \sum_{i=1}^n \chh{(\td_i-2)} \ind\set{\xi_i\leq m_n^\prime}=:\sum_i Y_i $ 
say, is a sum of independent random variables with  $|Y_i|\leq B\log{n}$. 
By \eqref{eqn:c1-c2-boundd-2}, for all large $n$, on the event $E_n$, 
$C_3 m_n^\prime < \sum_i \E_1(Y_i^2) < C_4 m_n^\prime$.
By Bennett's inequality \cite[Theorem 2.9]{boucheron2013concentration},
for all large $n$, on $E_n$, 
\begin{align}
	\pr_1\big(Q_n(m_n^\prime) > 0\big)
	=\pr_1\big(Q_n(m_n^\prime) - \bE_1 Q_n(m_n^\prime)> - \bE_1 Q_n(m_n^\prime)\big)
	\leq \exp\big(-{C_3 m_n' h(u_n)}/{(B\log{n})^2}\big),\label{eqn:bnnn}
\end{align}
where $h(u)= (1+u)\log{(1+u)}-u$, and $u_n = B \log{n}\cdot|\E_1(Q_n(m_n^\prime))|/(C_4m_n^\prime)$. Using \eqref{eqn:sn-expec-form} and the fact that $a\leq A_1/\log{n}$, we get, for large $n$, on $E_n$, $u_n\leq C_5:=2BC_2A_1/C_4$.
Now, there exists $C_6>0$ such that $h(u)\geq C_6 u^2$ for $u\in [0,C_5]$.  
Thus, for all large $n$, on $E_n$, 
\begin{align}\label{eqn:149}
	\pr_1\big(Q_n(m_n^\prime) > 0\big) \leq \exp\big(-C_7 a^2 m_n^\prime\big) = \exp(-C_8\theta_1 \log{n}).
\end{align}

Now consider the second term in \eqref{eqn:sn-tmn}. 
If $T_{m_n} < m_n^\prime$, then $\Upsilon_n(m_n/2) \geq m_n$. 
Note that $\Upsilon_n(m_n/2)$ is a sum of independent indicator random variables. 
Further,
$\E_1(\Upsilon_n(m_n/2)) = (1+o(1)){m_n}/{2}$.
By \cite[Theorem 6.12]{boucheron2013concentration}, for all large $n$,
$\pr_1(\Upsilon_n(m_n/2) \geq m_n) 
\leq 
\exp\left(-{3m_n}/{32}\right) \leq \exp\left(-3\theta_1\log^3 n/(32 A_1^2)\right)$. 
Combining this with \eqref{eqn:149} and \eqref{eqn:sn-tmn} and taking $\theta_1=(\kappa+3)/C_8$ yields \eqref{eqn:th-hit}.
\qed

\vskip5pt

\noindent{\bf Analysis of the diameter:} 
Let $\theta_2=2(\kappa+3)/C_2$, where $C_2$ is as in \eqref{eqn:vp-asymp}.
We will show that for all large $n$, on the event $E_n$,
\[
\pr_1\big(\diam\big(\cC^*(v(1))\big) \geq 2\theta_2 \log{n}/a \big)
\leq n^{-\kappa-2},
\]
which together with an argument similar to \eqref{eqn:999} will complete the proof.
Using Lemma \ref{lem:comp-size-bias} with $\kappa$ replaced by $\kappa+2$, we can choose $\theta_3$ so that on the event $E_n$,
\begin{equation}
	\label{eqn:cv-bd-ka}
	\pr_1(|\cC^*({v(1)})| > \theta_3 \log{n}/a^2)\leq n^{-\kappa-4}.
\end{equation}
As in the proof of Lemma \ref{lem:comp-size-bias}, we construct $\cC^*({v(1)})$ via a breadth-first exploration. 
For $r\geq 0$, let $\cF_r$ be the sigma-field generated by $\vd_n^*$ and this exploration process up to discovering the last vertex in the the $r$-th generation of the breadth-first tree.
Further, let $G_r$ be the number of current half-edges
attached to the vertices in generation $r$ right after the last vertex in generation $r$ has been discovered.
Then $(G_r,\, r\geq 0)$ is adapted to the filtration $(\cF_r,\, r\geq0)$. 
Suppose we show that there exists a constant $C_9> 0$ such that for all large $n$, for each $r\leq \theta_2\log{n}/a$, \ch{there exists a} set $\fA_r\in \cF_r$ such that 
\begin{enumeratei}
	\item 
	on $\fA_r\cap E_n$, 
	$\E(G_{r+1}|\cF_r) 
	\leq (1-a C_2)\left(1+C_9(\log{n})^{q+3}/\sqrt{n}\right) G_r
	$, and
	\item 
	on the event $E_n$, $\pr_1(\fA_r^c)\leq n^{-\kappa-4}$. (On $\fA_r^c$  {we have the} trivial bound $\E(G_{r+1}|\cF_r)\leq 2k$.)
\end{enumeratei}
Then this results in the following recursion: 
on the event $E_n$,
\[
\E_1(G_{r+1})
\leq 
(1-a C_2)\left(1+C_9(\log{n})^{q+3}/\sqrt{n}\right)\E_1(G_r) + 2k/n^{\kappa+4} .
\]
Iterating this starting at $r_n=\theta_2\log{n}/a$ and using the fact $a\geq A_2n^{-\delta}$,
a routine calculation yields, for all large $n$, on the event $E_n$,
\begin{align} \label{eqn:g-fin-bound}
	\E_1(G_{r_n}) 
	\leq 
	e^{-C_2 ar_n/2} +2k r_n/n^{\kappa+4}  
	\leq 
	n^{-\kappa-2}\, .
\end{align}
Since
$\pr_1(\diam(\cC^*(v(1))) \geq 2 r_n) \leq \E_1(G_{r_n})$, using \eqref{eqn:g-fin-bound} completes the proof.

Now we show the existence of the sets $\fA_r$. 
Let $\cA_r$ be the set of vertices reached by the exploration by \chh{generation} $r$, and $\cD_r = [n]\setminus \cA_r$. 
Let
$\fA_r = \{\sum_{v\in\cA_r} \td_v \leq \theta_3 B (\log{n})^2/a^2  \}$ with 
$\theta_3$ as in \eqref{eqn:cv-bd-ka}, and $B$ as in Condition \ref{ass:cm-degree-det}.
Then, by \eqref{eqn:cv-bd-ka} and Condition \ref{ass:cm-degree-det}~(a), on the event $E_n$,
$\pr_1(\fA_r^c) \leq n^{-\kappa-4}$.  
Recall that in the breadth-first exploration, new vertices are added in a size-biased manner. 
We now state an elementary lemma without proof.

\begin{lem}\label{lem:size-bias-dom}
	Let $(w_v,\, v\in [\ell])$ be a sequence of positive weights, and suppose 
	$(v(1),\ldots, v(\ell))$ is a size-biased permutation of $[\ell]$ using this weight sequence.
	Then $w_{v(1)}$ stochastically dominates $w_{v(i)}$ for $i=2,\ldots, \ell$. 
	Consequently, $\E(w_{v(1)})\geq \E(w_{v(i)})$, $i=2,\ldots, \ell$. 
\end{lem}
Using the above lemma, \chh{we get} 
$\E(G_{r+1}|\cF_r) \leq \alpha_r G_r$,
where $\alpha_r =\E(d_{v_r(1)}-1|\cF_r)$, {and} where $v_r(1)$ is a vertex selected via size biased sampling from $\cD_r$. 
Further,
\[
\E(d_{v_r(1)}-1|\cF_r) 
= 
\frac{\sum_{i=1}^n \td_i(\td_i-1) - \sum_{v\in \chh{\cA_r}} \td_v(\td_v-1)}{\sum_{i=1}^n \td_i -\sum_{v\in \chh{\cA_r}} \td_v}
\leq 
\frac{\sum_{i=1}^n \td_i(\td_i-1)}{2k-\sum_{v\in \chh{\cA_r}} \td_v} \, .
\]
Now, \eqref{eqn:vp-asymp} yields
$\alpha_r \leq (1-aC_2) (1+C_9(\log{n})^{q+3}/\sqrt{n})$ on the event $\fA_r\cap E_n$
for {an} appropriate constant $C_9$. \chh{This completes the proof of Theorem \ref{thm:CM-fixed-p}.}

\subsection{Properties of susceptibility functions}\label{sec:cm-prop-entrance}
Here we prove Theorem \ref{thm:config-bare-sibcrit}. 
Let $\boldsymbol{\cF}:=\big(\cF_t;\, t\geq 0\big)$ be the natural filtration of $\big(\CM_n(t);\, t\geq 0\big)$.  
For an $\boldsymbol{\cF}$-adapted semimartingale $J(t)$ of the form
\begin{equation}\label{eqn:semimart}
	d J(t) = \alpha(t) dt + d M(t) \ \ \text{ with }\ \
	\langle M, M \rangle(t) = \int_{0}^t \gamma(s) ds,
\end{equation}
write $\vdd(J)(t):= \alpha(t)$, $ \vv (J)(t) := \gamma(t)$, and $\vM(J)(t) {:=} M(t)$. 
We will write
$\Delta J(t) = J(t + \Delta t) - J(t)$. 
For nonenegative stochastic processes $\xi_n(\cdot)$, $n\geq 1$, and $T>0$, a term of the form $O_{T}(\xi_n(t))$ will represent stochastic processes $\big(\eps_n(t);\, 0\leq t\leq T\big)$, $n\geq 1$, for which there exists a deterministic constant $C_1> 0$ independent of $n$ such that whp, 
$|\eps_n(t)| \leq C_1\xi_n(t)$ for all $0\leq t\leq T$.
Note that by Lemma \ref{lem:s1-conc}, whp, 
$(2\mu)^{-1} \leq 1/\bars_1(t) \leq 2e^{2t_c}/\mu$ for all $t\leq t_c$. 
Thus, in the rest of this section, we will drop the factor $1/\bars_1(t)$ in a term of the form $O_{t_c}(\cdot)$.

Recall that $f_i(t)$ denotes the number of alive half-edges in $\cC_i(t)$, \chh{the $i$-th largest component at time $t$}. \chh{Note that $\max_{i\geq 1} f_i(t)\leq (B \log{n})|\cC_1(t)|$ by Condition \ref{ass:cm-degree-det}~(a); to simplify notation define} 
$I(t):=\chh{(B \log{n})}|\cC_1(t)|$. 
Recall {the} functions $\bars_l, \barg$, and $\barcd$ in \eqref{eqn:suscep-def-cm}. 
Define  
\begin{equation}
	\label{eqn:bars2-star-def}
	\bars_2^{\star}(t):= n^{-1}\chh{\sum_{i\geq 1}} |\cC_i(t)|^2 \, , \ \ \  t\geq 0.
\end{equation}

\begin{lem}\label{lem:ds-dg-dd}
	The processes $\bars_2, \bars_3,\barg$, $\barcd$, and $\bars_2^{\star}$ are $\boldsymbol{\cF}$-semimartingales satisfying the following:
	\begin{enumeratea}
		\item \label{item:s2} 
		$\vdd(\bars_2)(t)= F_2^{\vs}(\bars_1(t), \bars_2(t)) + O_{t_c}(I^2(t)\bars_2(t)/n) $, where
		\begin{equation*}
			\label{eqn:f2-def}
			F_2^{\vs}(s_1, s_2):= s_1^{-1}\left[2 s_2^2 + 4 s_1^2 - 8 s_2 s_1\right],
		\end{equation*}
		\item $\vdd(\bars_3)(t) = F_3^{\vs}(\bars_1(t), \bars_2(t), \bars_3(t)) + O_{t_c}(I^3(t)\bars_2(t)/n)$, where
		\begin{equation*}
			\label{eqn:f3-def}
			F_3^{\vs}(s_1,s_2,s_3):= s_1^{-1}{s_2}(6 s_3 - 12 s_2) + (24 s_2-12 s_3-8 s_1),
		\end{equation*}
		\item 
		$\vdd(\barg)(t) = F^{\vg}(\bars_1(t), \bars_2(t), \barg(t)) + O_{t_c}(I^2(t)\bars_2(t)/n)$, where
		\[	
		F^{\vg}(s_1,s_2,g):={s_1^{-1}}\left[2 g s_2 - 4g s_1\right],
		\]
		\item \label{item:bcd} $\vdd(\barcd)(t)=F^{\vD}(\bars_1(t),\bars_2(t),\barcd(t)) + O_{t_c}(\diam_{\max}(t) I^2(t)\bars_2(t)/n)$, where
		\begin{equation*}\label{eqn:fd-def}
			F^{\vD}(s_1,s_2,d):= {s_1^{-1}}\left[ \chh{4 d s_2 +  2 s_2^2 - 8 d s_1 - 4 s_2 s_1 - 2 s_1^2 } \right],
		\end{equation*}
		\item \label{item:s2-star-mean} 
		$\vdd(\bars_2^{\star})(t) = F_2^{\star}(\bars_1(t),\barg(t))+O_{t_c}(I^2(t)\bars_2(t)/n)$, where
		\[
		F_2^{\star}(s_1, g):= {2 {s_1^{-1}} g^2}\, ,
		\]
		\item $\vv(\bars_2)(t)= O_{t_c}(I^2(t)\bars_2^2(t)/n)$,
		\item \label{item:vv-d} 
		$\vv(\barcd)(t) = O_{t_c}(\bars_2^2\css{(t)} I^2(t) (\diam_{\max}(t))^2/n)$, and
		\item\label{item:s2-star-variance} 
		$\vv(\bars_2^\star)(t)=O_{t_c}\big(I^2(t) \barg^2(t)/n\big)$.
	\end{enumeratea}
\end{lem}
\noindent{\bf Proof:} We prove the assertions in
\eqref{item:bcd} and \eqref{item:vv-d} above. {The} remaining results follow in an identical fashion. 
For any half-edge $e$ alive at time $t$ in a component $\cC(t)$, let 
$\cD(e) := \sum_{e^\prime\in \cC(t), e^\prime \text{alive}} d_{\cC(t)}(e,{e^\prime})$. 
Fix $t>0$, two components $\cC_i(t)\neq \cC_j(t)$, and two alive half-edges $u_0\in \cC_i(t)$, $v_0\in \cC_j(t)$. 
The Poisson clock at the half-edge $u_0$ rings at rate one, and then $u_0$ forms a full edge by connecting to the half-edge $v_0$ with probability 
$(n\bars_1(t)-1)^{-1}=(n\bars_1(t))^{-1}+O_{t_c}(n^{-2})$, 
where the last step uses Lemma \ref{lem:s1-conc}.
Writing $\cC_i$ (resp. $\cC_j$) for $\cC_i(t)$ (resp. $\cC_j(t)$), 
the change in $\barcd_1(t)$ under this event satisfies
\begin{align*}
	&
	n (\Delta \barcd(t))  =
	2 \sum_{\substack{e \in \cC_i,\\ e\neq u_0}} \sum_{\substack{e^{\prime} \in \cC_j,\\ e^{\prime}\neq v_0}} 
	\big(d_{\cC_i}(e,u_0) + d_{\cC_j}(e^{\prime},v_0) + 1 \big) 
	-2\sum_{\substack{e \in \cC_i}} d_{\cC_i}(u_0,e) - 2\sum_{\substack{e^{\prime} \in \cC_j}} d_{\cC_j}(v_0,e^{\prime})  \\
	&
	=2\left[  \cD({u_0}) f_j(t) + f_i(t) \cD({v_0}) +  f_i(t) f_j(t) - \cD({u_0}) -f_i(t) - \cD({v_0}) - f_j(t) + 1  \right] - 2\cD({u_0}) - 2\cD({v_0}).
\end{align*}
If an alive half-edge  decides to connect to another alive half-edge in the same component, then the change can be bounded as
$n\Delta \barcd(t)\leq 2 f_i(t)(f_i(t)-1) \chh{\diam_{\max}}(t).$
Summing over all \chh{pairs of components and all pairs of} alive half-edges \chh{in these components} gives
\begin{align*}
	\vdd(\barcd)(t) 
	&= 
	F^{\vD}(\bars_1(t),\bars_2(t),\barcd(t)) 
	+O_{t_c}\Big(
	\frac{1}{n^2}F^{\vD}(\bars_1(t),\bars_2(t),\barcd(t)) 
	+
	\frac{\diam_{\max}(t)\sum_i f_i^4(t)}{n(n\bars_1(t)-1)}\Big)\\
	&=
	F^{\vD}(\bars_1(t),\bars_2(t),\barcd(t)) 
	+ 
	O_{t_c}\left(\chh{\diam_{\max}(t)}{I^2(t) \bars_2(t) }/n\right),
\end{align*}
where the last step uses the relation $\max_i f_i(t)\leq I(t)$. This proves \eqref{item:bcd}. 

To prove \eqref{item:vv-d}, note that if an alive half-edge in $\cC_i(t)$ gets paired with an alive half-edge in $\cC_j(t)$, then the change in $\barcd$ can be bounded by
$|\Delta \barcd(t)| \leq {12 f_i(t)f_j(t) \diam_{\max}(t)}/{n}.$
Thus,
\[
\vv(\barcd)(t) 
\leq 
\sum_{i,j} \Big(\frac{f_i(t)f_j(t)}{n\bars_1(t)-1} \Big)
\frac{(12 f_i(t)f_j(t) \diam_{\max}(t))^2}{n^2} 
\leq 
\frac{144 \bars_3^2(t) (\diam_{\max}(t))^2}{n\bars_1(t)-1} .
\]
Using \chh{$\bars_3(t) \leq \bars_2(t) I(t)$} completes the proof.
\qed

By Lemma \ref{lem:s1-conc} we {know that} for any $T>0$, $\sup_{0\leq t\leq T}|\bars_1(t) - \mu\exp(-2t)| = O_P(n^{-\eta})$ for any $\eta < 1/2$. 
Analogously, we would like to say that a function of the objects considered in Lemma \ref{lem:ds-dg-dd} is approximated by the same function applied to the solutions of certain ODEs. 
We first determine the boundary values.  
We have, $\barcd(0) = 0$.  
By Condition \ref{ass:cm-degree-det}, for all large $n$,
$|\bars_2(0) - (\nu+1)\mu|\leq {\log^q{n}}/{\sqrt{n}}$, and  $|\bars_3(0) - (\beta + \mu(3\nu+1))|\leq {\log^q{n}}/{\sqrt{n}}$.
On {the} interval $[0,t_c)$, define 
\begin{equation}\label{eqn:s2-form}
	s_2(t) = \frac{\mu e^{-2t}\left(2\nu - (\nu-1) e^{2t}\right)}{\nu-e^{2t}(\nu-1)}
	\ \ \text{ and }\ \ 
	s_3(t) = \frac{-\beta + e_3(t)}{[-\nu+(\nu-1)e^{2t}]^3}\, ,
\end{equation}
where
$e_3(t) = -4\nu^3 \mu-9\nu^2\mu e^{2t}+9\nu^3\mu e^{2t} -6\nu\mu e^{4t}+ 12 \nu^2\mu e^{4t}
- 6 \nu^3 \mu e^{4t} -\mu e^{6t} + 3\mu \nu e^{6t} -3\nu^2 \mu e^{6t} + \nu^3\mu e^{6t}$. Note {that}
$e_3(t_c) = 0$.  Define
\begin{equation}\label{eqn:gt-form}
	g(t)= \frac{\mu}{\nu-(\nu-1)e^{2t}} 
	\ \ \text{ and }\ \ 
	\cD(t):= \frac{\nu^2\mu (1-e^{-2t})}{(\nu - (\nu-1)e^{2t})^2}\, ,\ \ 0\leq t<t_c.
\end{equation}
Recall the functions $F_2^{\vs},F_3^{\vs},F^{\vg}, F^{\vD}$ from Lemma \ref{lem:ds-dg-dd}. The following result can be checked by a direct calculation; we omit the proof. 
\begin{lem}\label{lem:diff-eq-solution}
	Writing $s_1(t) = \mu\exp(-2t)$, the functions $s_2,s_3,g$, and $\cD$ are the unique solutions on $[0,t_c)$ of the {ODEs} 
	$s_2^\prime(t) = F_2^{\vs}(s_1(t), s_2(t)) $, 
	$s_3^\prime(t) = F_3^{\vs}(s_1(t),s_2(t),s_3(t))$, 
	$g^\prime(t) = F^{\vg}(s_1(t),s_2(t),g(t))$, and 
	$\cD^\prime(t)= F^{\vD}(s_1(t),s_2(t),\cD(t))$ 
	with respective boundary conditions $s_2(0)= (\nu+1)\mu$, $s_3(0) = \beta + \mu(3\nu+1)$, $g(0) =\mu $, and $\cD(0) = 0$.  
	Further, replacing $\bars_2, \bars_3,\barg$, and $\barcd$ in Theorem \ref{thm:config-bare-sibcrit} with $s_2,s_3,g$, and $\cD$, the  assertions of the theorem, namely \eqref{eqn:s2-asymp-cm}, \eqref{eqn:s3-asymp-cm}, and \eqref{eqn:d1-g-asymp-cm} are satisfied for all $\delta < 1/3$.
\end{lem}
Define the processes
$Y(\cdot) := {1}/{\bars_2(\cdot)}$, $Z(\cdot)= {\bars_3(\cdot)}/{(\bars_2(\cdot))^3}$, $ U(\cdot) = {\barg(\cdot)}/{\bars_2(\cdot)}$, and $ V(\cdot):={\barcd(\cdot)}/{(\bars_2(\cdot))^2}$.
Let $y(\cdot),z(\cdot), u(\cdot)$ and $v(\cdot)$ be the corresponding deterministic analogues, i.e., $y(\cdot) = 1/s_2(\cdot)$, $z(\cdot) = s_3(\cdot)/(s_2(\cdot))^3 $ etc. It is easy to check that as $t\uparrow t_c$,
\begin{equation}
	\label{eqn:yt-zt-near-tc}
	y(t) = \frac{2\nu}{\mu(\nu-1)}(t_c-t)(1+O(t_c-t) ), \qquad z(t)= \frac{\beta (1+O(t_c-t))}{\mu^3(\nu-1)^3},
\end{equation}
\begin{equation}\label{eqn:ut-vt-near-tc}
	u(t)=\frac{(1+O(t_c-t))}{\nu-1}\, ,
	\ \ \text{ and }\ \
	v(t) = \frac{\nu(1+O(t_c-t))}{\mu (\nu-1)^2}\, .
\end{equation}
The following result completes the proof of Theorem \ref{thm:config-bare-sibcrit}.
\begin{prop}\label{prop:suscep-close-det}
	Fix  $\delta \in (1/6,1/5)$ and let $t_n$ as in \eqref{eqn:tn-def}. Then
	\begin{gather}
		\sup_{0\leq t\leq t_n}\left|Y(t) - y(t)\right| =o_P\big(\log^5 n/ n^{(1-\delta)/2}\big) =o _P(n^{-1/3})\, ,
		\ \  \text{ and }\ \ 
		\label{eqn:yt-close}\\
		\sup_{0\leq t\leq t_n} \max\{ |Z(t)-z(t)|, |U(t)- u(t)|, |V(t) - v(t)|\} \weakc 0\, .\label{eqn:89}
	\end{gather}
\end{prop}

We will prove \eqref{eqn:yt-close}; the assertion in \eqref{eqn:89} follows similarly.  
(See also the proof of \cite[Theorem 3.2 (ii)]{bhamidi2014augmented}, where the same technique is employed to prove analogous results for bounded-size rules.)
The main tool is the following lemma.

\begin{lem}[{\cite[Lemma 6.10]{bhamidi2014augmented}}]\label{lem:conv-pure-jump}
	Suppose $\bD\subseteq\bR$ is open.
	Let $\varphi:[0,t_c)\times \bD \to \bR$ be a continuous function satisfying 
	\begin{equation} 
		\sup_{t \in [0, t_c)}|\varphi(t,w_1)-\varphi(t,w_2)| \leq C_4(\varphi)|w_1-w_2|, \; w_1, w_2 \in \bD, \label{eq:eq2127}
	\end{equation}
	for some  constant $C_4(\varphi) \in (0, \infty)$.
	Suppose $w(t)$, $t \in [0,t_c)$, is the unique solution of the differential equation
	$w^\prime(t) = \varphi(t,w(t))$,  $ w(0)=w_0\in\bD$.
	For $n\geq 1$, let $W_n$ be a semimartingale of the form \eqref{eqn:semimart} with $W_n(0)\in\bD$.
	Define $\tau_n:=\inf\{t\geq 0\, :\, W_n(t)\notin\bD\}$.
	Assume that there exist positive deterministic sequences $\{\theta_i(n)\}$, $i=1, 2, 3$, such that whp,
	
	\vskip3pt
	
	\noindent{\upshape (i)}
	$|W_n(0)-w_0|  \leq \theta_1(n)$,
	\hskip15pt
	{\upshape (ii)}		
	$  \int_0^{t_n\wedge\tau_n}\left| \vdd(W_n(t)) - \varphi(t,W_n(t))\right|dt \leq \theta_2(n)$, and
	
	\vskip3pt
	
	\noindent{\upshape (iii)}
	$\langle\vM(W_n), \vM(W_n)\rangle_{t_n\wedge\tau_n} \leq \theta_3(n)$.
	
	\vskip3pt
	
	Then whp,
	$ \sup_{0\leq t \leq t_n\wedge\tau_n}|W_n(t)-w(t)| \leq e^{C_4(\varphi)t_c}(\theta_1(n) + \theta_2(n) + \theta_4(n))$,
	where $\theta_4=\theta_4(n)$ is any sequence satisfying $ \sqrt{\theta_3(n)} = o(\theta_4(n))$.
\end{lem}

We start with the semimartingale decomposition of $Y(\cdot)$. Define
\begin{align}
	F^{\vy}(s_1,y)&:= -{s_1^{-1}}\left[2+4s_1^2 y^2 -8s_1 y \right]\label{eqn:fy-def}\, ,
\end{align}
and let
\begin{equation}
	\label{eqn:epsn-def}
	\eps_n(t)  := 8\max_i f_i(t)/(n\bars_2(t)) \leq 8I(t)/(n\bars_2(t)).
\end{equation}

\begin{lem}\label{lem:semi-mart-y-v}
	The process $Y(\cdot)$ satisfies
	\begin{equation}
		\label{eqn:y-decomp-semi}
		\vdd(Y)(t) = F^{\vy}(\bars_1(t),Y(t)) + O_{t_c}\Big(\frac{I^2(t) Y(t)}{n (1-\eps_n(t))}\Big), \;\; \vv(Y)(t):= O_{t_c}\Big(\frac{I^2(t) Y^2(t)}{n (1-\eps_n(t))^2}\Big).
	\end{equation}
\end{lem}

\noindent {\bf Proof:} 
The jumps of $Y(\cdot)$ are given by
\[
\Delta Y(t)= \chh{\frac{1}{\bars_2(t)+\Delta \bars_2(t)} - \frac{1}{\bars_2(t)} = - \frac{\Delta \bars_2(t)}{\bars_2^2(t)} + \frac{(\Delta \bars_2(t))^2}{\bars_2^2(t)(\bars_2(t)+\Delta \bars_2(t))}}\, . 
\]
If the change happens due to \chh{the} merger of \chh{a half-edge} in $\cC_i\css{(t)}$ with one in $\cC_j\css{(t)}$, then it is easy to check that $\Delta \bars_2(t) \geq 0$. If the change occurs owing to the merger in the same component $\cC_{i^*}\css{(t)}$ then
$|\Delta \bars_2(t)| = \left|{[(f_{i^*}(t)-2)^2 - f_{i^*}^2(t)]}/{n} \right|
\leq 
8\max_i f_i(t)/n.$
Thus,
$\chh{\bars_2(t)} + \Delta \bars_2(t) \geq \bars_2(t)(1-\eps_n(t))$, resulting in
\chh{\[\Delta Y(t) = - \frac{\Delta \bars_2(t)}{\bars_2^2(t)}  + O_{t_c}\Big(\frac{(\Delta \bars_2(t))^2}{\bars_2^3(t) (1-\eps_n(t))}\Big). \]}
Using Lemma~\ref{lem:ds-dg-dd} (a) and (e) yield the desired form of the infinitesimal mean $\vdd(Y)$. For the variance note that $(\Delta Y(t))^2\leq (\Delta \bars_2(t))^2/(\bars_2^4(t)(1-\eps_n(t))^2)$. 
Thus,
\[
\vv(Y)(t) 
\leq 
{\vv(\bars_2)(t)}/(\bars_2^4(t)(1-\eps_n(t))^2) 
= 
O_{t_c}\left({I^2(t) Y^2(t)}/\big(n (1-\eps_n(t))^2\big)\right),
\]
which gives the desired bound.
\qed

\vskip3pt

\noindent{\bf Proof of Proposition \ref{prop:suscep-close-det}:} We first show that
\begin{equation}
	\label{eqn:yt-close-15}
	\sup_{0\leq t\leq t_n} |Y(t) - y(t)| = O_P(n^{-1/5}).
\end{equation}
Recall that $\delta\in(1/6, 1/5)$.
Fix $\eta\in(1-3\delta, 1/2)$.
Then Theorem \ref{thm:config-largest-comp-diam}, Lemma \ref{lem:s1-conc}, and the relation $\bars_2(t)\geq \bars_1(t)$ together imply that there exists a constant $C_* > 0$ such that for all $t \leq t_n$, whp,
\begin{equation}\label{eqn:I-suscep-good}
	I(t) \leq \frac{C_* \css{\log^4{n}}}{(t_c-t)^2}\, , 
	\ \ 
	|\bars_1(t) - \mu\exp(-2t)| \leq n^{-\eta}\, ,
	\ \ \text{ and }\ \ 
	Y(t) \leq \frac{2\exp(2t_c)}{\mu}.
\end{equation}
Thus, 
\begin{equation}\label{eqn:fy-bars-s-bd}
	\sup_{t\leq t_n}\big|F^{\vy}(\bars_1(t), Y(t)) - F^{\vy}(s_1(t), Y(t)) \big| =O_P(n^{-\eta}),
	\ \ \text{ and }\ \ 
	\sup_{t\leq t_n} \eps_n(t) \to 0\, ,
\end{equation}
where $\eps_n$ is as in \eqref{eqn:epsn-def}.
Taking $\bD=(0, 2\exp(2t_c)/\mu)$, it is easy to check that the function 
$\varphi(t, y) := F^{\vy}(s_1(t), y)$ satisfies the Lipschitz condition in \eqref{eq:eq2127}. 
Further, $y'(t)=\varphi(t, y(t))$.
By Condition \ref{ass:cm-degree-det},
$|Y(0) - y(0)| =O(\log^q{n}/\sqrt{n})$.
Letting $\tau_n:=\inf\{t\geq 0\, :\, Y(t)\notin\bD\}$, and writing 
$\cE_n = \int_0^{t_n\wedge\tau_n}|\vdd(Y)(t) - F^{\vy}(s_1(t), Y(t))|dt$ and using Lemma \ref{lem:semi-mart-y-v}, \eqref{eqn:I-suscep-good}, and \eqref{eqn:fy-bars-s-bd} yields
\begin{equation}
	\label{eqn:fy-bound}
	\cE_n
	= 
	O_P\left(\frac{1}{n^{\eta}} + \int_0^{t_n} \frac{\log^8{n}}{n(t_c-t)^4}\chh{dt}\right) 
	= O\Big(\frac{\log^8 n}{n^{1-3\delta}}\Big).
\end{equation}
Using \eqref{eqn:y-decomp-semi} gives
\begin{equation}
	\label{eqn:my-bound}
	\langle \vM(Y), \vM(Y)\rangle_{t_n} 
	= O_P\big(n^{-1}\int_0^{t_n} {I^2(t)} dt\big)
	=O\Big(\frac{\log^8 n}{n^{1-3\delta}}\Big).
\end{equation}
Using Lemma \ref{lem:conv-pure-jump} with
$\theta_1(n) = {\log^{q+1}{n}}/{\sqrt{n}}$ and $\theta_2(n) =\theta_3(n)= \log^9 n/{n^{1-3\delta}}$, and using the fact that $\tau_n>t_n$ whp now completes the proof of \eqref{eqn:yt-close-15}.  

We now strengthen this estimate. 
Note that \eqref{eqn:yt-zt-near-tc} implies that 
$\inf_{0\leq t\leq t_n}y(t)=y(t_n)=\Theta(n^{-\delta})$. 
Since $\delta<1/5$, by \eqref{eqn:yt-close-15}, whp,
\begin{align}\label{eqn:346}
	y(t)/2\leq Y(t)\leq 2y(t)\ \ \text{ for all }\ \ t\in [0, t_n]\, .
\end{align} 
Redoing the above error bounds, we see that \eqref{eqn:fy-bound} can be replaced by
\begin{equation}\label{eqn:fy-new-bound}
	\cE_n 
	= 
	O_P\Big(\frac{1}{n^{\eta}} + \int_0^{t_n} \frac{\log^8 n\cdot y(t)}{n(t_c-t)^4}\chh{dt}\Big) 
	= 
	O\Big(\frac{1}{n^{\eta}}\Big)\, ,
\end{equation}
where we have used the first relation in \eqref{eqn:yt-zt-near-tc}.
Similarly, \eqref{eqn:my-bound} can be improved to 
\[
\langle \vM(Y), \vM(Y)\rangle_{t_n} 
= O_P\big(n^{-1}\int_0^{t_n} I^2(t)y^2(t) dt\big)
=O_P(\log^8 n/n^{1-\delta}).
\] 
Since $\delta\in(1/6, 1/5)$, another application of Lemma \ref{lem:conv-pure-jump} completes the proof of \eqref{eqn:yt-close}.

To prove \eqref{eqn:89}, write $Z(t)=\bars_3(t)Y^3(t)$, $U(t)=\barg(t)Y(t)$, and $V(t)=\barcd(t)Y^2(t)$.
One can write down the semimartingale decompositions of these processes using Lemma~\ref{lem:ds-dg-dd}, Lemma~\ref{lem:semi-mart-y-v}, and some direct calculations.
The claim in  \eqref{eqn:89} would then follow using \eqref{eqn:yt-close}, \eqref{eqn:I-suscep-good}, and Lemma~\ref{lem:conv-pure-jump}.
As mentioned before, we will omit these calculations.
\qed

\vskip5pt

We end this section with a result on $\bars_2^{\star}$ as defined in \eqref{eqn:bars2-star-def}. 
Recall the function $F_2^{\star}$ \chh{from} Lemma \ref{lem:ds-dg-dd} (g). 
Consider the deterministic analogue $s_2^{\star}$ solving the ODE
$(s_2^{\star})^\prime(t) = F_2^{\star}(s_1(t), g(t))$ with $s_2^{\star}(0)=1$.
This ODE can be explicitly solved on $[0,t_c)$ as
\[s_2^{\star}(t):= \Big(1-\frac{\mu}{\nu-1}\Big) + \frac{\mu}{(\nu-1)(\nu-e^{2t}(\nu-1))}\, , \qquad 0\leq t < t_c. \]
Note that at $t_n = t_c- \nu/[2n^{\delta}(\nu-1)]$, we have
$n^{-\delta}{s_2^{\star}(t_n)}\to {\mu}/{\nu^2}$, as  $n\to\infty$.
\begin{lem}\label{lem:s2-star-asymp}
	We have,
	$\sup_{0\leq t\leq t_n}\big|\bars_2^{\star}(t) - s_2^{\star}(t)\big|=o_P(n^{\delta})$.
	Consequently,
	$n^{-\delta}{\bars_2^{\star}(t_n)} \weakc {\mu}/{\nu^2}$. 
\end{lem}

\noindent{\bf Proof:} We only outline the argument.
Note that by \eqref{eqn:yt-close}, \eqref{eqn:346}, and \eqref{eqn:yt-zt-near-tc}, whp, 
\begin{align}\label{eqn:347}
	|\bars_2(t)-s_2(t)|
	=
	\big|\frac{1}{Y(t)}-\frac{1}{y(t)}\big|
	\leq
	\frac{\log^5 n}{  (t_c-t)^2\cdot n^{(1-\delta)/2} }
	\ \ \text{ for all }\ \ t\in[0, t_n]\, .
\end{align}
From Lemma~\ref{lem:ds-dg-dd} \eqref{item:s2-star-variance}, we see that
$\vv(\bars_2^\star)(t)=O_{t_c}\big(I^2(t)U^2(t)/(nY^2(t))\big)$, and consequently,
\begin{align}\label{eqn:348}
	\langle \vM(\bars_2^{\star}), \vM(\bars_2^{\star})\rangle_{t_n}
	=
	O_P\Big(
	\int_0^{t_n} \frac{I^2(t)}{n(t_c-t)^2}dt
	\Big)
	=
	O\Big(\frac{\log^8 n}{n^{1-5\delta}}\Big)\, ,
\end{align}
where the first step uses \eqref{eqn:89}, \eqref{eqn:346}, and  \eqref{eqn:yt-zt-near-tc}, and the second step uses \eqref{eqn:I-suscep-good}.
Finally, recall the expression for $\vdd(\bars_2^\star)$ from Lemma~\ref{lem:ds-dg-dd} \eqref{item:s2-star-mean}, and note that for all $t\in [0, t_n]$,
the quantity
$\big|F_2^{\star}(\bars_1(t), \barg(t))-F_2^{\star}(s_1(t), g(t))\big|$ is bounded above by
\begin{align}\label{eqn:54}
	2\Big|\frac{U^2(t)}{\bars_1(t)} - \frac{u^2(t)}{s_1(t)} \Big|\frac{1}{Y^2(t)}
	+
	2\frac{u^2(t)}{s_1(t)}\big|\bars_2(t)-s_2(t)\big|\Big(\frac{1}{Y(t)}+\frac{1}{y(t)}\Big)\, .
\end{align}
The expression in \eqref{eqn:54} can now be bounded using \eqref{eqn:346}, \eqref{eqn:yt-zt-near-tc}, \eqref{eqn:347}, \eqref{eqn:89}, and Lemma~\ref{lem:s1-conc}.
The desired claim will then follow from an application of Lemma \ref{lem:conv-pure-jump}.
\qed

\subsection{Modified process $\cG_n^{\modi}$}\label{sec:cm-mod-process-def}
From now on, we will work with a fixed $\delta\in (1/6, 1/5)$.
We will now describe a modification of the process $\CM_n(\cdot)$ that evolves like the multiplicative coalescent. 
Note that to obtain $\CM_n(t_c+\lambda/n^{1/3})$ starting from $\CM_n(t_n)$, we run the dynamic construction for an additional
$	r_n(\lambda):= {\nu}/{(2(\nu-1)n^{\delta})} + {\lambda}/{n^{1/3}}$
units of time. 
For $t\in (t_n, t_n+r_n(\lambda)]$, two components $\cC_i(t)$ and $\cC_j(t)$ form an edge 
at rate
\begin{equation}
	\label{eqn:fi-fj-edge-cm}
	f_i(t)\frac{f_j(t)}{n\bars_1(t)-1} + f_j(t)\frac{f_i(t)}{n\bars_1(t)-1} \approx \frac{2\nu f_i(t) f_j(t)}{n\mu(\nu-1)}\css{.}
\end{equation}
Thus, (a) {components} merge at rate (essentially) proportional to the product of {the number} of \emph{alive half-edges} in the respective components, and 
(b) when two components merge, the weight of the new component \css{(measuring the propensity of this component to merge)}  is $f_i(t)+f_j(t)-2$ (rather than {$f_i(t)+f_j(t)$}) since two half-edges were used \css{up in} the merger. {These observations motivate the following modification.}

\vskip3pt

\noindent{\bf Construction of $\cG_n^{\modi}$:} Recall that $s_1(t_c) = \mu(\nu-1)/\nu$. 
Let $\FR(t_n)$ denote the set of alive half-edges at time $t_n$.  
Recall that the connected components in $\CM_n(t_n)$ are called ``blobs.'' 
We will write $\fb_i$ for $\cC_i(t_n)$, and $f_i^\circ$ for $f_i(t_n)$. 
Consider independent Poisson processes $\cP_{\ve}$ indexed by ordered pairs $\ve = (u,v) \in \FR(t_n) \times \FR(t_n)$, each having rate 
\begin{equation}\label{eqn:rate-modi-def}
	\nu/(n\mu(\nu-1))=1/(n s_1(t_c)).
\end{equation} 
Every time $\cP_{\ve}$ rings, complete the edge $(u,v)$ but continue to consider $u,v$ as alive.  
Run this process for time $r_n(\lambda)$, and let $\cG_n^{\modi}(t_c+\lambda/n^{1/3})$ be the process at time $r_n(\lambda)$. 
The rate of creation of edges between two blobs $\fb_i$ and $\fb_j$ in the modified process is given by
$
2\nu f_i^\circ f_j^\circ/(n\mu(\nu-1))
$.
Compare this with \eqref{eqn:fi-fj-edge-cm}.  
Note that the construction of $\cG_n^{\modi}(t_c+\lambda/n^{1/3})$ differs from that of the space $\Gamma(\cG, \vx, \vM, \vX)$ described in Section \ref{sec:inter-blob-distance} in two ways:
In $\cG_n^{\modi}(t_c+\lambda/n^{1/3})$,
(i) 
some new edges with both endpoints within the same blob may be created, and
(ii) 
there may be blobs that are connected by two or more edges.
None of these affects the maximal components as we will see in Lemma \ref{lem:no-multi-intra} below.

Conditional on $\CM_n(t_n)$, the connectivity structure between different blobs in $\cG_n^{\modi}(t_c+\lambda/n^{1/3})$ has the following two step description:
(Step 1)
For $i\neq j$, connect $\fb_i$ and $\fb_j$ (independently across pairs of blobs) by a single edge with probability
\begin{equation}
	\label{eqn:cm-modi-initial-p}
	p_{ij} =1- \exp\Big(-f_i^\circ f_j^\circ\Big(\frac{1}{n^{1+\delta}}\frac{\nu^2 }{\mu(\nu-1)^2} + \frac{1}{n^{4/3}}\frac{2\nu}{\mu(\nu-1)}\lambda\Big) \Big)\, ,
\end{equation}
where each junction point in a blob $\fb$ is chosen according to the probability measure that assigns to each vertex $v\in\fb$ a mass proportional to the number of alive (at time $t_n$) half-edges attached to $v$.
By \eqref{eqn:cm-modi-initial-p}, conditional on $\CM_n(t_n)$, the blob-level superstructure constructed thus far has the same distribution as the random graph $\cG(\vx,q)$ (see Section \ref{sec:smc-def}), where
\begin{equation}
	\label{eqn:xi-q-modi}
	x_i = \frac{\beta^{1/3}}{\mu(\nu-1)} \frac{f_i^\circ}{n^{2/3}}\ ,\ \ \text{ and }\ \ 
	q= q_{\lambda}
	= n^{1/3-\delta}\frac{\mu\nu^2}{\beta^{2/3}} + \frac{2\mu(\nu-1)\nu}{\beta^{2/3}}\lambda\, .
\end{equation}
(Step 2)
Let $x_i$ and $q$ be as in \eqref{eqn:xi-q-modi}.
If an edge is created between $\fb_i$ and $\fb_j$ in Step 1, then add an additional $Y_{i, j}-1$ many edges between them, where $Y_{i, j}$ is distributed as a Poisson random variable with mean $qx_i x_j$ conditioned to be at least one; do this independently across different pairs of blobs that were connected by an edge in Step 1.

\subsubsection{Components of $\cG_n^{\modi}$ ranked by the number of alive half-edges at time $t_n$}
\label{sec:comp-free-cm-proof}  
Recall that $\fb_j$ denotes the $j$-th largest blob, and $f_j^\circ=f_j(t_n)$ denotes the number of alive (at time $t_n$) half-edges attached to the vertices in $\fb_j$.
Define the {\bf weight} of a component $\cC$ in $\cG_n^{\modi}(t_c+\lambda/n^{1/3})$ by $\sW(\cC) =  \sum_{i:\fb_i \subseteq \cC} f_i^\circ$.  
Write $\sC_{i}^{\modi}({t_c+\lambda/n^{1/3}})$ for the component in $\cG_n^{\modi}({t_c+\lambda/n^{1/3}})$ with the $i$-th largest weight.  
We will now use Theorem \ref{thm:aldous-review} to understand the maximal \css{weighted} components in $\cG_n^{\modi}({t_c+\lambda/n^{1/3}})$.
To use Theorem \ref{thm:aldous-review}, we use Theorems \ref{thm:config-bare-sibcrit} and \ref{thm:config-largest-comp-diam} to verify Condition \ref{ass:aldous-basic-assumption}. 
First, note that \eqref{eqn:s3-asymp-cm} implies
\begin{equation}
	\label{eqn:cm-first-aldous-ver}
	\frac{\sum_i x_i^3}{(\sum_i x_i^2)^3} = \frac{[\mu(\nu-1)]^3}{\beta} \css{\frac{\bars_3(t_n)}{[\bars_2(t_n)]^3}} \weakc 1\, .
\end{equation}
Next, using \eqref{eqn:s2-asymp-cm} and \eqref{eqn:xi-q-modi}, for any $\lambda\in\bR$ and a sequence $\lambda_n\to\lambda$,
\begin{align}
	\label{eqn:cm-sec-aldous-ver}
	q_{\lambda_n}-\frac{1}{\sum_i x_i^2} 
	= 
	q_{\lambda_n}- \frac{(\mu(\nu-1))^2}{\beta^{2/3}}\frac{n^{1/3}}{\bars_2(t_n)} \weakc \frac{2\mu\nu(\nu-1)}{\beta^{2/3}}\lambda
	=:\alpha(\lambda)\, .
\end{align}
Theorem \ref{thm:config-largest-comp-diam}, \eqref{eqn:s2-asymp-cm}, and $\delta \in (1/6,1/5)$ {implies}
\begin{align}\label{eqn:333}
	{x_{\max}}/{\sum_i x_i^2} = O_P\left({\log^4{n}}/{n^{1/3-\delta}}\right) \weakc 0\, .
\end{align}

\begin{lem}\label{lem:no-multi-intra}
	Fix $\lambda\in\bR$ and a sequence $\lambda_n\to\lambda$.
	Fix $i\geq 1$, and write $\fC_i^{\modi}$ for $\fC_i^{\modi}(t_c+\lambda_n/n^{1/3})$.
	Let $\fN_{\multi}$ be the number of pairs of (distinct) blobs $\fb_j, \fb_k\subseteq \fC_i^{\modi}$ that are connected by at least two edges.
	Let $\fN_{\intra}$ denote the number of edges created in the time interval $[t_n, t_c+\lambda_n n^{-1/3}]$ with both endpoints inside a blob that is contained in $\fC_i^{\modi}$.
	Then whp, $\fN_{\multi}=0$ and $\fN_{\intra}=0$.
\end{lem}

\noindent{\bf Proof:}
Let $q, \vx$ be as in \eqref{eqn:xi-q-modi}.
Let $\sC^{\modi; \sss\otimes}_i$ denote the graph obtained from $\sC^{\modi}_i$ by viewing each blob $\fb_j\subseteq \sC^{\modi}_i$ as a single vertex labeled $j$, and then ignoring the multiple edges between different blobs.
Set $\sW\big(\sC^{\modi; \sss\otimes}_i\big):=\sW\big(\sC^{\modi}_i\big)$.
Let $\fN_i^{\modi; \sss\otimes}$ be the number of surplus edges in $\sC^{\modi; \sss\otimes}_i$.
By Theorem \ref{thm:aldous-review}, \eqref{eqn:cm-first-aldous-ver}, \eqref{eqn:cm-sec-aldous-ver}, \eqref{eqn:333}, and the discussion around \eqref{eqn:cm-modi-initial-p}, 
\begin{equation}\label{eqn:cong-weight-free-prelim-version}
	\bigg(
	\Big(
	\frac{\beta^{1/3}}{\mu (\nu-1)n^{2/3}}\sW\big(\sC^{\modi; \sss\otimes}_i\big)\, ,\,
	\fN_i^{\modi; \sss\otimes}
	\Big)
	\, ;\,  i\geq 1\bigg) 
	\weakc 
	\mvXi\big(\alpha(\lambda)\big)
\end{equation}
with respect to the product topology on  $(\bR^2)^{\bN}$, where $\alpha(\lambda)$ is as in \eqref{eqn:cm-sec-aldous-ver}.
Write 
$\bE_1[\cdot]:=\bE[\, \cdot\, |\, \CM_n(t_n), \fC_i^{\modi; \sss\otimes}]$.
By the discussion in Step 2 below \eqref{eqn:xi-q-modi}, on the event $q x_{\max}^2\leq 1$,
\[
\bE_1\big[\fN_{\multi}\big]
\leq
\sum_{\{j, k\}\in E(\fC_i^{\modi; \sss\otimes})} e q_{\lambda_n} x_j x_k
\leq
eq_{\lambda_n}
\big(x_{\max}\sum_{j\in\fC_i^{\modi; \sss\otimes}} x_j
+
x_{\max}^2\cdot\fN_i^{\modi; \sss\otimes} \big)\, ,
\]
where the last step can be verified easily.
By \eqref{eqn:cong-weight-free-prelim-version}, $\fN_i^{\modi; \sss\otimes}=O_P(1)$ and $\sum_{j\in\fC_i^{\modi; \sss\otimes}} x_j=O_P(1)$.
It thus follows from \eqref{eqn:333}, \eqref{eqn:cm-sec-aldous-ver}, and the fact $x_{\max}\weakc 0$ that $\fN_{\multi}\weakc 0$.

Now, 
$\bE_1[\fN_{\intra}]
\leq 
q_{\lambda_n}\sum_{j\in\fC_i^{\modi; \sss\otimes}}x_j^2
\leq
q_{\lambda_n}x_{\max}\cdot\sum_{j\in\fC_i^{\modi; \sss\otimes}}x_j
=
q_{\lambda_n}x_{\max}\cdot O_P(1)
=
o_P(1)
$,
where the last step uses \eqref{eqn:333} and \eqref{eqn:cm-sec-aldous-ver}.
Consequently, 
$\fN_{\intra}\weakc 0$.
\qed

\vskip3pt

Lemma \ref{lem:no-multi-intra} shows that the number of surplus edges in $\sC^{\modi}_i(t_c+\lambda_n/n^{1/3})$ that were created in the time interval $[t_n, t_c+\lambda_n n^{-1/3}]$ equals $\fN_i^{\modi; \sss\otimes}$ whp.
Thus, \eqref{eqn:cong-weight-free-prelim-version} implies the following result.

\begin{prop}\label{prop:modi-free-max}
	Consider $\lambda\in\bR$ and a sequence $\lambda_n\to\lambda$.
	Let $\fN_i^{\modi}(t_c+\lambda_n/n^{1/3})$ be the number of surplus edges in $\sC^{\modi}_i(t_c+\lambda_n/n^{1/3})$ that were created in the time interval $[t_n, t_c+\lambda_n n^{-1/3}]$. 
	Then
	\begin{equation}\label{eqn:cong-weight-free}
		\bigg(
		\Big(
		\frac{\beta^{1/3}}{\mu (\nu-1)n^{2/3}}\sW\Big(\sC^{\modi}_i\big(t_c+\frac{\lambda_n}{n^{1/3}}\big)\Big),
		\fN_i^{\modi}\big(t_c+\frac{\lambda_n}{n^{1/3}}\big)
		\Big)
		\, ;\,  i\geq 1\bigg) 
		\weakc 
		\mvXi\big(\alpha(\lambda)\big)
	\end{equation}
	with respect to the product topology on  $(\bR^2)^{\bN}$, where $\alpha(\lambda)$ is as in \eqref{eqn:cm-sec-aldous-ver}.
\end{prop}
Recall from Section \ref{sec:gr-constr} that for a finite graph $\cG$, $|\cG|$ denotes the number of vertices in $\cG$.
\begin{prop}[Number of vertices]\label{prop:modi-max-comp-sizes}
	For any $\lambda\in\bR$ and any sequence $\lambda_n\to\lambda$,
	\begin{equation}\label{eqn:cong-comp-free}
		\bigg(
		\Big(
		\frac{\beta^{1/3}}{\mu n^{2/3}}\big|\sC^{\modi}_i\big(t_c+\frac{\lambda_n}{n^{1/3}}\big)\big|\, ,\,
		\fN_i^{\modi}\big(t_c+\frac{\lambda_n}{n^{1/3}}\big)
		\Big)
		\, ;\,  i\geq 1\bigg) 
		\weakc 
		\mvXi\big(\alpha(\lambda)\big)
	\end{equation}
	with respect to the product topology on  $(\bR^2)^{\bN}$, where $\alpha(\lambda)$ is as in \eqref{eqn:cm-sec-aldous-ver}.
\end{prop}

\noindent{\bf Proof:} 
The proof is similar to that of Proposition \ref{prop:moments-convergence-each-component}; we describe it briefly. 
Recall that in the proof of Proposition \ref{prop:moments-convergence-each-component}, a key role was played by a breadth-first exploration of the random graph $\cG(\vx,q)$ where the vertices appeared in a size-biased order with the sizes given by $\vx$, and where we further kept track of an associated weight sequence $\vu$ whose components were the average distances.

In the present setting, by the discussion around \eqref{eqn:xi-q-modi}, ignoring multiple edges between blobs and intra-blob edges, the blob-level superstructure of \css{$\cG_n^{\modi}(t_c+\lambda_n n^{-1/3})$} has the same distribution as the random graph $\cG(\vx,q_{\lambda_n})$, where $\vx$ and $q$ are as in \eqref{eqn:xi-q-modi}. 
Set
$\vu:=\big(|\fb_j|;\, j\geq 1\big)$. 
Writing $\sC^{\modi}_i$ for $\sC^{\modi}_i\big(t_c+\lambda_n/n^{1/3}\big)$, for any $i\geq 1$,
\begin{align}\label{eqn:1212}
	\frac{\sum_{j:\fb_j\subseteq\sC^{\modi}_i} f_j^\circ}{\sum_{j:\fb_j\subseteq\sC^{\modi}_i} |\fb_j|}
	\cdot
	\frac{\sum_j f_j^\circ |\fb_j|}{\sum_j (f_j^\circ)^2}
	=
	\frac{\sum_{j:\fb_j\subseteq\sC^{\modi}_i} x_j}{\sum_{j:\fb_j\subseteq\sC^{\modi}_i} u_j}
	\cdot
	\frac{\sum_j x_ju_j}{\sum_j x_j^2}
	\weakc 
	1\, ,
\end{align}
where the second step follows from an argument identical to the one used in the proof of \eqref{eqn:prop67-equiv}. 
(We omit the verification of the conditions in \eqref{eqn:ass-size-biased-partial-sum}, as they follow easily from Theorem \ref{thm:config-bare-sibcrit}, Theorem \ref{thm:config-largest-comp-diam}, and Lemma \ref{lem:s1-conc}.)
Now,
$
\sum_j (f_j^\circ)^2/\sum_j f_j^\circ |\fb_j|  
= 
\chhh{ \bar{s}_2(t_n) }/\chhh{\bar{g}(t_n)} \weakc (\nu - 1)
$,
where the last assertion follows from \eqref{eqn:s2-asymp-cm} and \eqref{eqn:d1-g-asymp-cm}.
This observation combined with \eqref{eqn:1212} shows that for any $i\geq 1$,
\begin{equation}
	\label{eqn:1214}
	\frac{\sW(\sC^{\modi}_i({t_c+\lambda_n/n^{1/3}}))}{|\sC^{\modi}_i({t_c+\lambda_n/n^{1/3}})|}
	\weakc 
	\nu-1\, ,\ \ \text{ as }\ \ \ n\to\infty\, .
\end{equation}
Combining \eqref{eqn:1214} with Proposition \ref{prop:modi-free-max} completes the proof. \qed

For $\lambda\in\bR$ and $i\geq 1$, endow $\sC_i^{\modi}(t_c+\lambda n^{-1/3})$ with the measure that assigns, to each vertex $v$, a mass that is equal to the number of half-edges of $v$ that were alive at time $t_n$, and let 
$\sC_{i; \free}^{\modi}(t_c+\lambda n^{-1/3})$ be the resulting metric measure space.
Then the next result follows from an application of Lemma \ref{lem:no-multi-intra} coupled with Theorem \ref{thm:aldous-gen-2} upon using Theorem \ref{thm:config-bare-sibcrit} and Theorem  \ref{thm:config-largest-comp-diam} to verify that Condition \ref{ass:aldous-gen-2} is met
(Condition \ref{ass:aldous-basic-assumption} has already been verified in \eqref{eqn:cm-first-aldous-ver}, \eqref{eqn:cm-sec-aldous-ver}, and \eqref{eqn:333}).

\begin{thm}\label{thm:modi-scaling-lim}
	Let $\zeta_1 = {\beta^{2/3}}/(\mu\nu)$ and $ \zeta_2:= {\beta^{1/3}}/(\mu(\nu-1)) $. 
	Then for any $\lambda\in\bR$ and any sequence $\lambda_n\to\lambda$, as $n\to\infty$,
	\[
	\Big(\scl\Big(\frac{\zeta_1}{n^{1/3}}, \frac{\zeta_2}{n^{2/3}}\Big)
	\sC_{i; \free}^{\modi}\bigg({t_c+\frac{\lambda_n}{n^{1/3}}}\bigg)\, ;\ i\geq 1\Big) \weakc \vCrit\Big(\frac{2\mu\nu(\nu-1)}{\beta^{2/3}}\lambda\Big)\, . \]
\end{thm}

\subsubsection{Scaling limit of components in $\cG_n^{\modi}$ endowed with the counting measure}
\label{sec:34}
We will denote the metric measure space obtained by endowing $\sC_i^{\modi}(t_c+\lambda n^{-1/3})$ with the counting measure on its vertices by the same notation.

\begin{thm}\label{thm:modi-scalin-lim-count}
	For any $\lambda\in\bR$ and any sequence $\lambda_n\to\lambda$, as $n\to\infty$,
	\[
	\Big(\scl\Big(\frac{\beta^{2/3}}{\mu\nu n^{1/3}}, \frac{\beta^{1/3}}{\mu n^{2/3}}\Big)
	\sC_i^{\modi}\left(t_c+\frac{\lambda_n}{n^{1/3}}\right);\, i\geq 1\Big) \weakc \vCrit\Big(\frac{2\mu\nu(\nu-1)}{\beta^{2/3}}\lambda\Big).  
	\]
\end{thm}

\noindent{\bf Proof:}
We only prove convergence of $\sC_1^{\modi}(t_c+\lambda_n/n^{1/3})$.
Throughout the rest of this section we will simply write $\sC_{1}^{\modi}$ to mean $\sC_1^{\modi}(t_c+\lambda_n/n^{1/3})$.
To derive the result in Theorem \ref{thm:modi-scalin-lim-count} for $\sC_{1}^{\modi}$ from Theorem \ref{thm:modi-scaling-lim}, it is enough to show that
\begin{equation}\label{eqn:enough-to-show-count-fr}
	n^{-2/3}\sum_{i:\fb_i \subseteq \sC_1^{\modi}} 
	\Big|
	\big|\fb_i \big| - \frac{f_i^\circ}{(\nu-1)}
	\Big|
	\weakc 0.
\end{equation}
As shown in the proof of Proposition \ref{prop:moments-convergence-each-component} (around \eqref{eqn:1438}), for any $\eps>0$ there exists $T_{\eps}>0$ such that in the breadth-first exploration of the graph $\cG(\vx, q)$, with probability at least $1-\eps$, the maximal component is completely explored before the $(T_{\eps}\sum_i x_i/\sum_i x_i^2)$-th step.
Using \eqref{eqn:s2-asymp-cm} and Lemma \ref{lem:s1-conc}, 
$\sum_i x_i/\sum_i x_i^2=\Theta_P(n^{2/3-\delta})$.
In view of the discussion around \eqref{eqn:cm-modi-initial-p} and the fact that the vertices appear in a size-biased order in the breadth-first exploration of $\cG(\vx, q)$, \eqref{eqn:enough-to-show-count-fr} will follow if we show that for any $T>0$, 
$
n^{-2/3}
\sum_{i=1}^{Tn^{2/3-\delta}}\Big|
\big|\fb_{v(i)}\big| - f_{v(i)}^\circ/ (\nu-1)\Big|
\weakc 0
$,
as $n\to\infty$,
where $(\fb_{v(1)}, \fb_{v(2)},\ldots)$ is a size-biased permutation of the blobs using the weights $(f_1^\circ, f_2^\circ, \ldots)$.
Now the following lemma completes the proof of \eqref{eqn:enough-to-show-count-fr}.
\qed

\begin{lem}\label{lem:comp-proof-modi-count}
	Let $(\fb_{v(1)}, \fb_{v(2)},\ldots)$ be as above. 
	Then for any $T>0$, as $n\to\infty$,
	\[n^{-2/3}
	\sum_{i=1}^{Tn^{2/3-\delta}}\Big|
	\big|\fb_{v(i)}\big| - f_{v(i)}^\circ/ (\nu-1)\Big|
	\weakc 0.\]
\end{lem}

\noindent {\bf Proof:}  
We will apply Lemma \ref{lem:size-biased-partial-sum} with 
$\ell=Tn^{2/3-\delta}$, $x_i = f_i^\circ$ and $u_i = \big| |\fb_i|-f_i^\circ/(\nu-1)\big|$.
Suppose we show that 
\begin{align}\label{eqn:3434}
	\sum_j f_j^\circ \big| |\fb_j|-f_j^\circ/(\nu-1)\big|
	=
	\Omega_P(n)\, .
\end{align} 
Then the conditions in \eqref{eqn:ass-size-biased-partial-sum} can be verified using \eqref{eqn:3434},
Theorem \ref{thm:config-largest-comp-diam}, Lemma \ref{lem:s1-conc}, Condition \ref{ass:cm-degree-det}~(a), and the fact $\delta<1/5$.
Thus, an application of Lemma \ref{lem:size-biased-partial-sum}  implies
\begin{align*}
	\sum_{i=1}^{Tn^{2/3-\delta}}
	\Big| \big|\fb_{v(i)}\big| - f_{v(i)}^\circ/ (\nu-1) \Big|
	\sim 
	Tn^{2/3-\delta} \frac{\sum_j f_j^\circ \big| |\fb_j|-f_j^\circ/(\nu-1)\big|}{\sum_j f_j^\circ}\, .
\end{align*}
Using the relation $\sum_j f_j^\circ=\Theta_P(n)$ and the Cauchy-Schwarz inequality yields
\begin{align*}
	\sum_{i=1}^{Tn^{2/3-\delta}}
	\Big| \big|\fb_{v(i)}\big| - \frac{f_{v(i)}^\circ}{\nu-1} \Big|
	= 
	\Theta_P\Big(n^{-\frac{1}{3}-\delta}\cdot
	\Big[ \Big(\sum_j (f_j^\circ)^2\Big)  \Big(\sum_j \Big(|\fb_j| - \frac{f_j^\circ}{\nu-1}\Big)^2\Big)\Big]^{1/2}
	\Big).
\end{align*}
The bound on the right can be written as
\[
n^{-\frac{1}{3}-\delta}\cdot
\sqrt{n \bars_2(t_n)} \Big[n\Big(
\bars_2^{\star}(t_n) + (\nu-1)^{-2}{\bars_2(t_n)} - 2 (\nu-1)^{-1}{\barg(t_n)}
\Big)\Big]^{1/2}.
\]
Using Theorem \ref{thm:config-bare-sibcrit} for $\bars_2$ and $\barg$, 
and Lemma \ref{lem:s2-star-asymp} for $\bars_2^{\star}$ yields the claim in Lemma \ref{lem:comp-proof-modi-count}.

Now we have to prove \eqref{eqn:3434}.
Pick $r\in\bZ_{>0}$ with $p_r>0$ and $r\neq \nu-1$.
Let $A_n:=\big\{i\in[n]\, :\, d_i=r \big\}$.
Let $A_n'$ be the set of those vertices in $A_n$ all of whose half-edges are alive at time $t_n$.
Then
\begin{align}\label{eqn:3434-a}
	\sum_j f_j^\circ \big| |\fb_j|-f_j^\circ/(\nu-1)\big|
	\geq
	\sum_{j\in A_n'} r \big| 1 - r/(\nu-1) \big|
	=
	\big| 1 - r/(\nu-1) \big|\cdot|A_n'|\, .
\end{align} 
Lemma \ref{lem:cmn-t-k-2k-equiv} implies that
$
\bE\big[ |A_n'|\, \big|\, \bars_1(t_n)\big]
=|A_n|\cdot (n\bars_1(t_n))_r/(n\mu^{\sss (n)})_r 
$,
where $(m)_r=m(m-1)\ldots (m-r+1)$.
Further,
$
\var\big[ |A_n'|\, \big|\, \bars_1(t_n)\big]
\leq
\sum_{i\in A_n}\var\big[ \ind_{\{i\in A_n'\}}\, \big|\, \bars_1(t_n)\big]
\leq n
$,
where the first step uses Theorem \ref{thm:negative-association}.
These observations together with the fact that $|A_n|\geq np_r/2$ for all large $n$ under Condition \ref{ass:cm-degree-det}, and Lemma \ref{lem:s1-conc} shows that $|A_n'|=\Omega_P(n)$, which combined with \eqref{eqn:3434-a} completes the proof of \eqref{eqn:3434}.
\qed

\subsection{Proof of Theorem \ref{thm:crit-main-res-cm}}\label{sec:proof-dynamic-cm-scaling-limit}
The proof proceeds in several steps.

\subsubsection{Coupling $\Perc_n(\cdot)$ and $\CM_n(\cdot)$}
Consider $\lambda\in\bR$ and a sequence $\lambda_n\to\lambda$.
Then using Condition \ref{ass:cm-degree-det}~(b), 
for any $\eta\in(1/3, 1/2)$, whp,
\begin{equation}\label{eqn:half-edge-perc-num}
	\Big|
	\#E(\Perc_n(p(\lambda_n)))- \frac{1}{2}\big({n\mu}/\nu + n^{2/3}\mu \lambda_n\big)
	\Big| 
	\leq 
	n^{1-\eta}.
\end{equation}

Simillarly, using Lemma \ref{lem:s1-conc}, for any $\lambda'\in\bR$, a sequence $\lambda'_n\to\lambda'$, and $\eta'\in(1/3, 1/2)$, 
\begin{equation}\label{eqn:half-edge-cm-num}
	\Big|
	\# E\big(\CM_n(t_c+\lambda'_n/n^{1/3})\big) -
	\Big(
	\frac{n\mu}{2\nu}+ n^{2/3}\mu\lambda'_n \frac{(\nu-1)}{\nu}
	\Big)
	\Big| 
	\leq 
	n^{1-\eta'}
\end{equation}
whp.
Comparing \eqref{eqn:half-edge-perc-num} and \eqref{eqn:half-edge-cm-num} and using Lemma \ref{lem:cmn-t-k-2k-equiv} and Proposition \ref{prop:perc-cm},  
we see that for any $\eta\in (1/3, 1/2)$, $\lambda\in\bR$, and a sequence $\lambda_n\to\lambda$, there exist couplings for all $n$ in which
\begin{equation}\label{eqn:couple-cm-perc}
	\CM_n\Big(t_c+ \frac{\nu}{2(\nu-1)}\frac{\lambda_n}{n^{1/3}} -\frac{1}{n^{\eta}}\Big)
	\subseteq 
	\Perc_n(p(\lambda_n))
	\subseteq 
	\CM_n\Big(t_c+ \frac{\nu}{2(\nu-1)}\frac{\lambda_n}{n^{1/3}}+\frac{1}{n^{\eta}}\Big)
\end{equation}
whp; similarly, there exist couplings for all $n$ such that whp,
\begin{equation}\label{eqn:couple-cm-perc-a}
	\Perc_n\Big(\frac{1}{\nu}+\frac{2(\nu-1)}{\nu}\frac{\lambda_n}{n^{1/3}}-\frac{1}{n^{\eta}}\Big)
	\subseteq 
	\CM_n\Big(t_c + \frac{\lambda_n}{n^{1/3}}\Big)
	\subseteq
	\Perc_n\Big(\frac{1}{\nu}+\frac{2(\nu-1)}{\nu}\frac{\lambda_n}{n^{1/3}}+\frac{1}{n^{\eta}}\Big).
\end{equation}

\subsubsection{Properties of maximal components of $\CM_n$}
We first state a result on the component sizes and surplus edges in $\Perc_n(\cdot)$.

\begin{thm}[\cite{dhara-hofstad-leeuwaarden-sen-cm-finite-third-moment}, Theorem 3.6]\label{thm:cm-perc-component-size-scaling}
	Consider $\lambda\in\bR$ and a sequence $\lambda_n\to\lambda$.
	Then under Condition \ref{ass:cm-degree-det},
	\[
	\bigg(
	\Big(
	\frac{\beta^{1/3}}{\mu n^{2/3}}|\cC_{\sss i,\Perc}(\lambda_n)|\, ,\, 
	\spls\big(\cC_{\sss i,\Perc}(\lambda_n)\big) 
	\Big)\, ; \,
	i\geq 1
	\bigg)
	\weakc
	\mvXi\Big(\frac{\mu\nu^2\lambda}{\beta^{2/3}}\Big)
	\]
	with respect to the product topology on $(\bR^2)^{\bN}$.
\end{thm}

\begin{rem}
	Let us clarify two points here.
	Firstly, in \cite[Theorem 3.6]{dhara-hofstad-leeuwaarden-sen-cm-finite-third-moment}, the description of the limiting sequence is slightly different.
	One can use Brownian scaling to rephrase this result in terms of the notation $\mvXi(\cdot)$; see \cite[Theorem 3.5]{addarioberry-sen}.
	Secondly, the edge retention probability in 
	\cite[Display (3.8)]{dhara-hofstad-leeuwaarden-sen-cm-finite-third-moment} is different from $p(\lambda_n)$.
	One can simply repeat the arguments in \cite{dhara-hofstad-leeuwaarden-sen-cm-finite-third-moment} to get the result for a sequence $\lambda_n\to\lambda$ (instead of considering a fixed $\lambda$).
	That $\nu_n$ in \cite[Display (3.8)]{dhara-hofstad-leeuwaarden-sen-cm-finite-third-moment} can be replaced by $\nu$ as in this paper follows from Condition \ref{ass:cm-degree-det}~(b).
\end{rem}

Note that Theorem \ref{thm:cm-perc-component-size-scaling} applied to both the random graph on the extreme left and the one on the extreme right in \eqref{eqn:couple-cm-perc-a} yields the same scaling limit.
This leads to the following result.

\begin{thm}	\label{thm:joseph}
	Let $N_i(t):=\spls(\cC_i(t))$, $i\geq 1$.
	Fix $\lambda\in\bR$ and a sequence $\lambda_n\to\lambda$.

	\begin{enumeratea}
		\item 
		The convergence in \eqref{eqn:cong-comp-free} continues to hold if we replace
		$\fC_i^{\modi}\big(t_c+\lambda_n n^{-1/3}\big)$ by $\cC_i\big(t_c+\lambda_n n^{-1/3}\big)$ and 
		$\fN_i^{\modi}\big(t_c+\lambda_n n^{-1/3}\big)$ by $N_i\big(t_c+\lambda_n n^{-1/3}\big)$.
		\item
		For any $K\geq 1$, in the coupling of \eqref{eqn:couple-cm-perc-a},
		\[
		\pr\Big(
		\cC_{\sss i, \Perc}\big(\fq_n^-\big)
		\subseteq
		\cC_i\big(t_c+\lambda_n n^{-1/3}\big)
		\subseteq
		\cC_{\sss i, \Perc}\big(\fq_n^+\big) \ \ \text{ for }\ \ 1\leq i\leq K
		\Big)\to 1\, ,
		\]
		where $\fq_n^\pm:=2(\nu-1)\lambda_n/\nu \pm n^{-\eta+1/3}$.
	\end{enumeratea}
\end{thm}

Any component $\cC$ of $\CM_n(t_c+\lambda n^{-1/3})$ is made up \chhh{of} a collection of blobs.
Analogous to the notation introduced in Section \ref{sec:comp-free-cm-proof}, let
$
\sW(\cC):= \sum_{i:\fb_i\subseteq\cC} f_i^\circ
$.

\begin{thm}\label{thm:cm-blob-level-surp-blob}
	Fix $\lambda\in\bR$ and a sequence $\lambda_n\to\lambda$.
	For simplicity, write $\cC_i$ for $\cC_i(t_c+\lambda_n /n^{1/3})$.
	Then the following hold: 
	
	\noindent {\upshape (a)}
	For any $i\geq 1$, 
	whp all surplus edges in $\cC_i$ are created in $[t_n, t_c+\lambda_n/n^{1/3}]$ and have endpoints in different blobs, and further, no two blobs in $\cC_i$ are connected by more than one edge.
	
	\noindent{\upshape (b)}
	$\sW(\cC_i)$ has the same asymptotics as $\sW\big(\fC_i^{\modi}\big(t_c+\lambda_n/n^{1/3}\big)\big)$ as given in Proposition \ref{prop:modi-free-max}, i.e.,
	\begin{align}\label{eqn:3456}
		\Big(
		\frac{\beta^{1/3}}{\mu (\nu-1)n^{2/3}} \sW(\cC_i)\, ;\,  i\geq 1
		\Big) 
		\weakc 
		\mvxi\Big(\frac{2\mu\nu(\nu-1)}{\beta^{2/3}}\lambda \Big) 
	\end{align}
	with respect to the product topology.
\end{thm}

For the rest of this section, $\lambda\in\bR$ and the sequence $\lambda_n\to\lambda$ will be fixed, and we will continue to write $\cC_i$ for $\cC_i(t_c+\lambda_n /n^{1/3})$.
We will use the following elementary result in the proof of Theorem \ref{thm:cm-blob-level-surp-blob} whose proof we omit.

\begin{lem}\label{lem:girth-auxiliary}
	{\upshape (a)}
	There exists $h:\bN\cup\{0\}\to\bN$ such that the following holds: Consider any finite connected graph $\cG$ with $\spls(\cG) \leq \fs$.
	Let $\Gamma_0=\big(u_1, u_2, \ldots, u_{\ell}\big)$ be a self-avoinding path in $\cG$ between the vertices $u_1$ and $u_{\ell}$.
	Then for any $k\geq 1$, $\#\{ v\in\Gamma_0\, :\, d_{\cG}(v, u_1)\leq k \}\leq k h(\fs)$.
	(For example, $h(\fs):=4(\fs+1)^2$ works.)
	
	\noindent{\upshape (b)} 
	Let $\cG$ be a finite connected graph.
	Let $\tilde\cG$ be the graph obtained from $\cG$ by adding an edge between two vertices $u$ and $v$.
	Then $\girth(\tilde\cG)=\min\{\girth(\cG),\, d_{\cG}(u, v) +1\}$.
\end{lem}

\noindent{\bf Proof of Theorem \ref{thm:cm-blob-level-surp-blob}~(a):}
Using Theorem \ref{thm:config-largest-comp-diam}, it follows that there exists $C_1>0$ such that $\max\{\girth(\fb_i)\, :\, \spls(\fb_i)\geq 1\} \leq C_1 n^{\delta}(\log n)^3$ whp.
Theorem \ref{thm:joseph}~(b) implies that in the coupling of \eqref{eqn:couple-cm-perc-a}, 
for any $i\geq 1$, 
$\girth(\cC_i)\geq \girth(\cC_{\sss i, \Perc} (\fq_n^+))$ whp.
Thus, the next result completes the proof.
\qed

\begin{lem}\label{lem:girth}
	For any $\eps>0$ and $i\geq 1$, $\girth(\cC_{\sss i, \Perc} (\fq_n^+))\geq n^{\frac{1}{3}-\eps}$ whp.
\end{lem}

\noindent{\bf Proof:}
Consider the exploration of a configuration model given in 
\cite[Algorithm1, page 9]{dhara-hofstad-leeuwaarden-sen-cm-finite-third-moment}.
We do not describe this exploration here; rather we collect some properties of this exploration that will be useful to us.
In this algorithm, each component is explored in a depth-first manner and a depth-first spanning tree of the component is built iteratively.
We refer to the first vertex to be explored in a component $\cC$ as the root of $\cC$, and for any $v\in\cC$, the path in the depth-first spanning tree of $\cC$ connecting $v$ to the root of $\cC$ will be called the ancestral line of $v$.
Under \cite[Assumption 3.1]{dhara-hofstad-leeuwaarden-sen-cm-finite-third-moment} on the degree sequence, the following results were obtained while proving
\cite[Theorem 3.3]{dhara-hofstad-leeuwaarden-sen-cm-finite-third-moment}:
(A) Fix $\eta>0$ and $K\geq 1$.
Then there exists $T>0$ such that with probability at least $1-\eta$, the $K$ maximal components are completely explored by the $Tn^{2/3}$-th step in the above depth-first exploration.
(B) For any $T>0$, the number of surplus edges found in the first $Tn^{2/3}$ steps of the exploration forms a tight sequence.

By Proposition \ref{prop:perc-cm}, conditional on its degree sequence, $\Perc_n(\fq_n^+)$ is distributed as a configuration model with that degree sequence.
Further, verifying that the degree sequence of $\Perc_n(\fq_n^+)$ satisfies 
\cite[Assumption 3.1]{dhara-hofstad-leeuwaarden-sen-cm-finite-third-moment} is routine; 
this is proved in \cite[Section 5]{SB-SD-vdH-SS-2020} for configuration models with a heavy-tailed degree sequence, and an identical argument works in the present setting.
Thus, for any $\eta>0$, we can choose $T_{\eta}>0$ and $\fs_{\eta}\in\bN$ such that 
with probability at least $1-\eta$, in the depth-first exploration of $\Perc_n(\fq_n^+)$,
$\cC_{\sss i, \Perc} (\fq_n^+)$ is completely explored by the $T_{\eta}n^{2/3}$-th step,
and at most $\fs_{\eta}$ many surplus edges are found in the first $T_{\eta}n^{2/3}$ steps.
Thus, it follows from Lemma \ref{lem:girth-auxiliary}~(b) that
\begin{align}\label{eqn:666}
	\pr\Big(
	\girth\big(\cC_{\sss i, \Perc} (\fq_n^+)\big)< n^{\frac{1}{3}-\eps}
	\Big)
	\leq
	\eta
	+
	\pr(\fA_n^c)
	+\sum_{j=1}^{T_{\eta}n^{2/3}}\pr(A_{n, j}\cap \fA_n)\, ,
\end{align}
where 
$
\fA_n:=
\{
\#E\big(\Perc_n(\fq_n^+)\big)\geq \mu n/(4\nu)
\}
$, and 
$A_{n,j}$ denotes the event that at most $\fs_{\eta}$ many surplus edges are found up to the $(j-1)$-th step of the exploration, and a surplus edge is found in the $j$-th step whose endpoints are within distance $n^{1/3-\eps}$ in the graph constructed up to the discovery of the vertex found in the $j$-th step.
Now, in the depth-first exploration, each surplus edge is of the form $\{v, v'\}$, where either $v'$ is a vertex on the ancestral line of $v$ or vice versa.
Applying Lemma \ref{lem:girth-auxiliary}~(a)  with $u_1$ as the vertex found in the $j$-th step of the exploration and $\Gamma_0$ as the ancestral line of this vertex and using Condition \ref{ass:cm-degree-det}~(a), we see that for all large $n$ and $1\leq j\leq T_{\eta}n^{2/3}$,
$
\pr(A_{n, j}\cap\fA_n)
\leq
n^{1/3-\eps} h(\fs_{\eta})(B\log n)^2\cdot
\Big(\mu n/(2\nu)  - T_{\eta} n^{2/3}\cdot (B\log n) \Big)^{-1}
$.
Combining this with \eqref{eqn:666} and the fact that 
$\pr(\fA_n^c)\to 0$ completes the proof.
\qed

We will next prove Theorem \ref{thm:cm-blob-level-surp-blob}~(b).
Using Lemma \ref{lem:cmn-t-k-2k-equiv}, conditional on $\CM_n(t_n)$, $\CM_n(t_c+\lambda_n/n^{1/3})$ can be constructed in the following steps:
Sample $k_n:=n\big(\bars_1(t_n)-\bars_1(t_c+\lambda_n/n^{1/3})\big)$.
Draw $k_n$ many half-edges without replacement from the set of $n\bars_1(t_n)=\sum_j f_j^{\circ}$ many half-edges alive at time $t_n$.
Let $a_i$ $(\leq f_i^{\circ})$ denote the number of half-edges sampled that are attached to the blob $\fb_i$.
Perform a uniform perfect matching of these $k_n$ half-edges to construct $\CM_n(t_c+\lambda_n/n^{1/3})$.
The next lemma collects some useful properties of the random variables $a_i$.

\begin{lem}\label{lem:a-i}
	We have,
	\begin{inparaenuma}
		\item
		$k_n=\sum_i a_i=\mu n^{1-\delta}+2\mu(\nu-1)\lambda_n n^{2/3}/\nu+O_P(n^{1-2\delta})$,
		\vskip3pt
		\noindent\item 
		$\max_i a_i=O_P(n^{\delta}\log^4 n)$, \ \ \ \
		\item
		$\sum_i a_i^2=2\mu n^{1-\delta}+6\mu(\nu-1)\lambda n^{2/3}/\nu+o_P(n^{2/3})$,
		\vskip4pt
		\noindent\item
		$\sum_i a_i^3\sim n\beta/\nu^3$, \ \ \ \ 
		\item
		$\sum_i a_i f_i^{\circ}\sim n\mu(\nu -1)/\nu$,\ \ \ \  and \ \ \ \ 
		\item
		$\sum_i a_i |\fb_i|\sim n\mu/\nu$.
	\end{inparaenuma}
\end{lem}

\noindent{\bf Proof:}
Let us first make a note of the following asymptotics:
\begin{align}\label{eqn:dd}
	\bars_1(t_n)=\mu\Big(\frac{\nu-1}{\nu}\Big)\big(1+O_P(n^{-\delta})\big)\, ,\  \text{ and }\ 
	\bars_2(t_n)=\mu\Big(\frac{\nu-1}{\nu}\Big)^2 n^{\delta}\big(1+o_P(n^{-\frac{1}{3}+\delta})\big)\, .
\end{align}
Here, the former relation uses Lemma \ref{lem:s1-conc}, and the latter follows from \eqref{eqn:s2-asymp-cm}. 
Now, the claim in (a) follows from Lemma \ref{lem:s1-conc}.
Next, write $a_i=\sum_{j=1}^{f_i^\circ } Y_{i, j}$, where
$Y_{i, j}$ is the indicator of the event that the $j$-th alive (at time $t_n$) half-edge attached to the blob $\fb_i$ is selected in the sampling process.
Let $\pr_1(\cdot)=\pr(\, \cdot\, |\,\CM_n(t_n), k_n)$, and similarly define $\bE_1[\cdot]$, $\var_1(\cdot)$, and $\cov_1(\cdot\, ,\cdot)$.
Note the following: 
\begin{inparaenumi}
	\item
	By the result in (a) and the first relation in \eqref{eqn:dd},
	$\bE_1[Y_{i, j}]=k_n/(n\bars_1(t_n))\sim \nu n^{-\delta}/(\nu-1)$,
	for $j=1,\ldots, f_i^\circ$.
	\item 
	By Theorem \ref{thm:negative-association}, $Y_{i, j}$, $j=1,\ldots, f_i^\circ$, are negatively associated under the measure $\pr_1$.
	\item 
	By Theorem \ref{thm:config-largest-comp-diam} and Condition \ref{ass:cm-degree-det}~(a), $\max_i f_i^\circ=O_P(n^{2\delta}\log^4 n)$.
\end{inparaenumi}
Combining these observations with the bound
$\pr_1(\max_i a_i\geq\ell)\leq \sum_i\pr_1(\sum_{j=1}^{f_i^\circ } Y_{i, j}\geq \ell)$, 
an application of Bennett's inequality 
\cite[Section 2.7]{boucheron2013concentration} for sums of negatively associated random variables \cite[Theorem 1]{qimanshao} yields the claim in (b).

To prove (c), note that
\begin{align*}
	\bE_1\big[\sum_i a_i(a_i-1)\big]
	=
	\sum_i \frac{f_i^\circ(f_i^\circ-1) k_n(k_n-1)}{n\bars_1(t_n)\big(n\bars_1(t_n)-1\big)}
	=
	\mu n^{1-\delta}+\frac{4\mu(\nu-1)}{\nu}\lambda n^{2/3}+o_P(n^{2/3})\, ,
\end{align*}
where the first step uses the expression for the second factorial moment of a hypergeometric random variable, and the second step follows from the result in (a), \eqref{eqn:dd}, the fact that $\delta\in(1/6, 1/5)$, and some routine calculations.
Thus, it is enough to show that $\var_1(\sum_i a_i^2)=O_P(n)$.
To this end, observe that
\begin{align*}
	\var_1(a_i^2)
	=
	\var_1\Big(
	\sum_{j=1}^{f_i^\circ} Y_{i, j}
	+\sum_{j_1\neq j_2}Y_{i, j_1}  Y_{i, j_2}
	\Big)
	\leq
	\sum_{j=1}^{f_i^\circ}\var_1\big(Y_{i, j}\big)
	+\sum_{j_1\neq j_2}\var_1\big(Y_{i, j_1}  Y_{i, j_2}\big)\\
	+2\sum_{j_1\neq j_2}\cov_1\big(Y_{i, j_1},\, Y_{i, j_1}  Y_{i, j_2}\big)
	+4\sum_{j_1, j_2, j_3\text{ distinct}}\cov_1\big(Y_{i, j_1}  Y_{i, j_2},\, Y_{i, j_1}  Y_{i, j_3}\big)\, ,
\end{align*}
since the other covariance terms are nonpositive by negative association.
Thus,
\begin{align*}
	&
	\var_1(\sum_i a_i^2)
	\leq
	\sum_i \var_1(a_i^2)
	\leq
	\sum_i\Big[\sum_{j=1}^{f_i^\circ}\bE_1(Y_{i, j})
	+
	\sum_{j_1\neq j_2}\bE_1\big(Y_{i, j_1}  Y_{i, j_2}\big)\\
	&\hskip120pt
	+2\sum_{j_1\neq j_2}\bE_1\big(Y_{i, j_1} Y_{i, j_2}\big)
	+4\sum_{j_1, j_2, j_3\text{ distinct}}\bE_1\big(Y_{i, j_1}  Y_{i, j_2} Y_{i, j_3}\big)\Big]\\
	&
	\leq 
	4\sum_i\Big[
	\sum_{j=1}^{f_i^\circ}\bE_1(Y_{i, j})
	+\sum_{j_1\neq j_2}\bE_1(Y_{i, j_1})\bE_1(Y_{i, j_2})
	+\sum_{j_1, j_2, j_3\text{ distinct}}\bE_1(Y_{i, j_1})\bE_1(Y_{i, j_2})\bE_1(Y_{i, j_3})
	\Big]\\
	&
	=O_P\Big(\sum_i\Big[
	f_i^\circ n^{-\delta}
	+\big(f_i^\circ\big)^2 n^{-2\delta}
	+\big(f_i^\circ\big)^3 n^{-3\delta}
	\Big]
	\Big)
	=
	O_P(n)\, ,
\end{align*}
where the first and the third steps use negative association, 
the fourth step uses the relation $\bE_1[Y_{i, j}]\sim \nu n^{-\delta}/(\nu-1)$,
and the last step follows from Theorem \ref{thm:config-bare-sibcrit}.
This completes the proof of (c).
The claims in (d), (e), and (f) can be proved in a similar way using Theorem \ref{thm:config-bare-sibcrit}.
We omit the details to avoid repetition.
\qed

\vskip3pt

For $t\geq t_n$, let $\CM_n^{\sss \oplus}(t)$ be the random graph obtained from  $\CM_n(t)$ by viewing each blob as a single vertex.
Then by Lemma \ref{lem:cmn-t-k-2k-equiv}, conditional on $\CM_n(t_n)$ and the sequence $a_1, a_2,\ldots$, $\CM_n^{\sss \oplus}(t_c+\lambda_n/n^{1/3})$ is distributed as a configuration model with degree sequence $a_1, a_2,\ldots$.
Consider the depth-first walk of a configuration model as defined around \cite[Display (5.2)]{dhara-hofstad-leeuwaarden-sen-cm-finite-third-moment} associated with the exploration described in \cite[Algorithm 1, page 9]{dhara-hofstad-leeuwaarden-sen-cm-finite-third-moment}.
Let $\big(Z_n^{\sss \oplus}(i);\, i=0,1,\ldots\big)$ be the depth-first walk of 
$\CM_n^{\sss \oplus}(t_c+\lambda_n/n^{1/3})$.
Recall the definition of $W_{\lambda}(t)$ from \eqref{eqn:parabolic-bm}.

\begin{prop}\label{prop:walk-convg-cm}
	Let $\alpha(\lambda)$ be as in \eqref{eqn:cm-sec-aldous-ver}.
	Then
	\begin{align}\label{eqn:25}
		\Big(n^{-1/3}  Z_n^{\sss \oplus}(\lfloor un^{2/3-\delta}\rfloor);\, u\geq 0\Big)
		\weakc
		\Big(
		\frac{\beta^{1/3}}{\nu}W_{\alpha(\lambda)}\Big(\frac{\beta^{1/3}u}{\mu\nu}\Big);\,  u\geq 0
		\Big)
	\end{align}
	with respect to the Skorohod $J_1$ topology.
\end{prop}

In view of the estimates in  Lemma \ref{lem:a-i}, Proposition \ref{prop:walk-convg-cm} can be proved using an argument identical to that used in the proof of 
\cite[Theorem 5.5]{dhara-hofstad-leeuwaarden-sen-cm-finite-third-moment}.
An outline is given in Appendix \ref{sec:appendix}.  
Finally, we need an elementary lemma to prove Theorem \ref{thm:cm-blob-level-surp-blob}~(b); we omit its proof.

\begin{lem}\label{lem:zzz}
	Suppose $X_1$ and $X_2$ are random variables defined on the same probability space and $X_1\leq X_2$ a.s.
	If, further, $X_1\equald X_2$, then $X_1=X_2$ a.s.
\end{lem}

\vskip3pt

\noindent \css{{\bf Proof of Theorem \ref{thm:cm-blob-level-surp-blob}~(b):}} 
For $t\geq t_n$ and a component $\cI^{\sss\oplus}$ in $\CM_n^{\sss \oplus}(t)$, we will write $\cI$ for the corresponding component in $\CM_n(t)$. 
Thus, $|\cI^{\sss\oplus}|$ equals the number of blobs in $\cI$.
Now, suppose we are working on a space on which the processes $\CM_n(\cdot)$, $n\geq 1$, and $W_{\alpha(\lambda)}(\cdot)$ are all defined, and the convergence in Proposition \ref{prop:walk-convg-cm} happens almost surely.
Then \cite[Lemma 7]{aldous1997brownian} implies that in this space, for all $n\geq 1$, there exist components $\cI_i^{\sss\oplus}$, $i\geq 1$, (some of which can be empty) of $\CM_n^{\sss\oplus}(t_c+\lambda_n/n^{1/3})$ such that
\begin{align}\label{eqn:4545}
	n^{\delta-2/3}\big( |\cI_i^{\sss\oplus}|\, ;\, i\geq 1\big)
	\convas
	\mu\nu\beta^{-1/3}\mvxi(\alpha(\lambda))
\end{align}
with respect to the product topology on $\bR^\bN$, and further, for any $K\geq 1$ and $\eps>0$, there exists $T_{K,\eps}>0$ such that 
\begin{align}\label{eqn:4646}
	\liminf_{n\to\infty}\,
	\pr\big(\cE_{n; K,\eps}\big)
	\geq 1-\eps\, ,
\end{align}
where $\cE_{n; K,\eps}$ denotes the event that $\cI_1^{\sss \oplus},\ldots, \cI_K^{\sss\oplus}$ are the $K$ maximal components among the components completely explored by the depth-first exploration of $\CM_n^{\sss\ \oplus}(t_c+\lambda_n/n^{1/3})$ up to the $T_{K, \eps}n^{2/3-\delta}$-th step.
Suppose we prove
\begin{align}\label{eqn:446}
	\Big(
	\frac{\beta^{1/3}}{\mu n^{2/3}}|\cI_i|\, ;\, i\geq 1
	\Big)
	\weakc
	\mvxi(\alpha(\lambda))\, ,\ \text{ and }\
	\Big(
	\frac{\beta^{1/3}\sW(\cI_i)}{\mu(\nu-1) n^{2/3}}\, ;\, i\geq 1
	\Big)
	\weakc
	\mvxi(\alpha(\lambda))\, ,
\end{align}
both with respect to the product topology on $\bR^\bN$.
Then comparing the the first convergence in \eqref{eqn:446} with Theorem \ref{thm:joseph}~(a), we see that $(n^{-2/3}|\cI_i|;\, i\geq 1)$ has the same distributional limit as $(n^{-2/3}|\cC_i|;\, i\geq 1)$.
Suppose there exists a subsequence $n_j$ and $\eta>0$ such that
$\pr(\cI_1\neq \cC_1)\geq\eta$ along the subsequence $n_j$.
Let $n_{j_r}$ be a further subsequence along which
$n^{-2/3}\big(|\cI_1|, |\cC_1|, |\cC_2|\big)$ converges in distribution to, say, $(Y_0, Y_1, Y_2)$.
Then $Y_0\equald Y_1$.
Further, $Y_0\leq Y_1$ a.s., since $|\cI_1|\leq |\cC_1|$.
Thus, Lemma \ref{lem:zzz} implies that $Y_0=Y_1$ a.s.
However, since $\eta \leq \pr(\cI_1\neq \cC_1) \leq \pr(|\cI_1|\leq |\cC_2|)$ along the subsequence $n_{j_r}$, 
we have $\eta\leq\pr(Y_0\leq Y_2)$.
Further, $Y_2<Y_1$ a.s., as $\gamma_2(\alpha(\lambda))<\gamma_1(\alpha(\lambda))$ a.s., and consequently, $\eta\leq\pr(Y_0< Y_1)$, which leads to a contradiction.
Hence, $\pr(\cI_1\neq \cC_1)\to 0$.
Iterating this argument and using the fact that for all $i\geq 1$ $\gamma_i(\alpha(\lambda))>\gamma_{i+1}(\alpha(\lambda))$ a.s., we conclude that
$\pr(\cI_i\neq \cC_i)\to 0$ for all $i\geq 1$.
Combining this with the second convergence in \eqref{eqn:446} completes the proof of 
Theorem \ref{thm:cm-blob-level-surp-blob}~(b).

Now we have to prove \eqref{eqn:446}.
As discussed at the beginning of 
\cite[Section 5.2]{dhara-hofstad-leeuwaarden-sen-cm-finite-third-moment},
in the depth-first exploration of $\CM_n^{\sss \oplus}(t_c+\lambda_n/n^{1/3})$, the vertices (which are blobs themselves) appear in a size-biased order, where the size of the blob $\fb_i$ is $a_i$.
To prove the first convergence in \eqref{eqn:446}, fix $\eps>0$ and $K\geq 1$, 
and apply Lemma \ref{lem:size-biased-partial-sum} with
$\vx=(a_1, a_2,\ldots)$, $\vu=(|\fb_1|, |\fb_2|,\ldots)$, and $\ell=T_{K,\eps}n^{2/3-\delta}$.
The conditions in \eqref{eqn:ass-size-biased-partial-sum} can be verified using
Lemma \ref{lem:a-i}~(a), (b), (f), and Theorem \ref{thm:config-largest-comp-diam}.
Then Lemma \ref{lem:size-biased-partial-sum} and \eqref{eqn:4646} imply that 
for $i=1,\ldots, K$,
\[
\sum_{j:\fb_j\in\cI_i^{\sss \oplus}}
\frac{|\fb_j| }{n^{2/3-\delta}(\sum_r a_r |\fb_r|/\sum_r a_r)}
-
\frac{|\cI_i^{\sss \oplus}|}{n^{2/3-\delta}}
\weakc 
0\, .
\]
Now the first convergence in \eqref{eqn:446} follows upon using \eqref{eqn:4545} together with 
Lemma \ref{lem:a-i}~(a), and (f), and observing that
$|\cI_i|
=
\sum_{j:\fb_j\in\cI_i^{\sss \oplus}} |\fb_j|
$.
The second convergence in \eqref{eqn:446} follows in a similar way.
Here, we apply Lemma \ref{lem:size-biased-partial-sum} with
$\vx=(a_1, a_2,\ldots)$, $\vu=(f_1^\circ, f_2^\circ,\ldots)$, and $\ell=T_{K,\eps}n^{2/3-\delta}$, and use the results in Lemma \ref{lem:a-i}~(a), (b), (e) together with \eqref{eqn:4545}.
We omit the details as no new idea is involved.
\qed

\subsubsection{Coupling $\cG_n^{\modi}$ and $\CM_n$}\label{sec:cm-coupling-mod-main}
Throughout this section we will work with a given $\lambda\in\bR$ and a sequence $\lambda_n\to\lambda$.
Now, starting at time $t_n$ with $\CM_n(t_n)$ we have considered two processes: the original process $\CM_n$ and the modified process $\cG_n^{\modi}$. 
We have already derived the scaling limit of $\cG_n^{\modi}$. 
We would now like to transfer this result from $\cG_n^{\modi}$ to $\CM_n$ using a suitable coupling.

For $k\geq 0$, let $\CM_n^{**}(k):=\CM_n(\tau_k)$, where $\tau_0:=0$, and for $k\geq 1$, $\tau_k\geq t_n$ denotes the time when the $k$-th full edge is added in $\CM_n(\cdot)$ {\bf after} time $t_n$. 
(Note the difference between $\CM_n^{**}$ and $\CM_n^*$ introduced in Section \ref{sec:proof-cm-max-diam}.)
Write $\big(\cG_n^{*,\modi}(k);\, k\geq 0\big)$ for the corresponding discrete time chain for $\cG_n^{\modi}$ also initialized at time zero with $\CM_n(t_n)$.

Observe the following properties of $\CM_n(\cdot)$:
{\bf(a)} 
The process $\big(\CM_n^{**}(k);\, k\geq 0\big)$ is obtained by sequentially selecting pairs of alive half-edges uniformly at random {\bf without replacement} and forming full edges.
{\bf (b)} 
For $t> t_n$, let $\fR_n[t]$ be the number of full edges formed in $\CM_n(\cdot)$ in the interval $[t_n,t]$.
Then by Lemma \ref{lem:a-i}~(a),
\begin{equation}
	\label{eqn:rn-bound}
	\fR_n[t_c+\lambda_n/n^{1/3}] 
	=
	\mu n^{1-\delta}/2+\mu(\nu-1)\lambda_n n^{2/3}/\nu+O_P(n^{1-2\delta})\, .
\end{equation}
Similarly, the process $~\cG_n^{\modi}$ satisfies the following:
{\bf (a$^\prime$)} 
The process $\big(\cG_n^{*,\modi}(k): k\geq 0\big)$ is obtained by sequentially selecting half-edges at random {\bf with replacement} and forming full edges.
{\bf (b$^\prime$)} Using \eqref{eqn:rate-modi-def}, the rate at which full edges are created is constant and is given by
\begin{equation}
	\label{eqn:alphan-full-rate-modi}
	\alpha_n
	:= 
	(n^2 \bars_1^2(t_n))/(n s_1(t_c))
	= 
	n\mu\Big[(\nu-1)/\nu+2n^{-\delta}+O_P(n^{-2\delta})\Big]\, ,
\end{equation}
where the last step follows from Lemma \ref{lem:s1-conc}.
Call an edge formed in the process $\big(\cG_n^{*,\modi}(k);\, k\geq 0\big)$ `good,' if none of its half-edges previously contributed to an edge added in the process;
otherwise, call the edge `bad.' 
Using {\bf (a)} and {\bf (a$^\prime$)},
there is an obvious coupling between $\CM_n^{**}$ and $\cG_n^{*,\modi}$ in which
\begin{equation}\label{eqn:cmn-modi-dom}
	\CM_n^{**}(k) \subseteq \cG_n^{*,\modi}(r_k)\ \ \ \mbox{ for all } k\geq 0\, ,
\end{equation}
where $r_k$ denotes the number of edges added in $\cG_n^{*,\modi}$ when we get the $k$-th good edge.
In this coupling, every good edge added in $\cG_n^{*,\modi}$ also registers in $\CM_n^{**}$, and the bad edges formed in $\cG_n^{*,\modi}$ do not register in $\CM_n^{**}$.
From now on, we will assume that $\CM_n$ and $\cG_n^{\modi}$ are coupled so that \eqref{eqn:cmn-modi-dom} holds.

For any $t\geq t_n$, conditional on $\CM_n(t_n)$, the number of full edges formed in $\cG_n^{\modi}$ in the interval $[t_n, t]$ is given by
$
\fR_n^{\modi}[t] \sim \mbox{ Poisson } (\alpha_n (t-t_n))
$.
Let
\begin{align}\label{eqn:def-T-+-n}
	T^+_n:=t_c+\lambda_n n^{-1/3}+n^{-2\delta}\log n \, .
\end{align}
Using \eqref{eqn:alphan-full-rate-modi}, a direct calculation shows that
\[
\alpha_n (T^+_n-t_n)
=
\frac{\mu n^{1-\delta}}{2}
+\frac{\mu(\nu-1)}{\nu}\lambda_n n^{2/3}
+\frac{\mu(\nu-1)}{\nu}n^{1-2\delta}\log n 
+O_P(n^{1-2\delta})\, .
\]
It follows from standard tail bounds for Poisson random variables that
\begin{align}\label{eqn:3435}
	\fR_n^{\modi} [T^+_n]
	=
	\frac{\mu n^{1-\delta}}{2}
	+\frac{\mu(\nu-1)}{\nu}\lambda_n n^{2/3}
	+\frac{\mu(\nu-1)}{\nu}n^{1-2\delta}\log n 
	+O_P(n^{1-2\delta})\, .
\end{align}
Note that conditional on $\CM_n(t_n)$, the probability that the $k+1$-th edge added in $\cG_n^{*,\modi}$ is bad is bounded by $4k/(n\bars_1(t_n))$.
Let $\fB_n$ denote the number of bad edges created in the interval $[t_n, T_n^+]$.
Then 
\begin{align}\label{eqn:343}
	\bE\big(\fB_n\, \big|\, \CM_n(t_n), \fR_n^{\modi} [T^+_n]\big)
	\leq
	\sum_{k=1}^{\fR_n^{\modi} [T^+_n]}4k/(n\bars_1(t_n))
	=
	O_P(n^{1-2\delta})\, ,
\end{align}
where the last step uses \eqref{eqn:3435}.
Combining \eqref{eqn:343} with \eqref{eqn:3435}, it follows that the number of good edges formed in $[t_n, T_n^+]$ exhibits the same asymptotics as that on the right side of \eqref{eqn:3435}.
Thus, \eqref{eqn:rn-bound} and \eqref{eqn:cmn-modi-dom} imply the following:

\begin{prop}\label{prop:12}
We have,
$
\CM_n(t_c+\lambda_n/n^{1/3}) \subseteq \cG_n^{\modi}(T_n^+)
$ 
whp.
\end{prop}

For the rest of this section we will write $\cC_i$ for $\cC_i(t_c+\lambda_n/n^{1/3})$.
Comparing Proposition \ref{prop:modi-free-max} with Theorem \ref{thm:cm-blob-level-surp-blob}~(b), we see that $\big(n^{-2/3}\sW(\cC_i); \geq 1\big)$ has the same distributional limit $\big(n^{-2/3}\sW(\sC_i^{\modi}(T_n^+)); \geq 1\big)$.
This fact coupled with Proposition \ref{prop:12} shows that 
\begin{align}\label{eqn:2929}
	\text{for each }\ \ i\geq 1,\ \ \
	\cC_i\subseteq \sC_i^{\modi}(T_n^+)\ \ \text{ whp.}
\end{align}
Now, Proposition \ref{prop:modi-max-comp-sizes} and Theorem \ref{thm:joseph}~(a) show that $\big(n^{-2/3}|\cC_i|; \geq 1\big)$ has the distributional limit as $\big(n^{-2/3}|\sC_i^{\modi}(T_n^+)|; \geq 1\big)$, which together with \eqref{eqn:2929} imply that for each $i\geq 1$,
$|\sC_i^{\modi}(T_n^+)|-|\cC_i|=o_P(n^{2/3})$.
If $\cC_i\subseteq \sC_i^{\modi}(T_n^+)$, and every surplus edge in $\cC_i$ is created in the time interval $[t_n, t_c+\lambda_n/n^{1/3}]$ (which occurs whp by Theorem \ref{thm:cm-blob-level-surp-blob}~(a)), then 
$N_i(t_c+\lambda_n/n^{1/3})\leq\fN_i^{\modi}(T_n^+)$.
However, by Theorem \ref{thm:joseph}~(a) and Proposition \ref{prop:modi-max-comp-sizes}, $N_i(t_c+\lambda_n/n^{1/3})$ and $\fN_i^{\modi}(T_n^+)$ have the same distributional limits, and consequently, for each $i\geq 1$, $N_i(t_c+\lambda_n/n^{1/3})=\fN_i^{\modi}(T_n^+)$ whp.
We collect these observations in the following corollary:

\begin{cor}\label{cor:surplus-cm}
	For each $i\geq 1$, the following hold:
	\begin{inparaenuma}
		\item
		Whp
		$\cC_i\subseteq \sC_i^{\modi}(T_n^+)$ , and
		$
		|\sC_i^{\modi}(T_n^+)|-|\cC_i|=o_P(n^{2/3})
		$.
		\item 
		Whp $N_i(t_c+\lambda_n/n^{1/3})=\fN_i^{\modi}(T_n^+)$, and 
		all surplus edges in $\cC_i$ are created in $[t_n, t_c+\lambda_n/n^{1/3}]$.
	\end{inparaenuma}
\end{cor}

Note that on the intersection of the events described in Corollary \ref{cor:surplus-cm}, ${\sC}_i^{\modi}(T_n^+)$ can be obtained from $\cC_i$ by attaching a finite number of connected graphs (that are tree superstructures whose vertices are blobs) each to a vertex of $\cC_i$ via a single edge.

\subsubsection{Completing the proof of Theorem \ref{thm:crit-main-res-cm}}
As before, we will write $\cC_i$ for $\cC_i(t_c+\lambda_n/n^{1/3})$.
For $i\geq 1$, define
\[
\bar{\cC}_i
:= 
\scl(n^{-1/3}, n^{-2/3}) \cC_i\, ,\ \ \text{ and }\ \ 
\bar{\sC}_i^{\modi}(T_n^+):=
\scl(n^{-1/3}, n^{-2/3}) \sC_i^{\modi}(T_n^+)\, .
\]
Then the claim in Theorem \ref{thm:crit-main-res-cm} follows upon combining the
next proposition with Theorem \ref{thm:modi-scalin-lim-count}.
\begin{prop}\label{prop:dhp-modi-orig-zero}
	For each $i\geq 1$, 
	$d_{\GHP}\big(
	\bar{\cC}_i \;,\; \bar{\sC}_i^{\modi}(T_n^+)
	\big) 
	\weakc 0$,
	as $n\to\infty$.	
\end{prop}

\noindent {\bf Proof:} 
Fix $i\geq 1$.
In view of Theorem \ref{thm:modi-scalin-lim-count} and Corollary \ref{cor:surplus-cm}, we can assume that we are working on a space $(\Omega, \cF, \pr)$ where $\CM_n(\cdot)$ and $\cG_n^{\modi}(\cdot)$ are defined for all $n$, and the following hold:
\begin{enumeratei}
	\item
	As $n\to\infty$, 
	$d_{\GHP}\big(\bar{\sC}_i^{\modi}(T_n^+),\, Z\big)\convas 0$ 
	for some random compact metric measure space $Z$ defined on $\Omega$.
	\item 
	Almost surely $\cC_i\subseteq \sC_i^{\modi}(T_n^+)$ eventually, and
	\begin{align}\label{eqn:36}
		n^{-2/3}\big(|\sC_i^{\modi}(T_n^+)|-|\cC_i|\big)\convas 0\, .
	\end{align}
	\item 
	There exist a $\bN$-valued random variable $\zeta$ defined on $\Omega$ such that the following holds for $\pr$-almost every $\omega\in\Omega$, and for all $n\geq\zeta(\omega)$:
	there exist distinct vertices $v_1,\ldots, v_{r_n}$ in $\cC_i$ and finite, connected, and rooted graphs $\sS_1^{\sss(n)}, \sS_2^{\sss(n)}, \ldots, \sS_{r_n}^{\sss(n)}$ such that $\sC_i^{\modi}(T_n^+)$ is the space obtained by identifying the root of $\sS_i^{\sss(n)}$ with $v_i$, $i=1,\ldots, r_n$.
	
\end{enumeratei}

From the construction of the spaces $\crit_i(\cdot)$, it follows that for any $\eta>0$, the infimum of the  measures of all balls of radius $\eta$ in $Z$ is positive almost surely.
Hence, for any $\eta>0$,
\begin{align}\label{eqn:27}
	\liminf_n \Big(\min_{v\in V({\sC}_i^{\modi}(T_n^+))} n^{-2/3}\big|\cB(v, \eta n^{1/3})\big|\Big)>0\ \ \ \
	\text{ a.s.}\, ,
\end{align}
where $\cB(v, a)$ denotes the ball of radius $a$ around $v$.
Note also that for $n\geq \zeta$,
\begin{align*}
	\frac{1}{n^{2/3}}\max_{1\leq j\leq r_n}(|\sS_j^{\sss(n)}|-1)
	\leq
	\frac{1}{n^{2/3}}\sum_{j=1}^{r_n}(|\sS_j^{\sss(n)}|-1)
	=
	\frac{1}{n^{2/3}}\big(|\sC_i^{\modi}(T_n^+)|-|\cC_i|\big)\, ,
\end{align*}
which combined with \eqref{eqn:36} and \eqref{eqn:27} implies that
\begin{align}\label{eqn:20}
	n^{-1/3}\max_{1\leq j\leq r_n}\diam(\sS_j^{\sss(n)} )
	\convas 0
	\text{ a.s.}
\end{align}
For $n\geq\zeta$, endow $\cC_i$ with the measure that assigns mass $|\sS_j^{\sss(n)}|$ to the vertex $v_j$, $j=1,\ldots, r_n$, and mass $1$ to every vertex $v\notin\{v_1,\ldots, v_{r_n}\}$, and call the resulting metric measure space $\cC_{i; \project}$.
For $n<\zeta$, set $\cC_{i; \project}=\cC_i$.
Let 
$
\bar{\cC}_{i; \project}
:= 
\scl(n^{-1/3}, n^{-2/3}) \cC_{i; \project}
$.
Using \eqref{eqn:36}, it follows that the total variation distance between the measures on 
$\bar{\cC}_{i; \project}$ and $\bar{\cC}_i$ goes to zero almost surely, and consequently,
$d_{\GHP}\big(\bar{\cC}_{i; \project},\, \bar{\cC}_i\big)\convas 0$.
Further, \eqref{eqn:20} implies that 
$d_{\GHP}\big(\bar{\cC}_{i; \project},\, \bar{\sC}_i^{\modi}(T_n^+) \big)\convas 0$.
We complete the proof by combining these last two observations.
\qed

\subsection{Proof of  Theorem \ref{thm:perc-cm}}\label{sec:proof-perc-cm}
Fix $\eta\in(1/3, 1/2)$, and define
\[
T_n^{\bullet}:= t_c+ \frac{\nu}{2(\nu-1)}\frac{\lambda_n}{n^{1/3}}+\frac{1}{n^{\eta}}.
\]
We will assume that $\Perc_n(\cdot)$ and $\CM_n(\cdot)$ are coupled as in \eqref{eqn:couple-cm-perc} so that
\begin{align}\label{eqn:45}
	\Perc_n(p(\lambda_n))
	\subseteq 
	\CM_n(T_n^\bullet)\ \ \text{ whp.}
\end{align}
Now, Theorem \ref{thm:crit-main-res-cm} implies that the claimed convergence in \eqref{eqn:1} holds if we replace $\cC_{\sss i,\Perc}(\lambda_n)$ by $\cC_i(T_n^{\bullet})$.
Thus, it is enough to show that for each $i\geq 1$,
\begin{align}\label{eqn:76}
	d_{\GHP}\Big(
	\scl(n^{-1/3}, n^{-2/3}) \cC_{\sss i,\Perc}(\lambda_n)
	,\,
	\scl(n^{-1/3}, n^{-2/3}) \cC_i(T_n^{\bullet})
	\Big)
	\weakc 0\, .
\end{align}

By Theorem~\ref{thm:joseph}~(a) and Theorem~\ref{thm:cm-perc-component-size-scaling}, 
$\big(n^{-2/3}|\cC_{\sss i,\Perc}(\lambda_n)|;\, i\geq 1\big)$ has the same distributional limit as $\big(n^{-2/3}|\cC_i(T_n^{\bullet})|;\, i\geq 1\big)$.
Thus, \eqref{eqn:45} implies that for each $i\geq 1$, 
$\cC_{\sss i,\Perc}(\lambda_n)\subseteq\cC_i(T_n^{\bullet})$ whp, and 
$|\cC_i(T_n^{\bullet})|-|\cC_{\sss i,\Perc}(\lambda_n)|=o_P(n^{2/3})$.
If $\cC_{\sss i,\Perc}(\lambda_n)\subseteq\cC_i(T_n^{\bullet})$, then
$\spls\big(\cC_{\sss i,\Perc}(\lambda_n)\big)\leq N_i(T_n^{\bullet})$.
Again, by Theorem \ref{thm:cm-perc-component-size-scaling} and Theorem \ref{thm:joseph}~(a), 
$\big(\spls\big(\cC_{\sss i,\Perc}(\lambda_n)\big);\, i\geq 1\big)$ and
$\big(N_i(T_n^{\bullet});\, i\geq 1\big)$ have the same distributional limit, which in turn implies that for each $i\geq 1$, 
$\spls\big(\cC_{\sss i,\Perc}(\lambda_n)\big) = N_i(T_n^{\bullet})$ whp.
Thus, \eqref{eqn:76} follows by an argument similar to the one used in the proof of Proposition \ref{prop:dhp-modi-orig-zero}.
\qed

\appendix

\section{}\label{sec:appendix}
We will now outline the steps in the proof of Proposition \ref{prop:walk-convg-cm}.
We once again ask the reader to recall the depth-first exploration of a configuration model given in \cite[Algorithm 1, page 9]{dhara-hofstad-leeuwaarden-sen-cm-finite-third-moment}, which we are not describing here.
For $j\geq 1$, let $\fb_{v(j)}$ be the blob found in the $j$-th step of the depth-first exploration of $\CM_n^{\sss \oplus}(t_c+\lambda_n/n^{1/3})$.
As mentioned in \cite[Display 5.2]{dhara-hofstad-leeuwaarden-sen-cm-finite-third-moment},
\begin{align}\label{eqn:23}
	Z_n^{\sss\oplus}(0)=0\, ,\ \ \text{ and }\ \
	Z_n^{\sss\oplus}(i)
	=
	\sum_{j=1}^i\big(a_{(j)}-2-2\fs_{(j)}\big)\ \ \text{ for }\ \ i\geq 1\, ,
\end{align}
where $a_{(j)}=a_{v(j)}$, and $\fs_{(j)}$ denotes the number of surplus edges found in the $j$-th step of the exploration.
\begin{lem}\label{lem:23}
	Let 
	$\tilde Z_n^{\sss\oplus}(0)=0$ and 
	$\tilde Z_n^{\sss\oplus}(i)
	=
	\sum_{j=1}^i\big(a_{(j)}-2\big)$.
	Then
	\begin{align}\label{eqn:24}
		\Big(n^{-1/3}  \tilde Z_n^{\sss \oplus}(\lfloor un^{2/3-\delta}\rfloor);\, u\geq 0\Big)
		\weakc
		\bigg(\Big(
		\sqrt{\frac{\beta}{\mu\nu^3}}\cdot B(u)
		+
		2\Big(\frac{\nu-1}{\nu}\Big)\lambda u
		-\frac{\beta u^2}{2\mu^2\nu^3}
		\Big);\,  u\geq 0
		\bigg)
	\end{align}
	with respect to the Skorohod $J_1$ topology, where $B(\cdot)$ is a standard Brownian motion.
\end{lem}

We will prove Lemma \ref{lem:23} shortly, but first note that by Brownian scaling, the right side of \eqref{eqn:24} has the same distribution as the process on the right side of \eqref{eqn:25}.
Thus, Proposition \ref{prop:walk-convg-cm} follows upon combining Lemma \ref{lem:23} and \eqref{eqn:23} with the next result.

\begin{lem}\label{lem:24}
	For any $T>0$,
	$\big(\sum_{j=1}^{Tn^{2/3-\delta}}\fs_{(j)}\big)$ is a tight sequence.
\end{lem}

\noindent{\bf Proof:}
Note that for any $i\geq 1$ and $j\leq i$,
$\tilde Z_n^{\sss\oplus}(j)\leq Z_n^{\sss\oplus}(j)+2\sum_{r=1}^i \fs_{(r)}$,
which in turn implies that
\begin{align}\label{eqn:34}
	Z_n^{\sss\oplus}(i)-\min_{ j\leq i}Z_n^{\sss\oplus}(j)
	=
	\tilde Z_n^{\sss\oplus}(i)- 2\sum_{j=1}^i \fs_{(j)}-\min_{ j\leq i}Z_n^{\sss\oplus}(j)
	\leq
	\tilde Z_n^{\sss\oplus}(i)-\min_{ j\leq i}\tilde Z_n^{\sss\oplus}(j)\, .
\end{align}
For $i\geq 0$, let $\cF_i^\star$ denote the $\sigma$-field generated by $\CM_n(t_n)$, $(a_1, a_2,\ldots)$, the information revealed up to the $i$-th step of the depth-first exploration, and $a_{(i+1)}$.
Now, in the $(i+1)$-th step of the depth-first exploration, each of the $a_{(i+1)}$ half-edges of the blob $\fb_{v(i+1)}$ can create a surplus edge with 
$\big(Z_n^{\sss\oplus}(i)-\min_{ j\leq i}Z_n^{\sss\oplus}(j) +a_{(i+1)}-1\big)$ 
many other half-edges.
Consequently, using Condition \ref{ass:cm-degree-det}~(a) and \eqref{eqn:34}, 
on the event $\cA_n:=\{\sum_j a_j\geq \mu n^{1-\delta}/2\}$,
for any $i\leq Tn^{2/3-\delta}$,
\begin{align}\label{eqn:26}
	\bE\big[
	\fs_{(i+1)}\, |\,\cF_i^\star
	\big]
	\leq
	\big(\tilde Z_n^{\sss\oplus}(i)-\min_{ j\leq i}\tilde Z_n^{\sss\oplus}(j) +B\log n\big)
	\cdot
	\frac{a_{(i+1)}}{\sum_j a_j- (i+1)B\log n}\, .
\end{align}
Using Lemma \ref{lem:23}, for any $\eps, T>0$, we can choose $K>0$ such that 
$\pr\big(
\cE_n^c\big)
\leq
\eps$,
where
$
\cE_n:=
\big\{
\max_{0\leq i\leq Tn^{2/3-\delta}}|\tilde Z_n^{\sss\oplus}(i)|
\leq 
Kn^{1/3}
\big\}
$.
On the event $\cE_n$, 
\begin{align}\label{eqn:28}
	\sum_{j=1}^{Tn^{2/3-\delta}} a_{(j)}
	=
	\tilde Z_n^{\sss\oplus}\big(Tn^{2/3-\delta}\big)+2Tn^{2/3-\delta}
	\leq
	(K+2T)n^{2/3-\delta}\, .
\end{align}
Combining \eqref{eqn:26} and \eqref{eqn:28} gives
\begin{align*}
	\bE\Big[
	\sum_{j=1}^{Tn^{2/3-\delta}} \fs_{(j)}\ind_{\cE_n\cap\cA_n}
	\Big]
	\leq
	\frac{(2K+B)n^{1/3}\cdot (K+2T)n^{2/3-\delta}}{\mu n^{1-\delta}/2-BTn^{2/3-\delta}\log n}
	\leq
	\frac{4}{\mu}(2K+B)(K+2T)
\end{align*}
for all large $n$.
Since $\pr(\cA_n^c)\to 0$ by Lemma~\ref{lem:a-i}~(a), this completes the proof.
\qed

\vskip3pt

We now turn to the proof of Lemma \ref{lem:23}, which follows from an application of the martingale functional central limit theorem.
For $i\geq 0$, let $\cF_i$ denote the $\sigma$-field generated by $\CM_n(t_n)$, $(a_1, a_2,\ldots)$, and the information revealed up to the $i$-th step of the depth-first exploration.
Let 
$\cA_n(i):=
\sum_{j=1}^i\bE\big[a_{(j)}-2\, |\, \cF_{j-1}\big]
$, and
$\cB_n(i):=
\sum_{j=1}^i\big(\bE\big[a_{(j)}^2\, |\, \cF_{j-1}\big]
-
\bE^2\big[a_{(j)}\, |\, \cF_{j-1}\big]\big)
$.
Then as observed around
\cite[Display (5.26a)]{dhara-hofstad-leeuwaarden-sen-cm-finite-third-moment}, it is enough to show that
\begin{gather}
	\sup_{u\leq t}\Big|
	\frac{1}{n^{1/3}}\cA_n\big(\lfloor un^{2/3-\delta}\rfloor\big)
	+
	\frac{\beta u^2}{2\mu^2\nu^3}
	-
	\frac{2(\nu-1)\lambda u}{\nu}
	\Big|
	\weakc
	0\, ,\ \ \text{ and} \label{eqn:56}\\
	n^{-2/3}\cB_n\big(\lfloor tn^{2/3-\delta}\rfloor\big)
	\weakc
	\beta t /(\mu \nu^3)\, ,\label{eqn:57}
\end{gather}
for any $t>0$.
(The analogues of 
\cite[Displays (5.26c) and (5.26d)]{dhara-hofstad-leeuwaarden-sen-cm-finite-third-moment}
follow trivially in our setting from Lemma \ref{lem:a-i}~(b).)
We prove \eqref{eqn:56} and \eqref{eqn:57} below to complete the proof of Lemma \ref{lem:23}.

\vskip3pt

\noindent{\bf Proof of \eqref{eqn:56}:}
As explained in \cite[Section 5.2]{dhara-hofstad-leeuwaarden-sen-cm-finite-third-moment}, conditional on $\big(\CM_n(t_n), (a_1, a_2,\ldots)\big)$, $(\fb_{v(1)}, \fb_{v(2)},\ldots)$ is a size-biased permutation of the blobs with weight sequence $(a_1, a_2,\ldots)$.
Thus,
\begin{align}\label{eqn:10}
	\bE\big[
	a_{(j)}-2\, |\, \cF_{j-1}
	\big]
	=
	\frac{\sum_k a_k (a_k-2)-\sum_{r=1}^{j-1}a_{(r)}(a_{(r)}-2)}{\sum_k a_k- \sum_{r=1}^{j-1}a_{(r)}}\, .
\end{align}
Fix $t>0$, and apply Lemma \ref{lem:size-biased-partial-sum} with 
$\ell=tn^{2/3-\delta},$ $x_i=a_i$, and $u_i=a_i(a_i-2)$.
Then Lemma \ref{lem:a-i}~(a), (b), (c), and (d) imply that the conditions in \eqref{eqn:ass-size-biased-partial-sum} are satisfied, and further,
\begin{align}\label{eqn:11}
	c_n=\sum_i x_i u_i/\sum_i x_i\sim \beta n^{\delta}/(\mu\nu^3)\, .
\end{align}
Now, Lemma \ref{lem:size-biased-partial-sum} implies that
\begin{align}\label{eqn:12}
	\max_{j\leq tn^{2/3-\delta}}
	\Big|
	\sum_{r=1}^{j-1}a_{(r)}(a_{(r)}-2)-(j-1)c_n
	\Big|
	=
	o_P(n^{2/3})\, .
\end{align}
Finally, by Lemma \ref{lem:a-i}~(a), (b), and (c),
$\sum_k a_k(a_k-2)=2\mu(\nu-1)\lambda n^{2/3}/\nu+o_P(n^{2/3})$,
$\sum_k a_k\sim \mu n^{1-\delta}$,
and 
$\sum_{r=1}^{j-1}a_{(r)}= jO_P(n^{\delta}\log^4 n)$.
Combining these observations with \eqref{eqn:10}, \eqref{eqn:11}, and \eqref{eqn:12}, a direct calculation yields \eqref{eqn:56}.
\qed

\vskip3pt

\noindent{\bf Proof of \eqref{eqn:57}:}
Here,
\begin{align}\label{eqn:10-a}
	\bE\big[
	a_{(j)}^2\, |\, \cF_{j-1}
	\big]
	=
	\frac{\sum_k a_k^3-\sum_{r=1}^{j-1}a_{(r)}^3}{\sum_k a_k- \sum_{r=1}^{j-1}a_{(r)}}\, .
\end{align}
Fix $t>0$, and apply Lemma \ref{lem:size-biased-partial-sum} with 
$\ell=tn^{2/3-\delta},$ $x_i=a_i$, and $u_i=a_i^3$.
Then Lemma \ref{lem:a-i}~(a),  (b), and (d) imply that 
$c_n=\sum_i x_i u_i/\sum_i x_i$ satisfies
\begin{align}\label{eqn:11-a}
	c_n\geq \frac{\sum_k a_k^3}{\sum_k a_k}=\Omega_P(n^{\delta})\, , \ \ \text{ and }\ \ 
	c_n\leq \frac{\max_j a_j\cdot\sum_k a_k^3}{\sum_k a_k}=O_P(n^{2\delta}\log^4 n)\, .
\end{align}
The conditions in \eqref{eqn:ass-size-biased-partial-sum} can be verified using the first relation in \eqref{eqn:11-a} and Lemma \ref{lem:a-i}.
Now, the second relation in \eqref{eqn:11-a} and Lemma \ref{lem:size-biased-partial-sum} imply that
$
\sum_{r=1}^{\lfloor tn^{2/3-\delta}\rfloor}a_{(r)}^3
=
O_P(n^{2/3+\delta}\log^4 n)
$.
Combining this observation with \eqref{eqn:10-a}, an application of Lemma \ref{lem:a-i} shows that
\[
n^{-2/3}
\sum_{j=1}^{\lfloor tn^{2/3-\delta}\rfloor}
\bE\big[
a_{(j)}^2\, |\, \cF_{j-1}
\big]
\weakc
\beta t /(\mu \nu^3)\, .
\]
Finally, using Lemma \ref{lem:a-i}, it is straightforward to show that
$
\sum_{j=1}^{\lfloor tn^{2/3-\delta}\rfloor}
\bE^2\big[a_{(j)}\, |\, \cF_{j-1}\big]
=
O_P(n^{2/3-\delta})
$.
This completes the proof of \eqref{eqn:57}.
\qed

\section*{Acknowledgments}
We thank Amarjit Budhiraja and Steve Evans for many stimulating conversations. 
We are also grateful to an anonymous referee for a detailed reading of a long paper which led to significant improvements in the exposition and proofs.
SB has been partially supported by NSF-DMS grants 1105581, 1310002, 1613072, 1606839, SES grant 1357622 and {ARO grant W911NF-17-1-0010}.
SS has been supported in part by NSF grant DMS-1007524, Netherlands Organisation for Scientific Research (NWO) through the Gravitation Networks grant 024.002.003, MATRICS grant MTR/2019/000745 from SERB, and by the  Infosys foundation, Bangalore.
XW has been supported in part by the National Science Foundation (DMS-1004418, DMS-1016441), the Army Research Office (W911NF-0-1-0080, W911NF-10-1-0158) and the US-Israel Binational Science Foundation (2008466).

\bibliographystyle{plain}
\bibliography{scaling}

\end{document}